\documentclass[oneside, 12pt, a4paper]{amsart}

\usepackage{amscd, amssymb, amsmath, mathrsfs}
\usepackage[main=english]{babel}
\usepackage{tikz-cd}
\usepackage[pdftex, colorlinks=true,  citecolor=blue, linkcolor=blue, linktocpage=true, backref=page]{hyperref}
\usepackage[all,cmtip]{xy}

\newcommand\mylabel[1]{\label{#1}}

\newtheorem{theorem}{Theorem}[section]
\newtheorem*{maintheorem}{Theorem}
\newtheorem{lemma}[theorem]{Lemma}
\newtheorem{proposition}[theorem]{Proposition}
\newtheorem{corollary}[theorem]{Corollary}

\theoremstyle{definition}
\newtheorem{definition}[theorem]{Definition}

\newtheorem*{acknowledgement}{Acknowledgement}

\theoremstyle{remark}

\DeclareFontFamily{U}{wncy}{}
\DeclareFontShape{U}{wncy}{m}{n}{<->wncyr10}{}
\DeclareSymbolFont{mcy}{U}{wncy}{m}{n}
\DeclareMathSymbol{\Sh}{\mathord}{mcy}{"58}

\newcommand{\NN}	{\mathbb{N}}
\newcommand{\ZZ}	{\mathbb{Z}}
\newcommand{\QQ}	{\mathbb{Q}}

\newcommand{\PP}	{\mathbb{P}}
\renewcommand{\AA}	{\mathbb{A}}
\newcommand{\GG}	{\mathbb{G}}

\newcommand{\ideal}[1]{\mathfrak #1}

\newcommand{\formal}[1]{\mathfrak #1}
\newcommand{\ul}[1]{\underline{#1}}

\newcommand  {\shA}     {\mathscr{A}}

\newcommand  {\shB}     {\mathscr{B}}
\newcommand  {\shC}     {\mathscr{C}}

\newcommand  {\shE}     {\mathscr{E}}
\newcommand  {\shF}     {\mathscr{F}}

\newcommand  {\shI}     {\mathscr{I}}

\newcommand  {\shM}     {\mathscr{M}}

\newcommand  {\shN}     {\mathscr{N}}
\newcommand  {\shL}     {\mathscr{L}}

\newcommand  {\shT}     {\mathscr{T}}

\newcommand  {\ad}      {\operatorname{ad}}

\newcommand  {\Aff}     {\text{{\rm Aff}}}

\newcommand  {\Alb}     {\operatorname{Alb}}

\newcommand  {\Aut}     {\operatorname{Aut}}

\newcommand  {\can}     {{\rm \text{can}}}

\newcommand  {\Cl}      {\operatorname{Cl}}

\newcommand  {\contr}   {\operatorname{contr}}

\renewcommand{\cong}    {\equiv}

\newcommand  {\Der}	{\operatorname{Der}}

\newcommand  {\depth}   {\operatorname{depth}}

\newcommand  {\edim}    {\operatorname{edim}}
\newcommand  {\Enr}      {{\text{\rm Enr}}}

\newcommand  {\Ext}     {\operatorname{Ext}}

\newcommand  {\fppf}    {{\text{{\rm fppf}}}}

\newcommand  {\Frac}    {\operatorname{Frac}}

\newcommand  {\GL}      {\operatorname{GL}}

\newcommand  {\Hom}     {\operatorname{Hom}}

\newcommand  {\id}      {{\operatorname{id}}}

\newcommand  {\I}  	{\text{\rm I}}
\newcommand  {\II}  	{\text{\rm II}}
\newcommand  {\III}  	{\text{\rm III}}
\newcommand  {\IV}  	{\text{\rm IV}}

\newcommand  {\Kernel} 	{\operatorname{Ker}}

\newcommand  {\length}	{\operatorname{length}}
\newcommand  {\Lie}     {\operatorname{Lie}}

\newcommand  {\lra}     {\longrightarrow}

\newcommand  {\maxid}   {\mathfrak{m}}

\newcommand  {\MW}      {\operatorname{MW}}

\newcommand  {\NS}      {\operatorname{NS}}

\newcommand  {\primid}  {\mathfrak{p}}
\renewcommand{\O}       {\mathscr{O}}

\newcommand  {\pdeg}    {\operatorname{pdeg}}

\newcommand  {\Pic}     {\operatorname{Pic}}
\newcommand  {\pd}   	{\operatorname{pd}}
\newcommand  {\PGL}     {\operatorname{PGL}}

\newcommand  {\Proj}    {\operatorname{Proj}}

\newcommand  {\quadand} {\quad\text{and}\quad}

\newcommand  {\ra}      {\rightarrow}

\newcommand  {\rank}    {\operatorname{rank}}
\newcommand  {\red}     {{\operatorname{red}}}

\newcommand  {\Sch}     {{\text{\rm Sch}}}

\newcommand  {\Sing}    {\operatorname{Sing}}

\newcommand  {\Spec}    {\operatorname{Spec}}

\newcommand  {\Tor}     {\operatorname{Tor}}

\newcommand{\xleftrightarrow}[2][]{\ext@arrow 3359\leftrightarrowfill@{#1}{#2}}
\newcommand{\xdashrightarrow}[2][]{\ext@arrow 0359\rightarrowfill@@{#1}{#2}}
\newcommand{\xdashleftarrow}[2][]{\ext@arrow 3095\leftarrowfill@@{#1}{#2}}
\newcommand{\xdashleftrightarrow}[2][]{\ext@arrow 3359\leftrightarrowfill@@{#1}{#2}}
\def\rightarrowfill@@{\arrowfill@@\relax\relbar\rightarrow}
\def\leftarrowfill@@{\arrowfill@@\leftarrow\relbar\relax}
\def\leftrightarrowfill@@{\arrowfill@@\leftarrow\relbar\rightarrow}
\def\arrowfill@@#1#2#3#4{%
  $\m@th\thickmuskip0mu\medmuskip\thickmuskip\thinmuskip\thickmuskip
   \relax#4#1
   \xleaders\hbox{$#4#2$}\hfill
   #3$%
}

\def\mydate{\number\day\space\ifcase\month \or January\or February\or March\or 
April\or May\or June\or July\or
August\or September\or October\or November\or December\fi \space\number\year}

\newcommand{\fpb}	{{(2/\PP^1)}}

\DeclareFontFamily{U}{wncy}{}
\DeclareFontShape{U}{wncy}{m}{n}{<->wncyr10}{}
\DeclareSymbolFont{mcy}{U}{wncy}{m}{n}
\DeclareMathSymbol{\Sh}{\mathord}{mcy}{"58}

\begin{document}

\title[K3-like coverings]{Enriques surfaces with normal K3-like coverings}

\author[Stefan Schr\"oer]{Stefan Schr\"oer}
\address{Mathematisches Institut, Heinrich-Heine-Universit\"at, 40204 D\"usseldorf, Germany}
\email{schroeer@math.uni-duesseldorf.de}

\subjclass[2010]{14J28, 14J27, 14B05, 14L15, 14G17}

\dedicatory{Second revised version, 16 May 2019}

\begin{abstract}
We analyze the structure of simply-connected Enriques surface in characteristic two whose K3-like covering is normal,
building on the work of   Ekedahl, Hyland and Shepherd-Barron.
We develop general methods to construct  such surfaces and the resulting twistor lines
in the moduli stack of Enriques surfaces, including
the case that the K3-like covering is a normal rational surface rather
then a normal K3 surface. 
Among other things, we show that   elliptic double points
indeed do occur. In this case,
there is only one singularity.
The main idea is to apply   flops to    Frobenius pullbacks of   rational elliptic surfaces,
to get the desired K3-like covering.
Our results hinge on Lang's classification of rational elliptic surfaces,
the determination of their Mordell--Weil lattices by Shioda and Oguiso,
and the behavior of unstable fibers under Frobenius pullback via Ogg's Formula.
Along the way, we develop a general theory of Zariski singularities in arbitrary dimension,
which is tightly interwoven with the theory of   height-one group schemes actions
and restricted Lie algebras. Furthermore, we determine  under what conditions
tangent sheaves are locally free, and introduce a theory of canonical coverings 
for arbitrary proper algebraic schemes.
\end{abstract}

\maketitle
\tableofcontents

\section*{Introduction}
\mylabel{Introduction}

A central result in the classification of algebraic surfaces is that
there are only four types  of algebraic surfaces $S$ with   first Chern class $c_1(S)=0$.
These are  the abelian surfaces, the bielliptic surfaces,
the K3 surfaces and the Enriques surfaces. They are distinguished by their second Betti numbers,
which take the respective values
$$
b_2(S)=6,\quad b_2(S)=2,\quad b_2(S)=22\quadand b_2(S)=10.
$$
This holds regardless of the characteristic. Over the complex numbers,  
the bielliptic surfaces are quotients of abelian surfaces
by a finite free group action, and that Enriques surfaces $Y=X/G$ are quotients of K3 surfaces
$X$ by   free involutions. 

This no longer holds true literally when the characteristic $p\geq 0$ of 
the ground field $k$ divides the order   of the group $G$.
A  fundamental insight of Bombieri and Mumford \cite{Bombieri; Mumford 1977}
was that much carries over if   one replaces finite groups  by \emph{finite group schemes}. This comes with
the price that  singularities appear  on the covering $X$, although the
quotient $Y=X/G$ is smooth. In some sense, the non-smoothness of $X$ and $G$ cancel each other.
This phenomenon has numerous applications.
For example, I have  used this effect
to give the correct Kummer construction in characteristic $p=2$,
by replacing an abelian surface by the self-product of  a rational cuspidal curve
\cite{Schroeer 2007}.

For Enriques surfaces $Y$, the numerically trivial part $P=\Pic^\tau_{Y/k}$
of the Picard scheme is a group scheme of order two.
By the Tate--Oort classification of finite group schemes of prime order
\cite{Tate; Oort 1970}, there are three possibilities in characteristic $p=2$.
In the case where the group scheme $P$ is unipotent, one says that $Y$ is
a \emph{simply-connected Enriques surface}.
These surface come along with a torsor $\epsilon:X\ra Y$ with respect
to the Cartier dual $G=\underline{\Hom}(P,\GG_m)$, which here indeed is a local group scheme. 
The torsor takes over the role
of the universal covering, and is therefore called the \emph{K3-like covering}.
This is  very adept terminology, because the K3-like covering has the same cohomological
properties of a K3 surface, yet is never a smooth surface.

The goal of this paper is to analyze the structure of  K3-like coverings $X$, and the ensuing
simply connected Enriques surfaces $Y$.
These where already investigated, among others,  by Bombieri and Mumford \cite{Bombieri; Mumford 1976},
Blass \cite{Blass 1982}, Lang \cite{Lang 1983, Lang 1988},
Cossec and Dolgachev \cite{Cossec; Dolgachev 1989}, Ekedahl and Shepherd-Barron
\cite{Ekedahl; Shepherd-Barron 2004} in the non-normal case,
Ekedahl, Hyland and Shepherd Barron \cite{Ekedahl; Hyland; Shepherd-Barron 2012} in the normal case,
Katsura and Kondo \cite{Katsura; Kondo 2015}, and Liedtke \cite{Liedtke 2015}.
There are, however, many open  foundational questions.
In this paper, we shall concentrate on   \emph{normal K3-like coverings}, although
we take care that our theory works in general. 
An open question was whether non-rational singularities could appear on a K3-like covering $X$.
One main result of this paper is that they do:

\begin{maintheorem}
{\rm (Compare Theorem \ref{existence K3-like II})}
There are   normal K3-like coverings $X$
containing an elliptic double point, which is obtained from the contraction of a rational cuspidal $(-1)$-curve.
\end{maintheorem}

We actually give a systematic way to produce normal K3-like coverings, not based on techniques using
explicit equations and rational vector fields, but on a procedure combining Frobenius base-changes   with flops,
a method borrowed from  the birational geometry of higher-dimensional varieties.
If a non-rational singularity appears, the situation is rather special:

\begin{maintheorem}
{\rm (Compare Theorem \ref{unique elliptic singularity})}
If   a normal K3-like covering $X$ contains  a non-rational singularity,
then there are no other singularities on $X$.
\end{maintheorem}

Each elliptic fibration  $\varphi:Y\ra\PP^1$ on the Enriques surface
comes with a jacobian fibration, which is a rational elliptic surface  $J\ra\PP^1$.
According to the result of Liu, Lorenzini and Raynaud \cite{Liu; Lorenzini; Raynaud 2005}, the Kodaira types
of all fibers $Y_b$ and $J_b$ coincide. The latter were classified by
Persson \cite{Persson 1990} and Miranda \cite{Miranda 1990} over the complex numbers, 
which was extended by Lang in \cite{Lang 1994, Lang 2000} to characteristic two.
There are 110 families with only reduced fibers.
We can show that, with a few uncovered cases, 
all these arise from simply-connected Enriques surfaces:

\begin{maintheorem}
{\rm (Compare Theorem \ref{good mutations exist})}
For each  rational elliptic surface $J\ra\PP^1$
with only reduced fibers, with the possible exception of six cases,
there is a simply-connected Enriques surface $Y$ whose K3-like covering $X$
is birational to the Frobenius pullback $X'=J^\fpb$, 
and having an elliptic fibration $\varphi:Y\ra\PP^1$ whose jacobian is $J\ra\PP^1$.
Moreover, the induced fibration $f:X\ra\PP^1$ induces an injection 
$H^0(X,\Theta_{X/k})\subset H^0(\PP^1,\Theta_{\PP^1/k})$.
\end{maintheorem}

This relies on Shioda's theory of Mordell--Weil lattices \cite{Shioda 1990}
and their classification for rational elliptic surfaces by Oguiso and Shioda \cite{Oguiso; Shioda 1991}.
The passage 
$$
X'\longleftarrow S\lra X
$$ 
from the Frobenius pullback $X'=J^\fpb$
to a K3-like covering $X$ is the fundamental idea of the paper.
It should be seen as a \emph{flop}, a notion coming from the minimal model program for the
classification of higher-dimensional algebraic varieties. Roughly speaking, it describes
the passage from one minimal model to another, without changing the canonical class.

Ekedahl, Hyland and Shepherd Barron   showed 
in their  groundbreaking paper \cite{Ekedahl; Hyland; Shepherd-Barron 2012}
that for simply-connected Enriques surfaces $Y$
whose K3-like covering $X$ has only rational double points,
the tangent sheaf $\Theta_{X/k}=\underline{\Hom}(\Omega^1_{X/k},\O_X)$ is locally free,
in fact isomorphic to $\O_{X}\oplus\O_X$. We generalize this to arbitrary normal K3-like coverings.
Indeed, the main  role of   the flop $X'\leftarrow S\ra X$ is to make the tangent sheaf
trivial.

It turns out that in the restricted Lie algebra $\ideal g= H^0(X,\Theta_{X/k})$,
each vector is $p$-closed, that is, each line is invariant under the $p$-map
$x\mapsto x^{[p]}$. Using the correspondence between finite-dimensional restricted Lie algebras
$\ideal l$ 
and finite height-one groups schemes $G=\Spec(U^{[p]}(\ideal l)^\vee)$, we get 
a rational map 
$$
\xymatrix{
\PP(\ideal g) \ar@{-->}[r]	& \shM_\Enr,	& \ideal l\ar@{|->}[r]	& X/G
}
$$
from the \emph{twistor curve} $\PP(\ideal g)$ into the moduli stack of Enriques surfaces $\shM_\Enr$.
It was introduced in  \cite{Ekedahl; Hyland; Shepherd-Barron 2012}
to study this  moduli stack.
We like to call it the \emph{twistor construction}, because it resembles the 
2-sphere of complex hyperk\"ahler manifolds obtained by rotating the
complex structure while keeping the underlying topological space.
Whether this is more than an  analogy remains to be seen.

The twistor map  is not everywhere defined. 
This is  explained by our  general theory of \emph{Zariski singularities},
which is crucial to understand the non-smoothness of the K3-like coverings
$\epsilon:X\ra Y$ and their generalizations.
We therefore develop a theory of Zariski singularities in full generality for arbitrary dimensions
$n\geq 2$ and characteristics $p>0$.
These are hypersurface singularities $A=R/(g)$, where $R=k[[x_1,\ldots,x_n,z]]$
where  the power series   can be chosen of the form $g=z^p-f(x_1,\ldots,x_n)$.
Using methods from commutative algebra involving free resolutions and projective dimensions,
we get the following criterion for the freeness of the tangent module:

\begin{maintheorem}
{\rm (See Theorem \ref{theta and pd})}
An hypersurface singularity $A=R/(g)$ has free tangent module $\Theta_{A/k}$ of rank $n$ if
and only if $A$ is geometrically reduced, and the module of first-order deformations
$T^1_{A/k}$ has projective dimension
$\pd\leq 2$.
\end{maintheorem}

This automatically holds for the 2-dimensional geometrically reduced Zariski singularities.
Note that it is usually difficult to recognize a Zariski singularity, if it is not
given in the usual form $g=z^p-f(x_1,\ldots,x_n)$.
We therefore develop a homological characterization of Zariski singularities, which
is analogous to Serre's characterization of regularity in local rings.
Roughly speaking:
An hypersurface singularity $A=R/(g)$ Zariski singularity if and only if there
is an $\alpha_2$-action so that the structure sheaf $\O_Z$ of the
closed orbit $Z\subset\Spec(A)$  has finite projective dimension.

Furthermore, we develop easy  tools, involving Tjurina numbers, local class groups, local fundamental groups,
and length formulas    to recognize Zariski singularities.
The following are the Zariski singularities among the rational double points:
$$
A_1, D_{2n}^0, E_7^0, E_8^0\; (p=2),\quad A_2, E_6^0, E_8^0\; (p=3),\quad 
A_4, E_8^0\; (p=5),\quad  A_{p-1}\; (p\geq 7).
$$
Finally, we develop a theory of \emph{canonical coverings} for arbitrary proper algebraic 
schemes $Y$ with $h^0(\O_Y)=1$,
using the Raynaud Correspondence $ H^1(Y_\fppf, G_Y) = \Hom(\hat{G},\Pic^\tau_{Y/k})$.
There is  indeed a canonical choice: Let $G$ be the Cartier dual to the Frobenius kernel
$\Pic_{Y/k}[F]$. To understand the singularities on the resulting canonical coverings
$\epsilon:X\ra Y$, we restrict the covering to closed subscheme $A\subset Y$
that have certain singularities, or on which the $G$-torsor trivializes.
For example, for ADE-configurations $A$ on smooth surfaces $Y$ we introduce
the \emph{fundamental singular locus}, which consist of the singular locus
of the fundamental cycle, and show:
If $G$ is non-reduced, then the $G$-torsor $X\ra Y$ is singular over
the fundamental singular locus of each ADE-configuration.
Under the assumption that $X$ is normal, this gives some crucial restrictions on the geometry
of the curves, in particular if $Y$ is an Enriques surface.

\medskip
The paper is organized as follows:
In Section \ref{Zariski singularities}, we  develop the
theory of Zariski singularities in arbitrary dimensions. Section \ref{Homological characterization}
contains a homological characterization of such singularities.
This is applied to rational double points in Section \ref{Rational double points}.
In Section \ref{Canonical coverings}, we give a general theory of canonical coverings,
based on the Raynaud correspondence.
In Section \ref{Simply-connected enriques}, we review the notion of 
simply-connected Enriques surfaces and discuss their basic properties.
In Section \ref{Trivial tangent sheaf}, we introduced six  conditions
that normal K3-like coverings satisfy, and that indeed allow the construction of
K3-like coverings. The resulting twistor curves are discussed in Section \ref{Twistor curves}.
In Section \ref{K3-like coverings fibrations}, we relate K3-like coverings $X$
to elliptic fibrations  on Enriques surfaces $Y$ and rational surfaces $J$.
To put this to work, we use Ogg's Formula to understand the behavior
of Kodaira types under Frobenius base-change in Section \ref{Ogg's formula}.
In Section \ref{Tate algorithm}, we give a geometric interpretation of the Tate Algorithm,
which computes minimal Weierstra{\ss} equations, in order to describe certain
elliptic singularities.
All this is put into action in Section \ref{Frobenius pullback}, which deals
with the Frobenius pullback of rational elliptic surfaces.
In Section \ref{Lang classification}, we discuss Lang's Classification of
rational elliptic surfaces in characteristic $p=2$, which will be crucial for our goals.
 Section \ref{K3-like with edp} the construction of normal K3-like coverings with
an elliptic singularity.
The absence of other singularities is established in Section \ref{Uniqueness edp}.
In Section \ref{K3-like with rdp} we study normal K3-like coverings with rational singularities.

\begin{acknowledgement}
This paper resulted from  many stimulating discussions with
Shigeyuki Kondo and Toshiyuki Katsura, to whom the
author is very grateful. Large parts  of the paper were
written during two stays of the author at the    Nagoya University, financed by
Shigeyuki Kondo's grant 
JSPS, Grant-in-Aid for Scientific Research (S) No.\ 15H05738.
The author greatly benefited from the ensuing mathematical interchange with Shigeyuki Kondo,
and wishes to thank the Graduate School of Mathematics, Nagoya University    for its hospitality.

I also wish to thank Matthias Sch\"utt for   comments,
Yuya Matsumoto   suggestions and for pointing out some mistakes in the first version,
and the referee for   thorough reading and valuable remarks.
This research was also conducted in the framework of the   research training group
\emph{GRK 2240: Algebro-geometric Methods in Algebra, Arithmetic and Topology}, which is funded
by the DFG.
\end{acknowledgement}

\section{Zariski singularities}
\mylabel{Zariski singularities}

In this section we systematically develop a general theory of Zariski singularities,
a class of local rings introduced by Ekedahl, Hyland and Shepherd--Barron in dimension two
\cite{Ekedahl; Hyland; Shepherd-Barron 2012}. 

Let $k$ be a   ground field, and $A$ be a complete local noetherian ring
such that the canonical map  $k\ra A/\maxid_A$ is bijective. We then get an identification $k=A/\maxid_A$ of
the ground field with the residue field. 
Let $\dim(A)\leq \edim(A)$ be the dimension and embedding dimension.
Equality holds if and only if the 
local noetherian ring $R$ is regular.
We are interested in the following situation:

\begin{proposition}
The   two conditions below are equivalent:
\begin{enumerate}
\item 
The ring $A$ is equidimensional, contains no embedded associated primes, and we have 
$\edim(A)\leq\dim(A)+1$.
\item 
The ring $A$ is isomorphic to $R/(g)$, with $R=k[[x_1,\ldots,x_n,z]]$ and some power series
$g\in R$ that is neither a unit nor zero.
\end{enumerate}
If this holds, then  $ \dim(A)=n$, and the ring $A$ is  Cohen--Macaulay and Gorenstein.
\end{proposition}

\proof
If  (ii) holds, then    Krull's Principal Ideal Theorem ensures that every
irreducible component of $\Spec(A)$ is $n$-dimensional.
Moreover, since $R$ is Cohen--Macaulay and Gorenstein, the same holds for $A=R/(g)$.
In particular, $A$ contains no embedded associated prime.

To see (i)$\Rightarrow$(ii), set $n=\dim(A)$.
Using  the completeness and the assumption on the embedding dimension, we write
$A=R/\ideal a$ for some  ideal $\ideal a$ in the formal power series ring
$R=k[[x_1,\ldots,x_n,z]]$, which is neither the zero nor the unit ideal. According to  \cite{EGA IVd}, Proposition 21.7.2
the subscheme $\Spec(A)\subset\Spec(R)$ is a Weil divisor.
Then $\ideal a$ is generated  by a non-zero non-unit $g\in R$, because the ring $R$ is factorial.
\qed

\medskip
It is convenient to say that the local scheme $\Spec(A)$,
and also the local ring $A$ as well, is  a  \emph{hypersurface singularity} if
the equivalent conditions of the Proposition hold. 
Note that this includes the case where the local ring $A$ is regular, when
we may choose $g=z$.
An isomorphism $A\simeq R/(g)$
as above is called a \emph{presentation}.
Note also that in light of condition (i),
the ring $R'=R\otimes_kk'$ is a hypersurface singularity  if and only if $R$ is a hypersurface
singularity, where $k\subset k'$ is a finite field extension.
Let us write 
$$
J=(g_{x_1},\ldots,g_{x_n},g_z)\subset A
$$
for the  \emph{jacobian ideal} in $A=R/(g)$ generated by residue classes  of the partial derivatives
$g_{x_i}=\partial g/\partial x_i$, $1\leq i\leq n$ and $g_z=\partial g/\partial z$. 
The primes $\primid\subset A$ not containing the jacobian ideal are exactly those
for which the local noetherian rings $A_\primid$ are geometrically  regular over $k$.
By abuse of language, we call the  closed subset   $\Spec(A/J)\subset\Spec(A)$   the \emph{locus of non-smoothness}.
By the preceding observation, this subset does not depend on the choice of presentation $A=R/(g)$.
Note that  
$$
A/J=R/(g,g_{x_1},\ldots,g_{x_n},g_z),
$$
where the ideal on the right is the \emph{Tjurina ideal} in $R$. We observe that 
the ring $A$ is reduced or integral if and only
if the power series $g\in R$ is squarefree or  irreducible, respectively. This holds  because 
$R=k[[x_1,\ldots,x_n,z]]$ is factorial.
The ring $A$ is regular if and only if one of the partial derivatives of $g(x_1,\ldots,x_n,z)$ is a unit,
and then it is even geometrically regular.
We say that the $R$ is a  \emph{geometrically isolated hypersurface singularity} if 
$\dim(A/J)=0$, such that $\Spec(A)$ is geometrically regular outside the closed point.

Now suppose that the ground field $k$ has  positive characteristic $p>0$.
Choose some algebraically closed extension field $\Omega$,
and let $k'=k^{1/p}$ be the intermediate extension comprising all elements $\omega\in \Omega$
with $\omega^p\in k$. Then $k\subset k'$ is a purely inseparable  extension of height one,
and the vector space dimension
$$
\pdeg(k) =  
\dim_k(\Omega^1_{k'/k}) \in \NN\cup\left\{\infty\right\}
$$
is called the \emph{$p$-degree} of $k$. This is also known as the \emph{degree of imperfection},
and can   also be seen as the cardinality of
a $p$-basis  for  $k\subset k^{1/p}$ or equivalently  $k^p\subset k$. 
If $k$ is the field of fractions for a polynomial ring or formal power series ring
in $d$ indeterminates over a perfect field $k_0$, then $\pdeg(k)=n$.
In what follows, we suppose that \emph{our ground field $k$ has finite $p$-degree}, which 
holds in particular for perfect fields.
This assumption will save us some troubles when it comes to K\"ahler differentials
of formal power series, which becomes apparent as follows:

Let $kR^p\subset R$ be the join subring generated by the subrings $k$ and $R^p$.
Clearly, we have an equality 
$R^p=k^p[[x_1^p,\ldots,x_n^p,z^p]]$ inside $R$.
Since the field extension $k^p\subset k$ is finite, the canonical maps
$$
k\otimes_{k^p} k^p[[x_1^p,\ldots,x_n^p,z^p]]\lra kR^p\subset k[[x_1^p,\ldots,x_n^p,z^p]] 
$$
are bijective. It follows that the module of K\"ahler differentials
$$
\Omega^1_{R/k} = \Omega^1_{R/kR^p}
$$
is finitely generated, in fact free, with basis  $dx_1,\ldots,dx_n,dz$.
For ground fields of characteristic zero, or of infinite $p$-degree, the module of K\"ahler differentials is usually
not finitely generated, not even $\maxid_R$-adically separated.

In our situation, the standard exact sequence $gR/g^2R\ra\Omega^1_{R/k}\otimes_RA\ra\Omega^1_{A/k}\ra 0$
becomes
$$
A\lra A^{n+1}\lra\Omega^1_{A/k}\lra 0,
$$
where the map on the right sends the standard basis vectors to the differentials $dx_1,\ldots,dx_n,dz$
and the map on the left is the $(n+1)\times 1$-matrix of partial derivatives given by the transpose of 
$(g_{x_1},\ldots,g_{x_n},g_z)$.
The latter is non-zero if and only if  the ring $A$ is  geometrically reduced. Since $A$ contains no embedded 
associated primes,
the map is then injective. Dualizing the short exact sequence we get an exact sequence
\begin{equation}
\label{theta sequence}
0\lra\Theta_{A/k}\lra A^{\oplus n+1}\lra A\lra\Ext^1(\Omega^1_{A/k},A) \lra 0.
\end{equation}
The term on the right is also called $T^1_{A/k}=\Ext^1(\Omega^1_{A/k},A)$,
the \emph{module of first order deformations}. By the exact sequence, this module coincides with the residue class ring $A/J$, 
whose spectrum is the locus of non-smoothness, and the tangent module
$\Theta_{A/k}=\Hom(\Omega^1_{A/k},A)$ on the left is a syzygy for $A/J$.
This   reveals that the scheme structure on the locus of non-smoothness does not depend on the presentation $A=R/(g)$.
We also immediately get:

\begin{theorem}
\mylabel{theta and pd}
Suppose that $A$ is a  hypersurface singularity. 
Then the tangent module $\Theta_{A/k}$
is free of rank $n=\dim(A)$ if and only if  $A$ is  geometrically reduced and the $A$-module $T^1_{A/k}=A/J$ 
has projective dimension $\leq 2$.
\end{theorem}

The above property is  rather peculiar: Write $A=R/(g)$ with some
power series $g(x_1,\ldots,x_n,z)$. Consider the
\emph{Frobenius power} $J^{[p]}=(g_{x_1}^p,\ldots, g_{x_n}^p,g_z^p)\subset A$ of the jacobian ideal
$J=(g_{x_1}, \ldots,g_{x_n},g_z)\subset A$.
The following condition   was already used in \cite{Schroeer 2008} in dimension two
with the computer algebra system  Magma \cite{Magma}:

\begin{proposition}
\mylabel{theta surface}
Suppose that $A$ is  geometrically isolated   hypersurface singularity
of dimension $n=\dim(A)$.
Then the tangent module $\Theta_{A/k}$  is
free  of rank $n$ if and only if $n\leq 2$ and  the length formula $\length(A/J) = p^n\cdot\length(A/J^{[p]})$ holds.
\end{proposition}

\proof
For geometrically isolated hypersurface singularities, the jacobian ideal $J\subset A$ is $\maxid_A$-primary.
According to a result of Miller (\cite{Miller 2003}, Corollary 5.2.3),
an $\maxid_A$-primary ideal $\ideal a\subset A$ has finite projective dimension if
and only if the length formula $\length(A/\ideal a) = p^n\cdot\length(A/\ideal a^{[p]})$ holds.

It follows that  our condition is sufficient: If $n\leq 2$ and the length formula holds, $A/J$ has finite projective
dimension, which a priori is $\pd(A/J)\leq n\leq 2$, hence $\Theta_{A/k}$ is free of rank $n$ by Proposition \ref{theta and pd}.
The condition is also necessary: If the tangent module
$\Theta_{X/k}$ is free, then the  module $A/J$ has finite projective dimension $\leq 2$,
and the length formula holds by Miller's result. Moreover, the Auslander--Buchsbaum Formula
$\pd(A/J) +\depth(A/J) =\depth(A) $ 
immediately gives $n=\depth(A)=\pd(A/J)\leq 2$.
\qed

\medskip
We say that our complete local ring $A$  is a \emph{Zariski singularity} if it is a hypersurface
singularity  admitting a presentation $A=R/(g)$ with a  power series   of the form 
$$
g(x_1,\ldots,x_n,z)= z^p-f(x_1,\ldots,x_n),
$$
where $f\in k[[x_1,\ldots,x_n]]$ is not a unit. This term was coined   by 
Ekedahl, Hyland and Shepherd-Barron \cite{Ekedahl; Hyland; Shepherd-Barron 2012}
in the case of  isolated surface singularities,
apparently referring to the theory of \emph{Zariski surfaces}, compare Blass and J.\ Lang \cite{Blass; Lang 1987}.

Clearly, the Zariski singularity  $A$ is reduced if and only if $f$ is not a $p$-power. 
Using \cite{AC 8-9}, Chapter VIII, \S7, No.\ 4, Proposition 7, we immediately get:

\begin{proposition}
\mylabel{hilbert-samuel multiplicity}
The  Hilbert--Samuel multiplicity $e(A)\geq 1$ of a Zariski singularity satisfies the 
inequality $ e(A)\leq p$.
\end{proposition}

Now set $A_0=k[[x_1,\ldots,x_n]]$ and $J_0=(f_{ x_1},\ldots,f_{x_n})\subset A_0$.
Then $J=J_0A$ is the jacobian ideal for $A$. This has the following consequence:

\begin{proposition}
\mylabel{theta and depth}
Suppose that $A$ is a geometrically reduced Zariski singularity.
Then the tangent module $\Theta_{A/k}$
is free of rank $n=\dim(A)$ if and only if the ring
$A_0/J_0$ has depth $\geq n-2$. 
\end{proposition}

\proof
Clearly, $A$ is a finite flat algebra over $A_0$ of degree $p$, so  the projective dimension of 
$A/J=A/J_0A=A_0/J_0\otimes_{A_0}A$ as $R$-module
coincides with the projective dimension of $A_0/J_0$ as an $A_0$-module.
By the Auslander--Buchsbaum Formula, we have
$$
\pd(A_0/J_0) + \depth(A_0/J_0) = \depth(A_0) = n,
$$
and the result follows from Proposition \ref{theta and pd}.
\qed

\medskip
In dimension $n=2$, the condition becomes vacuous, and we get:

\begin{corollary}
\mylabel{theta for zariski surface}
The tangent module $\Theta_{A/k}$ of any two-dimensional geometrically reduced Zariski singularity is free of rank two.
\end{corollary}

If $A$ is a geometrically isolated hypersurface singularity, that is, the $A$-module $T^1_{A/k}$ has finite length,
its length $\tau=\length T^1_{A/k}$ is called the the \emph{Tjurina number}.
This is also the colength of the jacobian ideal $J\subset A$, or equivalently the colength of the
Tjurina ideal in $R$. For Zariski singularities, we get:

\begin{proposition}
\mylabel{tjurina for isolated zariski}
If $A$ is  a geometrically isolated Zariski singularity given by a formal power series
$g=z^p-f(x_1,\ldots,x_n)$, then the Tjurina number satisfies   
$$
\tau=p\cdot \length k[[x_1,\ldots,x_n]]/(f_{x_1},\ldots, f_{x_n}).
$$
In particular, the Tjurina number $\tau\geq 0$ is a multiple of the characteristic $p>0$.
\end{proposition}

Finally, we state some useful facts on  Zariski-singularities pertaining to     fundamental groups
and class groups:

\begin{proposition}
\mylabel{coverings for normal zariski}
Let  $A$ be a normal Zariski singularity, and 
$A\subset B$ be a finite ring extension, with $\Spec(B)$ normal and connected. Assume that  $\Spec(B)\ra\Spec(A)$
is \'etale in codimension one. 
Then there is a finite separable field extension $k\subset k'$ with $B=A\otimes_kk'$.
\end{proposition}

\proof
Write $A=R/(g)$ as usual, and set $A_0=k[[x_1,\ldots,x_n]]$.
The resulting morphism $\Spec(A)\ra\Spec(A_0)$ is a universal homeomorphism that is 
finite and flat of degree $p$.
Let $U_0\subset\Spec(A_0)$ be the open subset corresponding to the regular locus $U\subset\Spec(A)$.
According to \cite{SGA 1}, Expos\'e IX, Theorem 4.10 every finite \'etale covering $U'\ra U$ is the
pullback of some finite \'etale covering $U_0'\ra U_0$.
The latter extends to $\Spec(A_0)$, by the Zariski--Nagata Purity Theorem
\cite{SGA 2}, Expos\'e X, Theorem 3.4.
Since the complete local ring $A_0$ is henselian, the extension must have the form $A_0'=A_0\otimes_kk'$
for some finite separable field extension $k\subset k'$. Consider the  pullback $A'=A\otimes_{A_0} A_0'=A\otimes_kk'$.
Then $\Spec(A')$ is regular in codimension one and Cohen--Macaulay, thus normal.
By construction, the two morphisms $\Spec(B)\ra\Spec(A)$ and $\Spec(A')\ra\Spec(A)$ coincide
over $U\subset\Spec(A)$, whose complement has codimension $d\geq 2$.
Now \cite{Hartshorne 1994}, Theorem 1.12 ensures that there is an isomorphism $B\simeq A'$ of $A$-algebras.
\qed

\begin{corollary}
\mylabel{components rational}
Suppose that $A$ is a normal two-dimensional Zariski singularity, with separably closed ground field $k$
and resolution of singularities $X\ra \Spec(A)$. Then each irreducible component $E_i$
of the exceptional divisor $E\subset X$ is a rational curve.
\end{corollary}

\proof
Fix some prime $l\neq p$, and suppose that some $E_i$ is non-rational.
Consequently one  finds an invertible sheaf $\shL_E\not\simeq\O_E$ on $E$ with $\shL_E^{\otimes l}\simeq\O_E$.
The short exact sequences $0\ra\O_E(-nE)\ra\O_{(n+1)E}^\times\ra\O_{nE}^\times\ra 1$
induces long exact sequences
$$
H^1(E,\O_E(-nE))\lra \Pic((n+1)E)\lra \Pic(nE) \lra 0
$$
where the term on the left is a  $p$-group.
It follows that there is an  invertible sheaf $\shL$ on $X$, restricting to $\shL_E$, with  $\shL^{\otimes l}\simeq \O_X$.
Then the ensuing $\O_X$-algebra $\shA=\O_X\oplus\shL^{\otimes -1}\oplus\ldots\oplus\shL^{\otimes 1-l}$ of rank $l>1$
defines a ring extension $A\subset B$ of degree $l$ as in the Proposition.
However, since $k$ is separably closed, we must have $B=A$, contradiction.
\qed

\begin{proposition}
\mylabel{class group for normal zariski}
Let  $A$ be a normal Zariski singularity.
Then the class group $\Cl(A)$ is an abelian $p$-group,
and every element is annihilated by $p^{d+n-1}$, where $d=\pdeg(k)$ and $n=\dim(A)$.
\end{proposition}

\proof
Write $A=k[[x_1,\ldots,x_n,z]]/(g)$ for some power series $g=z^p-f(x_1,\ldots,x_n)$,
and consider the subring $A_0=[[x_1,\ldots,x_n]]$ and the finite ring extension $A_0\subset B$
with $B=k^{1/p}[[x_1^{1/p},\ldots,x_n^{1/p}]]$. The   inclusion $A_0\subset B$ is finite
and flat, of degree $p^{n+d}$. Clearly, the power series $f\in A_0$ admits a $p$-th root $h\in B$,
which gives an integral homomorphism $A\ra B$ of $A_0$-algebras via $z\mapsto h$. This map must be injective,
because $A$ is integral and $\dim(A)=\dim(B_0)$. The composite inclusion $A_0\subset A\subset B$
has degree $p^{d+n}$, and the first inclusion has degree $p$, whence $A\subset B$ has
degree $d=p^{d+n-1}$. Let $U\subset\Spec(A)$ be the regular locus, and $V\subset \Spec(B)$ be its
preimage. Then  $\Cl(A)=\Pic(U)$ and $ \Cl(B)=\Pic(V)=\Pic(B)=0$.
The composition $\Pic(U)\ra\Pic(V)$ with the norm map $\Pic(V)\ra\Pic(U)$
is multiplication by $\deg(V/U)=p^{d+n-1}$, and the result follows.
\qed

\begin{corollary}
\mylabel{intersection matrix for surface zariski}
Suppose that $A$ is a normal two-dimensional Zariski singularity,   
with resolution of singularities $X\ra \Spec(A)$. Let $E\subset S$ be the exceptional divisor,
$E=E_1\cup\ldots\cup E_r$ be the irreducible components, and $N=(E_i\cdot E_j)$ the   
intersection matrix. If the ground field 
$k$ is algebraically closed, then the Smith group $\ZZ^r/N\ZZ^r$ is annihilated by $p$.
In particular, the  absolute value of $\det(N)$ is a $p$-power.
\end{corollary}

\proof
The Smith group $\ZZ^r/N\ZZ^r$ is finite, because the intersection matrix is negative-definite.
Any effective Weil divisor on $\Spec(A)$ can be seen as an effective Weil divisor 
$D\subset X$ that contains no exceptional curve $E_i$,
and we get intersection numbers $(D\cdot E_i)$.
This yields a well-defined homomorphism $\Cl(A)\ra\ZZ^r/N\ZZ^r$.
In light of the Proposition, we merely have to check that the map is surjective.
Each element of the Smith group is represented by a difference of vectors with
non-negative entries.
Suppose we have a vector  of non-negative integers $(n_1,\ldots,n_r)\in\ZZ^r$.
As the ground field $k$ is algebraically closed, there is an effective Cartier divisor $D_0\subset E$
such that $\shL_0=\O_E(D)$ has intersection numbers $(\shL_0\cdot E_i)=n_i$.
Since $A$ is complete and in particular henselian, we may extend it to some invertible sheaf $\shL$ on $X$
by \cite{EGA IVd}, Corollary 21.9.12.
This yields an element in $\Cl(A)$ having the same class as $(n_1,\ldots,n_r)$
in the Smith group $\ZZ^n/N\ZZ^r$.
\qed

\section{Homological characterization}
\mylabel{Homological characterization}

Since  a complete local noetherian ring $A$ may arise in very many different ways,
it is a priori  not clear how to decide whether or not   it is  a  Zariski singularity,
even if it comes with a presentation $A=R/(g)$.
In this section, we establish a homological characterization
using $p$-closed derivations and actions of height one group schemes
that is both   intrinsic and practical.
We start by recalling the equivalence between finite groups schemes of height one
and finite-dimensional restricted Lie algebras. For details, we refer to the monograph of
Demazure and Gabriel \cite{Demazure; Gabriel 1970}, Chapter II, \S7.
A very brief summary   is contained in  \cite{Schroeer 2007}, Section 1.

Fix a ground field $k$  of characteristic $p>0$.
Let  $\ideal g$ be a \emph{restricted Lie algebra}, with Lie bracket $[x,y]$ and
$p$-map $x^{[p]}$. The latter satisfies the semilinearity property  $(\lambda x)^{[p]}=\lambda^px^{[p]}$,
and the two Jacobson Formulas  
$$
\ad(x^{[p]})=\ad(x)^p\quadand (x+y)^{[p]} = x^{[p]} + y^{[p]} + \sum_{r=1}^{p-1} s_r(x,y),
$$
where $\ad(x)(y)=[x,y]$ is the adjoint representation, and $s_r(x,y)$ is some universal expression in Lie brackets of $x$ and $y$,
for example explained in \cite{Demazure; Gabriel 1970}, Chapter II, \S7, Definition 3.3.
For more details, see the book of Strade and Farnsteiner \cite{Strade; Farnsteiner  1988}.

Let $U(\ideal g)$ be the universal enveloping algebra,  
$U^{[p]}(\ideal g)$ be the quotient by the ideal generated by 
the elements $x^p-x^{[p]}\in U(\ideal g)$ with $x\in\ideal g$, and $U^{[p]}(\ideal g)^\vee$ the dual  vector space.
The  diagonal  map $\ideal g\ra \ideal g\oplus \ideal g$
induces the comultiplication $\Delta: U^{[p]}(\ideal g)\ra U^{[p]}(\ideal g)\otimes U^{[p]}(\ideal g)$
and dually a  multiplication in $U^{[p]}(\ideal g)^\vee$. The latter becomes a commutative $k$-algebra,
and we may form
$$
G=\Spec(U^{[p]}(\ideal g)^\vee).
$$
As explained in loc.\ cit.\ Proposition 3.9, this acquires the structure of  a group scheme of height one, 
with $\Lie(G)=\ideal g$. Its $R$-valued points $G(R)$ is the group of elements
$g\in U^{[p]}(\ideal g)\otimes_k R$ satisfying $\Delta(g)=g\otimes g$. 
Such elements $g$ are also called \emph{group-like}. Furthermore, the correspondences
$$
G\longmapsto \Lie(G)\quadand \ideal g\longmapsto \Spec(U^{[p]}(\ideal g)^\vee)
$$
are adjoint equivalences between the categories of finite  groups schemes
of height one and the category of finite-dimensional restricted Lie algebras,
by loc.\ cit.\ Proposition 4.1.

If $X$ is a scheme, the right $G$-actions
$X\times G\ra X$ correspond to 
the homomorphisms $\ideal g=\Lie(G)\ra\Der_k(\O_X,\O_X)$   
of restricted Lie algebras. Here 
$$
\Der_k(\O_X,\O_X)=\Hom_{\O_X}(\Omega^1_{X/k},\O_X) = H^0(X,\Theta_{A/k})
$$
is the restricted Lie algebra of   derivations $D:\O_X\ra \O_X$, where the 
$p$-map  is  the $p$-fold composition $D^p=D\circ\ldots\circ D$ in the associative algebra
of all differential operators $\O_X\ra\O_X$.
In the case that  $X=\Spec(A)$ is affine, this boils down to  specify a    derivation  $D:A\ra A$.

If the restricted Lie algebra $\ideal g$ is one-dimensional, with basis vector 
$x\in\ideal g$, then the Lie bracket must be trivial,
and the $p$-map is determined via $x^{[p]}=\lambda x$ by some  unique scalar $\lambda\in k$.
Let us write $\ideal g_\lambda$ and $G_\lambda$ for the resulting restricted Lie algebra
and height-one group scheme. 
The  right actions of $ G_\lambda$ on the scheme $X$ thus correspond to 
derivations $D:A\ra A$ satisfying $D^p=\lambda D$. Such derivations are called
\emph{$p$-closed}. 

To make this more explicit, write $U^{[p]}(\ideal g_\lambda)= k[x]/(x^p-\lambda x)$, and
choose as a vector space basis the divided powers $u_i=x^i/i!$ with $0\leq i\leq p-1$.
This basis has the advantage that the comultiplication is given by $\Delta(u_i)=\sum_{r+s=i}u_r\otimes u_s$.
Denote the dual basis by  $t_i\in U^{[p]}(\ideal g_\lambda)^\vee$.
Inside the ring $U^{[p]}(\ideal g_\lambda)^\vee$, we get $t_r\cdot t_s=t_{t+s}$, and in particular
$t_1^p=0$. With $t=t_1$, one obtains
$U^{[p]}(\ideal g_\lambda)^\vee=k[t]/(t^p)$.
Consequently, the action $  X\times G_\lambda\ra X$   takes the explicit form
\begin{equation}
\label{action formula}
A\lra A\otimes_k U^{[p]}(\ideal g_\lambda)^\vee,\quad
f\longmapsto \sum_{i=0}^{p-1}\frac{D^i(f)}{i!}\otimes t^i.
\end{equation}
Now suppose we have a right $G_\lambda$-action on an affine scheme  $X=\Spec(A)$, 
and let  $a\in\Spec(A)$ be a $k$-rational point, corresponding to 
a maximal ideal $\maxid \subset A$. Let 
$\ideal a\subset A$ be the  intersection of the kernels for the $k$-linear maps
\begin{equation}
\label{action kernels}
\maxid \stackrel{D^i}{\lra} A\lra A/\maxid ,\quad 1\leq i\leq p-1.
\end{equation}
One easily checks that this  vector subspace $\ideal a\subset A$ is   an ideal. We call it the \emph{orbit ideal}
for the   $G_\lambda$-action, or the corresponding $p$-closed derivation $D:A\ra A$.

\begin{proposition}
\mylabel{fixed and free points}
The closed subscheme
$X'=\Spec(A/\ideal a)$ is a $G_\lambda$-invariant Artin scheme
consisting  only of the point $a\in X$, and its length is either $l=1$ or $l=p$.
In the latter case, $X'$ is a trivial $G_\lambda$-torsor over $\Spec(k)$, and their is 
an open neighborhood $U\subset X$ of the point $a\in X$ on which the $G_\lambda$-action is free.
\end{proposition}

\proof
For each $f\in\maxid $, we obviously have $f^p\in\ideal a$, whence $X'=\{a\}$ holds as a set.
From \eqref{action formula}  we infer that $X'\subset X$ is invariant under the $G_\lambda$-action.
To compute its length, we may assume that $\ideal a=0$, by replacing  $X$ with $X'$.
Then the  composite map 
$$
A\lra k\otimes U^{[p]}(\ideal g_\lambda)^\vee=k[t]/(t^p)
$$
obtained from \eqref{action formula}  is  injective, and thus $\length(A)\leq p$.
Furthermore, we may regard $A$ as a subalgebra of $B=k[t]/(t^p)$. 
If $\maxid_A=0$, then $\length(A)=1$.
Now suppose that $\maxid_A\neq 0$. Fix some non-zero
$f\in\maxid_A$. Since $A\subset B$, there is some $0\leq i\leq p-1$ so that $D^i(f)$
does not vanish in $A/\maxid_A=k$. Hence $D^i(f)\in A$ is a unit.
This already ensures  that $\Spec(A)\ra\Spec(A^D)$ is a torsor, for example
by \cite{Restuccia; Schneider 2003}, Theorem 4.1.
Here $A^D$ is  the kernel of the derivation $D:A\ra A$, which is a $k$-subalgebra $A^D\subset A$.
It follows that $\length(A)\geq p$. In light of the inclusion $A\subset B$,
we must have $\length(A)=p$.

Now suppose that $\length(A)=p$.  Then $X\times G\ra X$ is a $G$-torsor
containing a rational point, so the torsor is trivial. 
To show that there is an open neighborhood on which the action is free,
we revert to the original situation $X=\Spec(A)$, with $X'=\Spec(A/\ideal a)$.
Let $A^D=\Kernel(D)$ be the kernel of the additive map $D:A\ra A$, which  is a subring of $A$. 
Since $\ideal a\neq\maxid$, there is some $f\in\maxid$ so that  $D^i(f)$   is non-zero in $A/\maxid$,
whence a unit in the local ring $A_\primid$. Replacing $A$ by some suitable localization 
$A_g$, we may assume that $D^i(f)\in A$ is a unit.
This already ensures  that $\Spec(A)\ra\Spec(A^D)$ is a torsor, again
by \cite{Restuccia; Schneider 2003}, Theorem 4.1.
\qed

\medskip
The kernel of the derivation $D:A\ra A$ is a subring, and we call this subring
$A^D\subset A$  the \emph{ring of invariants}.
Then $X/G_\lambda=\Spec(A^D)$ is a categorical quotient in the category $(\Aff/k)$ of
all affine  schemes, and in the category of all  schemes $(\Sch/k)$ as well. 
If the action is free,
that is, the canonical map $X\times G\ra X\times_{X/G_\lambda} X$ given by $(a,g)\mapsto (a,ag)$
is an isomorphism, and
the canonical morphism $X\ra X/G_\lambda$ is a $G_\lambda$-torsor.

We say that a rational point $a\in X=\Spec(A)$ is a \emph{fixed point} if
the orbit ideal $\ideal a\subset A$ is maximal. Otherwise, the  $G_\lambda$-action is free in an open neighborhood.
Being an $A$-module, the orbit ideal has a projective dimension $\pd(\ideal a)\geq 0$, which may be finite or infinite.
From this we get the desired \emph{homological characterization} of Zariski singularities:

\begin{theorem}
\mylabel{homological characterization}
Suppose that $A$ is a complete local noetherian ring with $A/\maxid_A=k$.
\begin{enumerate}
\item 
If $A$ is a Zariski singularity, there is a derivation $D:A\ra A$ satisfying $D^p=0$  such
that the orbit ideal $\ideal a\subset A$ is not the maximal ideal, and has finite projective dimension.
\item
If there is a $p$-closed derivation $D:A\ra A$ such that the orbit ideal
is not the maximal ideal and has finite projective dimension, then $A\otimes_kk'$ is
a Zariski singularity for some finite separable extension $k\subset k'$.
\end{enumerate}
\end{theorem}

\proof
(i)  If $A=R/(g)$ is a Zariski singularity, say given by the 
power series $g=z^p-f(x_1,\ldots,x_n)$,
we may take the standard derivation $D=D_z$, which satisfies $D^p=0$. 
Clearly, $x_1,\ldots,x_n,z^p$ are contained in the orbit ideal $\ideal a\subset A$,
but not $z$, because $D(z)=1$ is nonzero in $A/\maxid_A$.
Thus the orbit ideal is not maximal, and by Proposition \ref{fixed and free points} we must have
$\ideal a=(x_1,\ldots,x_n,z^p)$.
This ideal is induced from the maximal ideal for the finite flat extension $k[[x_1,\ldots,x_n]]\subset A$, whence
has finite projective dimension.

(ii) Write $D^p=\lambda D$ for some scalar $\lambda\in k$. Choose some $\mu=\lambda^{1/(p-1)}$ in some
algebraic closure. This gives a finite separable extension
$k'=k(\mu)$. Replacing $A$ by $A\otimes_kk'$ we may assume $k=k'$.
If $\lambda\neq 0$, we replace $D$ by $\mu D$.
The upshot is that either $D^p=D$ or $D^p=0$.
Then the corresponding height one group scheme  is   $G=\mu_p$ or $G=\alpha_p$, respectively.
According to Proposition \ref{fixed and free points}, the $G$-action on $X$ is free.
Since the ring of invariants $A^D$ contains the ring of $p$-th powers $A^p$,
the projection $X\ra X/G$ is a universal homeomorphism, whence the ring of invariants $A^D$ is local.
According to the Eakin--Nagata Theorem, or equivalently by faithful flatness, it is also noetherian.
It then also must be complete. It follows from  \cite{Schroeer 2007}, Proposition 2.2
that the invariant ring is also regular. 

Choose some regular system of parameters $x_1,\ldots,x_n\in A^D$, and write the invariant ring as
$A^D=k[[x_1,\ldots,x_n]]$.
Now set $S=\Spec(A^D)$, and suppose that $G=\mu_p$. 
From the short exact sequence $0\ra G\ra\GG_m\stackrel{F}{\lra}\GG_m\ra 0$,
we get an exact sequence
$$
H^0(S,\GG_m)\lra H^1(S,G)\lra H^1(S,\GG_m).
$$
The term on the right vanishes, because $\Pic(S)=0$. Whence the total space $X=\Spec(A)$ of the torsor
takes the form $A=A^D[z]/(z^p-f)$
for some invertible formal power series $f(x_1,\ldots,x_n)$. Its constant term $\alpha\in k$
must be a $p$-power, because  $A^D\subset A$ has trivial residue field extension.
Replacing $z$ by $ z-\beta$ with $\beta^p=\alpha$ reveals that $A$ is a Zariski singularity.
In the case $G=\alpha_p$, one uses the short exact sequence
$0\ra G\ra\GG_a\stackrel{F}{\ra}\GG_a\ra 0$ and argues similarly.
\qed

\medskip
According to Proposition \ref{theta for zariski surface}, the tangent module
of a 2-dimensional geometrically reduced Zariski singularity is free.
It is possible to give an explicit basis:

\begin{proposition}
\mylabel{basis tangent module}
Let $A=k[[x,y,z]]/(g)$ be a two-dimensional geometrically isolated Zariski singularity, 
given by an formal power series $g=z^p-f(x,y)$. Then the derivations
\begin{equation}
\label{basis tangent}
D_z\quadand f_yD_x - f_xD_y
\end{equation}
form an $A$-basis for the tangent module $\Theta_{A/k}=\Der_k(A,A)=\Hom(\Omega^1_{A/k},A)$. 
\end{proposition}

\proof
In our situation, the exact sequence \eqref{theta sequence} becomes
$$
0\lra\Theta_{A/k}\lra A^{\oplus 3}\stackrel{(f_x,f_y,0)}{\lra} A.
$$
Clearly, the transpose of the  vectors $(0,0,1)$ and $(f_y,-f_x,0)$ lie in the kernel of the map on the right,
so the   derivations in \eqref{basis tangent} indeed can be regarded as elements of the tangent module $\Theta_{A/k}$.
Moreover, $D_z\in\Theta_{A/k}$ generates a free  direct summand,
and  the kernel of the linear map $(f_x,f_y):A^2\ra A$ gives a complement.

By assumption, $A$ is not geometrically regular, such that $f_x,f_y\in\maxid_A$.
Regard $A$  as a finite flat algebra over $A_0=k[[x,y]]$. Since the jacobian ideal $(f_x,f_y)$
is $\maxid_{A_0}$-primary, the elements $f_x,f_y\in A_0$ form
a regular sequence. Thus the Koszul complex
$$
0\lra A_0\stackrel{\binom{f_y}{-fx}}{\lra} A_0^2 \stackrel{(f_x,f_y)}{\lra} A_0
$$
is exact. Tensoring with the flat algebra $A$ shows that the kernel of $(f_x,f_y):A^2\ra A$
is free of rank one, with generator $f_yD_x - f_xD_y$.
\qed

\begin{proposition}
\mylabel{free actions}
Let $A=k[[x,y,z]]/(g)$ be a two-dimensional geometrically isolated Zariski singularity,
given by an formal power series $g =z^p-f(x,y)$. Let 
$$
D=uD_z+v(f_yD_x - f_xD_y)\in\Theta_{A/k}
$$
be a $p$-closed derivation, for some coefficients $u,v\in A$.
Then the corresponding $G$-action on $\Spec(A)$ is free if and only if $u\in A^\times$.
If this is the case, the orbit ideal is given by $\ideal a=(x,y)$, and the ring of invariants
$A^D$ is regular.
\end{proposition}

\proof
The orbit ideal $\ideal a\subset A$ obviously contains $z^p$.
Since the local ring $A$ is geometrically singular, we have $f_x,f_y\in\maxid_A$. Thus also $x,y\in\ideal a$.
If the coefficient $u\in A$ is not invertible, that is, $u\in\maxid_A$, we also have $z\in\ideal a$,
and the $G$-action on $\Spec(A)$ is not free.
Now suppose that $u\in A^\times$. Then $D(z)=u$ does not vanish in the residue field, and thus $z\not\in\ideal a$.
Consequently, the $G$-action on $\Spec(A)$ is free, and the inclusion $(x,y,z^p)\subset\ideal a$ is an equality.
In light of the relation $z^p=f(x,y)$ in $A$, we get $\ideal a=(x,y)$, and this has finite projective dimension.
It follows from  \cite{Schroeer 2007}, Proposition 2.2
that the invariant ring $A^D$ is   regular. 
\qed

\medskip
This is quite remarkable, because the  ideal $\ideal a=(x,y)$ in the ring $A$ neither depends on 
the choice of the regular system of parameters  $x,y,z\in R$, nor the choice of the formal power series  $g=z^p-f(x,y)$, 
nor the choice of the  $p$-closed derivation $D\in\Theta_{A/k}$ giving a free group scheme action.
We simply call it the \emph{orbit ideal}
of the two-dimensional  geometrically reduced Zariski singularity. It defines a
canonical zero-dimensional subscheme  of embedding dimension one on the singular scheme $X=\Spec(A)$,
called the \emph{orbit subscheme}.
In turn, it yields a unique tangent vector 
\begin{equation}
\label{orbit direction}
\Spec k[\epsilon]\subset \Spec k[[x,y,z]]/(x,y,z^p)\subset\Spec(A)
\end{equation}
that indicates the direction of any $p$-closed vector field corresponding to  a free action. 
We call it the
\emph{orbit direction}. Here and throughout, $\epsilon $ denotes
an indeterminate satisfying $\epsilon^2=0$.

\begin{theorem}
\mylabel{canonical line}
Let $A$ be a two-dimensional  geometrically reduced Zariski singularity that is not regular,
and $\ideal g\subset\Theta_{A/k}$ be a restricted Lie subalgebra with $\dim_k(\ideal g)=2$
for which every element is $p$-closed, 
and $A\ideal g=\Theta_{A/k}$. Then there is a line $\ideal l\subset \ideal g$ so that
for each element $D\in\ideal g$, the corresponding group scheme
action on $\Spec(A)$ is free if and only if $D\not\in\ideal l$.
If these conditions holds, the invariant ring $A^D$ is regular.
\end{theorem}

\proof
Suppose $A$ is given by the formal power series $g=z^p-f(x,y)$. Choose a $k$-basis $D_1,D_2\in \ideal g$.
Then $D_1,D_2\in \Theta_{A/k}$ form an $A$-basis as well.
According to Proposition \ref{free actions}, the $D_z,f_yD_x - f_xD_y\in\Theta_{A/k}$ yield another
$A$-basis. Write
$$
D_z= r_{11}D_1 + r_{12}D_2\quadand f_yD_x - f_xD_y= r_{21}D_1 + r_{22}D_2
$$
for some base-change matrix $(r_{ij})\in\GL_2(A)$.
Let $\bar{r}_{ij}\in A/\maxid_A=k$ be the residue classes of the matrix entries,
and let $\ideal l\subset\ideal g$ be the line generated by the vector $\bar{r}_{21}D_1 + \bar{r}_{22}D_2$.
Each   derivations $D=\lambda D_1+\mu D_2\in\ideal g$ with coefficients $\lambda,\mu\in k$ can be written as 
$D=uD_z+v(f_yD_x - f_xD_y)$ with coefficients $u,v\in A$.
Obviously, $D\not\in\ideal l$ means that $u\in A^\times$.
The assertion thus follows from Proposition \ref{free actions}.
\qed

\medskip
In the above situation, we call $\ideal l\subset\ideal g$ the
\emph{canonical line}.
This also can be regarded in the following way:
The composite  map $\ideal g \subset\Theta_{A/k}\ra \Theta_{A/k}\otimes k$ is bijective.
Under the canonical map
$$
\Theta_{A/k}\otimes k=\Hom_A(\Omega^1_{A/k},A)\otimes_Ak \lra \Hom_k(\Omega^1_{A/k}\otimes k,k) = \Hom_A(\Omega^1_{A/k},k),
$$
the elements $D\in\ideal g$ are turned into derivations $D:A\ra k$,
which correspond to homomorphisms $A\ra k[\epsilon]$.
Note, however, that the vector spaces above have dimension two and three,
and that the canonical line $\ideal l\subset\ideal g$ is the kernel of the above linear map.
Hence the image is one-dimensional, and indeed induces the orbit direction
$\Spec k[\epsilon] \subset\Spec(A)$.

In Theorem \ref{canonical line}, the 2-dimensional restricted Lie algebra $\ideal g$
has the property that every vector is $p$-closed.
Up to isomorphism, there are exactly two  such restricted Lie algebras of dimension two:

\begin{lemma}
\mylabel{classification lie algebras}
Suppose $k$ is algebraically closed.
Then the 2-dimensional restricted Lie algebras $\ideal g$ for which every vector is $p$-closed
are up to isomorphism given by $\ideal g=kx\oplus ky$ with
$$
[x,y] = x^{[p]}=y^{[p]}=0 \qquad\text{or}\qquad [x,y]=y,\, x^{[p]}=x,\, y^{[p]}=0.
$$
\end{lemma}

\proof
This follows from the classification of 2-dimensional restricted Lie algebras \cite{Wang 2013}, Proposition A3.
\qed

\section{Rational double points}
\mylabel{Rational double points}

Let $k$ be an algebraically closed  ground field of characteristic $p>0$, and $A$ be a complete local ring
that is normal and of dimension two. We say that $A$ is a \emph{rational double point}
if it is a rational singularity of Hilbert--Samuel multiplicity $e(A)=2$.
Equivalently, it is a rational Gorenstein singularity.
Recall that in characteristic $p=2$, the \emph{rational double points} are classified into
$$
A_n,\quad  D_{2n}^r, D_{2n+1}^r,\quad  E_6^0,E_6^1,\quad  E_7^0,\ldots,E_7^3,\quadand E_8^0,\ldots, E_8^4,
$$
with $0\leq r\leq n-1$, according to Artin's analysis \cite{Artin 1977}.
We now can determined which of them are Zariski:

\begin{theorem}
\mylabel{zariski rdp}
In characteristic $p=2$, the rational surface singularities that are Zariski singularities are
precisely the rational double points of type $A_1$, $D_{2m}^0$, $E_7^0$ and $E_8^0$.
\end{theorem}

\proof
It follows from the defining equations in \cite{Artin 1977} that the given rational double points
are Zariski.
The converse is not at all obvious.
Let $A=k[[x,y,z]]/(g)$ be a rational  Zariski singularity, given by the formal power series $g=z^2-f(x,y)$. 
Being a hypersurface singularity, it  is Gorenstein, hence a rational double point.
According to Corollary \ref{intersection matrix for surface zariski},
the intersection matrix $N=(E_i\cdot E_j)$ for the exceptional divisor 
on a resolution of singularities has $\det(N)=\pm 2^\nu$ for some $\nu\geq 0$, 
whence  only type $A_n$ with $n =2^\nu-1$, $D_n$, $E_7$ and $E_8$ are possible.
Furthermore, this Smith group $\Phi=\ZZ^r/N\ZZ^r$ is anihilated by $p=2$, which rules
out $A_n$ with $n\geq 2$, and $D_n$ with $n=2m+1$.

Clearly, the module of K\"ahler differentials  $\Omega^1_{A/k}$ is generated by $dx,dy,dz$,
and $Adz\subset\Omega_{A/k}$ is an invertible direct summand.
According to \cite{Schroeer 2008}, proof of Theorem 6.3 such an invertible direct summand
does not exist for $D_n^r$ with $r>0$. Likewise, in the proof for loc.\ cit.\ Theorem 6.4
the same is  shown for the singularities $E_7^r$ and $E_8^r$ with $r>0$.
\qed

\medskip
Using Proposition \ref{intersection matrix for surface zariski}, one shows:

\begin{proposition}
\mylabel{zariski rdp other p}
In characteristic $p=3$, the rational surface singularities that are Zariski singularities
are the rational double points of type $A_2, E_6^0,E_8^0$.
In characteristic $p=5$, these are $A_4$ and $E_8^0$.
For $p\geq 7$, there is only $A_{p-1}$.
\end{proposition}

A rational double point that is also a Zariski singularity is called a \emph{Zariski rational double point}.
Going through their normal form \cite{Artin 1977}, one obtains the following observation:

\begin{proposition}
\mylabel{tjurina rdp}
The Tjurina numbers for   Zariski rational double points 
of type $A_n$, $D_n$ or $E_n$ in characteristic $p>0$ are $\tau=pn/(p-1)$. 
In particular, $\tau=2n$ in characteristic two.
\end{proposition}

\section{Canonical coverings of algebraic schemes}
\mylabel{Canonical coverings}

The goal of this section is to  establish some general facts about ``canonical coverings'' of proper algebraic schemes.
We will apply this later to the K3-like covering of simply-connected Enriques surface, but the
underlying principles have a  far more general relevance.

Fix a ground field $k$, and let $Y$ be a proper $k$-scheme
with $k=H^0(Y,\O_Y)$, and $\Pic_{Y/k}$ be its Picard scheme.
Let $G$ be a finite commutative group scheme, and $\hat{G}=\ul \Hom(G,\GG_m)$ be the Cartier dual. 
This  is also a finite commutative group scheme, of the same order  $h^0(\O_G)=h^0(\O_{\hat{G}})$,
and the biduality map  
$$
G\lra\ul\Hom(\ul\Hom(G,\GG_m),\GG_m),\quad g\longmapsto (f\mapsto f(g))
$$ 
is an isomorphism of group schemes. 

Set $S=\Spec(k)$,   let $f:Y\ra S$ be the structure morphism, and write $G_Y=f^*(G)$ for the induced
relative group scheme over $Y$.
According to \cite{Raynaud 1970}, Proposition 6.2.1 applied to the Cartier dual $\hat{G}$, we have a a canonical bijection 
$$
R^1f_*(G_Y) \lra \underline{\Hom}(\hat{G},\Pic_{Y/k}) = \underline{\Hom}(\hat{G},\Pic^\tau_{Y/k})
$$
of abelian sheaves on the site $(\text{Sch}/k)$ of all $k$-schemes, endowed with the fppf-topology.
Taking global sections yields an identification
$$
H^0(S, R^1f_*(G_Y)) = \Hom(\hat{G},\Pic^\tau_{Y/k}).
$$
The term on the left sits in the  exact sequence coming from the Leray--Serre spectral sequence
$$
0\lra H^1(S, f_*(G_Y) ) \lra H^1(Y, G_Y) \lra H^0(S, R^1f_*(G)) \lra H^2(S, f_*(G_Y)).
$$
If the ground field $k$ is algebraically closed,  every non-empty $k$-scheme of finite type has a rational point, by
Hilbert's Nullstellensatz, and it follows that $H^r(S,F)=0$ for every degree $r\geq 1$ and every abelian fppf-sheaf $F$.
Thus the outer terms in the preceding sequence vanish. This yields:

\begin{proposition}
\mylabel{raynaud correspondence}
If $k$ is algebraically closed, then the canonical map 
$$
H^1(Y, G_Y) \lra \Hom(\hat{G},\Pic^\tau_{Y/k})
$$
is bijective. In other words, the isomorphism classes of $G$-torsors $\epsilon:X\ra Y$ correspond to homomorphisms of
group scheme  $\hat{G}\ra\Pic_{Y/k}$.
\end{proposition}

Often there are \emph{canonical choices} for $G$ with respect to $Y$ and the resulting numerical
trivial part $P=\Pic^\tau_{Y/k}$.
For example, if $Y$ is Gorenstein   with numerically trivial dualizing sheaf $\omega_Y$,
we may take for $\hat{G}\subset P$ the discrete subgroup generated by the class of $\omega_Y$.
One may call the resulting $\epsilon:X\ra Y$ the \emph{$\omega$-canonical covering}.
If the Picard scheme is zero-dimensional, then one may choose $G$ with $\hat{G}=P$
or $\hat{G}=P^0$, 
and the resulting torsors are called the \emph{$\tau$-canonical covering} and \emph{$0$-canonical covering},
respectively.
For every $n\neq 0$ that is prime to the characteristic exponent $p\geq 1$,
we furthermore get the \emph{$n$-canonical coverings} with $\hat{G}=P[n]$,
the kernel of the multiplication by $n$-map. 
In characteristic $p>0$, we may also use the kernel of the relative Frobenius map $F:P\ra P^{(p)}$,
and get the \emph{$F$-canonical covering}.
Later, we will apply this to simply-connected Enriques surfaces $Y$. Then we have
$\Pic^0_{Y/k}=\Pic^\tau_{Y/k}=\Pic_{Y/k}[F]$,  and the ensuing torsor $\epsilon:X\ra Y$ is called
the \emph{K3-like covering}.

In what follows, we will assume that the ground field $k$ is algebraically closed,
such that we have  $H^1(Y, G_Y)=\Hom(\hat{G},\Pic^\tau_{Y/k})$.
This identification is  natural in $Y$ and $G$, and the following obvious consequence  will play a key role throughout:

\begin{lemma}
\mylabel{pullback trivial}
Let  $\epsilon:X\ra Y$ be  a $G$-torsor, corresponding to a homomorphism $\hat{G}\ra\Pic_{Y/k}$.
Let    $Z\ra Y$ be a proper morphism.
If $h^0(\O_Z)=1$, and the composite  homomorphism $\hat{G}\ra\Pic_{X/k}\ra\Pic_{Z/k}$ is zero,
then the induced torsor $X\times_Y Z\ra Z$ is trivial.
\end{lemma}

\proof
The condition $h^0(\O_Z)=1$ ensures that $H^1(Z,G)=\Hom(\hat{G},\Pic_{Z/k})$,
and the naturality of this identification implies that  the induced torsor  becomes trivial.
\qed

\medskip
As a consequence, we get a statement on the singularities of the total space of   torsors:

\begin{proposition}
\mylabel{singularity total space}
Let $\epsilon:X\ra Y$ be   a $G$-torsor with $G$ non-reduced, $Z\subset Y$ be a Weil divisor with
$h^0(\O_Z)=1$, and assume  that the composite map $\hat{G}\ra\Pic^\tau_{Y/k}\ra\Pic_{Z/k}$ is zero.
Let  $x\in X$ be a closed point mapping to a   point $z\in Z$ where the local ring $\O_{Z,z}$ is singular.
Then $\O_{X,x}$ is singular as well.
\end{proposition}

\proof
Set $n=\dim(\O_{X,x})$, such that the local ring $\O_{Z,z}$ on the Weil divisor has dimension $n-1$.
Since this local ring is not regular, its embedding dimension is $\edim(\O_{Z,z})\geq n$.
The induced torsor is  $T\simeq Z\times G$. The local ring  at the origin $0\in G$
has embedding dimension $\edim(\O_{G,0})\geq 1$, because the group scheme $G$ is non-reduced. Thus we have 
$$
\edim(\O_{X,x})\geq \edim(\O_{T,x})=\edim(\O_{Z,z})+\edim(\O_{G,0})\geq  n+1,
$$
hence the local ring $\O_{X,x}$ is not regular.
\qed

\medskip
Sometimes, one may verify the assumption by looking at cohomology groups:

\begin{corollary}
\mylabel{singularity via cohomology}
Assumptions as in the previous proposition. Suppose further   $\hat{G}$ is local, 
that $Y$ is normal  of dimension $n\geq 1$, and that the canonical mapping  
$H^{n-1}(Y,\omega_Y)\ra H^{n-1}(Y,\omega_Y(Z))$ is injective.
Let  $x\in X$ be a closed point mapping to a   point $z\in Z$ where $\O_{Z,z}$ is singular.
Then $\O_{X,x}$ is singular as well.
\end{corollary}

\proof
The short exact sequence $0\ra\O_Y(-Z) \ra\O_Y\ra\O_Z\ra 0$ gives a long exact sequence
$$
H^1(\O_X(-Z))\lra H^1(X,\O_X) \lra H^1(Z,\O_Z).
$$
The map on the left hand side is surjective, because it is  Serre dual to
the injection $H^{n-1}(Y,\omega_Y)\ra H^{n-1}(Y,\omega_Y(Z))$.
Whence the map on the right is zero. The latter is the tangent map 
for the restriction   $\Pic_{Y/k}\ra\Pic_{Z/k}$.
Suppose for a moment that $\hat{G}$ has height one. Then this group scheme is 
determined by its restricted Lie algebra
$\Lie(\hat{G})$,
the composite map $\hat{G}\ra\Pic_{Y/k}\ra\Pic_{Z/k}$ is trivial, and the assertion follows from the
proposition. In the general case, the Frobenius kernel $\hat{H}=\hat{G}[F]$  has height one
and corresponds to a surjection $G\ra H$. The result thus holds for the associated $H$-torsor $X\times_GH$,
and therefore also on $X$.
\qed

\medskip
The following connection to singularity theory follows immediately
from Theorem \ref{homological characterization}, our homological characterization of Zariski singularities,
or a direct local computation:

\begin{proposition}
\mylabel{canonical covering zariski}
Let $\epsilon:X\ra Y$ be a $G$-torsors. Suppose that $Y$ is smooth, and that $\hat{G}$ is unipotent of order $p$,
in other words either  $\hat{G}=\ZZ/p\ZZ$ or $\hat{G}=\alpha_p$. Then for each closed point
$x\in X$, the complete local ring $\O_{X,x}^\wedge$ is a Zariski singularity.
\end{proposition}

Now suppose that $Y$ is  a smooth surface.
A curve $E\subset Y$ is called \emph{negative-definite} if the intersection
matrix $N=(E_i\cdot E_j)$ attached to the integral components $E_1,\ldots,E_r\subset E$
is negative-definite.
This are precisely the \emph{exceptional divisors} for some contraction $Y\ra Y'$, 
where $Y'$ is  a normal 2-dimensional algebraic space.
To each such negative-definite curve we have the \emph{fundamental cycle} $Z=\sum n_iE_i$
with certain coefficients $n_i\geq 1$. It can be defined as the smallest such cycle
for which all intersection numbers are $(Z\cdot E_i)\leq 0$, and we have Artin's Algorithm
\cite{Artin 1966}
to compute it from the intersection matrix $N$.
We call the singular locus $\Sing(Z)$ of the fundamental cycle
the \emph{fundamental-singular locus} of $E\subset Y$. Note that 
$$
\Sing(Z)  = \bigcup_{i\neq j} (E_i\cap E_j)\quad  \cup \quad\bigcup_{n_i>1} E_i.
$$
This   coincides with the singular locus of the reduced curve $E$ if the
fundamental cycle is reduced, but usually is much larger.

We define the \emph{fundamental-genus} of the negative-definite curve $E\subset Y$
as the natural number $h^1(\O_Z)\geq 0$.
A negative-definite curve $E\subset Y$ is called of  \emph{rational type} if
it has fundamental-genus $h^1(\O_Z)=0$. In other words, it contracts
to a rational singularity.  In the special case of rational double points,
the intersection matrix $N$ then 
corresponds to one of  the Dynkin diagrams $A_n$, $D_n$, $E_6$, $E_7$ or $E_8$,
and we call the curve an \emph{ADE-configuration}.

\begin{proposition}
\mylabel{singularity over fundamental-singular}
Let $Y$ be a smooth surface, and  $\epsilon:X\ra Y$ be a $G$-torsor for some non-reduced $G$.
Then the local ring $\O_{X,x}$ is singular for each point $x\in X$ mapping
to the fundamental-singular locus $\Sing(Z)\subset X$ of a negative-definite curve of rational type $A\subset Y$.
\end{proposition}

\proof
Since $G$ is non-reduced, there is a   subgroup scheme $H\subset G$ so that
$G'=G/H$ is either $\mu_p$ or $\alpha_p$.
Set $X'=X/H$. By construction, $X\ra X'$ is a $H$-torsor, and 
$X'\ra Y$ is a $G'$-torsor, in fact coming from the
composite $\hat{G'}\subset \hat{G}\ra\Pic_{Y/k}$.
According to \cite{EGA IVc}, Proposition 6.7.4 it thus suffices to treat the case $G=G'$.
Then the Cartier dual $\hat{G}$ is isomorphic to $\ZZ/p\ZZ$ or $\alpha_p$, respectively. The Picard scheme
of the fundamental cycle is of the form $\Pic_{Z/k}=\ZZ^{\oplus r}$, where $r\geq 1$ is the number
of irreducible components of $Z$. In turn, the homomorphism $\hat{G}\ra\Pic_{Z/k}$ must be trivial,
and the assertion follows from Proposition \ref{singularity total space}.
\qed

\medskip
We say that a curve $A$ is \emph{semi-normal}
at some closed point $a\in A$ if  
the irreducible components in $\Spec(\O_{A,a}^\wedge)$ are normal and meet
like coordinate axes in $\AA^n$.
For example, this holds if $\O_{A,a}^\wedge$ is isomorphic to $k[[x,y]]/(xy)$,
that is, if the curve has \emph{normal crossings} at $a\in A$.

\begin{proposition}
\mylabel{singularity over nc}
Let $Y$ be a normal surface, and $\epsilon:X\ra Y$ be a $G$-torsor for some non-reduced $G$.
Let $A\subset Y$ be a curve whose normalization is a disjoint union of copies of $\PP^1$.
Then the local ring $\O_{X,x}$ is singular for each point $x\in X$ 
that maps to a point in $a\in A$ where  $A$ is singular but semi-normal.
\end{proposition}

\proof
Again it suffices to treat the case that $G$ is either $\mu_p$ or $\alpha_p$.
Clearly, we may assume that each irreducible component of $A$ contains the point $a\in A$.
Let $A'\ra A$ be the morphism that is the normalization outside $a\in A$, and an isomorphism
over an open neighborhood of $a\in A$.  We have  $h^0(\O_{A'})=1$, because the curve $A'$ is connected
and reduced. Moreover, its   Picard scheme then sits in a short exact sequence
$$
0\lra\GG_m^{\oplus r}\lra\Pic_{A'/k}\lra\ZZ^{\oplus s}\lra 0
$$
for some integers $r,s\geq 0$. Since $\hat{G}$ is isomorphic to $\ZZ/p\ZZ$ or $\alpha_p$,
there is only the zero homomorphism $\hat{G}\ra\Pic_{A'/k}$, whence the torsor $X\times_YA'\ra A'$
is trivial. It follows that the torsor $X\times_YA\ra A$ is trivial on some open neighborhood of $a\in A$.
By Proposition \ref{singularity total space}, the local ring $\O_{X,x}$ is singular.
\qed

\begin{proposition}
\mylabel{singularity over line}
Suppose $Y$ is a normal Gorenstein surface, and $\epsilon:X\ra Y$ is a $G$-torsor for some non-reduced $G$.
Let $A\subset Y$ be a curve isomorphic to $\PP^1$. If  the invertible sheaf $\omega_X$ is nef,
then the preimage $\epsilon^{-1}(A)$ passes through some  point $x\in X$ where the local ring
$\O_{X,x}$ is not factorial, in particular not regular.
\end{proposition}

\proof
It suffices to treat the case that $G$ is either $\mu_p$ or $\alpha_p$. Since $\Pic_{A/k}$ is torsion free and reduced,
the torsor is trivial over $A$, whence the preimage $\epsilon^{-1}(A)$ is isomorphic to $A\times G$.
In turn, $C=\epsilon^{-1}(A)_\red$ is isomorphic to $\PP^1$.
Seeking a contradiction, we assume that this preimage lies in the   locus where $X$ is locally factorial.
The conormal sheaf   of $C$ in $\epsilon^{-1}(A)=A\times G$ obviously is $\shI/\shI^2=\O_A$, and this is also the conormal sheaf of $C\subset X$.
It follows that $C^2=0$.
The Adjunction Formula gives $-2=(K_X\cdot C)  + C^2=(K_X\cdot C)$, contradicting $(K_X\cdot C)\geq 0$.
\qed

\medskip
We finally turn to some global invariants. 
Let us formulate the following fact:

\begin{proposition}
\mylabel{properties covering}
Suppose that $Y$ is normal and proper, and that $P=\Pic^\tau_{Y/k}$ is a finite unipotent
group scheme of order $p$. Let $\epsilon:X\ra Y$ be the $G$-torsor with $\hat{G}=P$.
Then
$$
h^0(\O_X)=1,\quad \chi(\O_X)=p\chi(\O_Y)\quadand \omega_X=\epsilon^*(\omega_Y).
$$
\end{proposition}

\proof
The assumption on $P=\Pic^\tau_{Y/k}$ means that the group scheme $G$ of order $p$ is local.
In turn, $\epsilon:X\ra Y$ is a universal homeomorphism. Let $\eta\in Y$ be the generic point.
Then $\O_{X,\eta}$ is an Artin local ring, and has degree $p$ as an algebra over the function field
$k(Y)=\O_{Y,\eta}$. Moreover, the field extension $k(Y)\subset k(X)$ is purely inseparable,
hence its degree $d\geq 1$ is either $d=1$ or $d=p$.

Seeking a contradiction, we assume $h^0(\O_X)>1$. Since the ground field $k$ is algebraically closed, the
algebra $H^0(X,\O_X)$ and whence the scheme $X$ is non-reduced. 
Since $X$ contains no embedded associated points, the Artin local ring $\O_{X,\eta}$ is non-reduced.
We infer that $d=1$, such that $X_\red\ra Y$ is birational.
By Zariski's Main Theorem, the finite birational $X_\red\ra Y$ is an isomorphism, hence
the $G$-torsor $\epsilon:X\ra Y$ admits a section, contradiction.

The statement on the Euler characteristic is a special case of  \cite{Mumford 1970}, 
Theorem 2 on page 121. 
It remains to establish the assertion on the dualizing sheaves.
The relative dualizing sheaf for the finite dominant morphism $\epsilon:X\ra Y$
is given by $\omega_{X/Y}=\underline{\Hom}_{\O_Y}(\O_X,\O_Y)$.
In light of $\omega_X=\epsilon^*(\omega_Y)\otimes\omega_{X/Y}$, it suffices
to show that $\omega_{X/Y}\simeq\O_X$.
Since $\epsilon:X\ra Y$ is a $G$-torsor, the finite flat $\O_Y$-algebra $\O_X$ of degree $p$ locally admits
a $p$-basis consisting of a single element.
It follows that the sheaf of K\"ahler differentials $\Omega^1_{X/Y}$ is invertible.
According to \cite{Kiehl; Kunz 1965}, Satz 9 we have $\omega_{X/Y}=(\Omega^1_{X/Y})^{\otimes (1-p)}$, so it suffices
to check that $\Omega^1_{X/Y}\simeq\O_X$.

Recall from \cite{EGA IVd}, Section 16.3 that the sheaf of K\"ahler differentials can be defined by the exact sequence
$$
0\lra\Omega^1_{X/Y} \lra\O_{X^{(1)}}\lra\O_X\lra 0,
$$
where $X=X\times_XX$ is the diagonal inside $X\times_YX$, and $X^{(1)}$ is its first infinitesimal neighborhood.
Since the morphism $\can:X\times G\lra X\times_YX$ given by $(x,g)\mapsto (x,xg)$ is an isomorphism,
we get a commutative diagram
$$
\begin{CD}
X	@>>> 	X\times G\\
@V\Delta VV		@VV\can V\\
X\times_XX	@>>>	X\times_YX,
\end{CD}
$$
where the upper map is given by $x\mapsto (x,e)$. Using that the underlying scheme of $G$
is isomorphic to the spectrum of $k[t]/(t^p)$, we infer 
that the conormal sheaf of $X\subset X\times_YX$ is isomorphism to $\O_X$,
in other words $\Omega^1_{X/Y}\simeq\O_X$.
\qed

\begin{corollary}
\mylabel{invariants covering}
Assumptions as in the proposition.
Suppose furthermore that $Y$ is Gorenstein, of even dimension $n\geq 2$, with $\omega_Y$ numerically trivial,
and such that 
$$
\sum_{i=3}^n (-1)^i h^i(\O_Y)\geq 0\quadand \sum_{j=2}^{n-1}(-1)^jh^j(\O_X)\leq 0.
$$
Then the above inequalities are equalities, and we have $p=2$, $\omega_X=\O_X$, $h^1(\O_X)=0$ and $h^0(\O_X)=h^n(\O_X)=1$.
\end{corollary}

\proof
The invertible sheaf $\omega_X=\epsilon^*(\omega_Y)$ is numerically trivial, so Serre Duality gives $h^n(\O_X)\leq 1$.
If the Picard scheme is non-smooth, we have $\Pic^\tau_{Y/k}=\alpha_p$, thus $h^1(\O_Y)=1$ and
the obstruction group $H^2(\O_Y)$ must be non-zero, according to \cite{Mumford 1966}, 
Corollary on page 198. 
If the Picard scheme is smooth, we have $\Pic^\tau_{Y/k}=\ZZ/p\ZZ$, consequently $h^1(\O_Y)=0$.
In both case, the inequality $1-h^1(\O_Y)+h^2(\O_Y)\geq 1$ holds. This gives a chain
of inequalities
$$
p\leq p(1-h^1(\O_Y)+h^2(\O_Y))\leq p\chi(\O_Y) =\chi(\O_X) \leq 1-h^1(\O_X) + h^n(\O_X) \leq 2-h^1(\O_X).
$$
The only possibility is $p=2$, $h^1(\O_X)=0$ and $h^n(\O_X)=1$. The latter shows that $\omega_X$ has
a non-zero global section $s:\O_X\ra\omega_X$. This map is necessarily bijective, because $\omega_X$ is numerically trivial.
\qed

\section{Simply-connected Enriques surfaces and   K3-like coverings}
\mylabel{Simply-connected enriques}

Fix an algebraically closed ground field $k$ of characteristic $p=2$,
and let $Y$ be an \emph{Enriques surface}. By definition, this is a smooth proper connected surface
with 
$$
c_1(Y)=0\quadand b_2(Y)=10.
$$
This implies that $Y$ is minimal, and that the dualizing sheaf  $\omega_Y$ is numerically trivial.
They  form one of the four classes of surfaces with $c_1=0$, the
other being  the abelian, bielliptic and K3-surfaces, which respective Betti number
$b_2=6$, $b_2=2$ and $b_2=22$, according to the   classification of surfaces.

For Enriques surfaces, the numerically trivial part $P=\Pic^\tau_{Y/k}$   is a group scheme of order two, and its group
of rational points $P(k)$ is generated by $\omega_Y$.
The Enriques surfaces come in three types: \emph{Ordinary, classical  and supersingular},
which means that $P$ is isomorphic to the respective group schemes
$\mu_2$, $\ZZ/2\ZZ$ and $\alpha_2$. Note that there are several other
designations in the literature.

Let $G=\ul\Hom(P,\GG_m)$ be the Cartier dual, such that $\hat{G}=P$,
and denote by  $\epsilon:X\ra Y$ the resulting $G$-torsor, as discussed in Section \ref{Canonical coverings}.
If $P=\mu_2$ is diagonalizable, this is an \'etale covering,
and $X$ is a K3 surface endowed with an free involution, which can also be viewed
as the universal covering of $Y$.
If $P=\ZZ/2\ZZ,\alpha_2$ is unipotent, then the Cartier dual $G=\mu_2,\alpha_2$
is local. In this case we say that $Y$ is a \emph{simply-connected Enriques surface},
and the $G$-torsor $\epsilon:X\ra Y$, which is a universal homeomorphism, is called the \emph{K3-like covering}.
Let us say that a connected reduced surface $X$ is a \emph{K3-like covering} if it is isomorphic
to the K3-like covering of some Enriques surface $Y$. 
We record:

\begin{proposition}
\mylabel{enriques quotient}
Suppose that a simply-connected reduced surface $X$ with   numerically trivial $\omega_X$
admits a free action of a local group scheme $G$
of order two whose quotient $Y=X/G$ is smooth. Then $Y$ is an Enriques surface,
the quotient map $\epsilon:X\ra Y$ is the K3-like covering, and $\Pic^\tau_{Y/k}\simeq\hat{G}$.
\end{proposition}

\proof
First note that since $X$ is reduced, the $G$-torsor $\epsilon:X\ra Y$ is non-trivial,
hence the corresponding homomorphism $\hat{G}\ra \Pic^\tau_{Y/k}$ is non-zero.
This already ensures that $Y$ is not a K3 surface.
The scheme $X$ is Cohen--Macaulay and Gorenstein, because this holds for $Y$,
and $\epsilon:X\ra Y$ is a $G$-torsor. In particular, $\omega_X$ is invertible.
The proof for Proposition \ref{properties covering} shows that $\epsilon^*(\omega_Y)=\omega_X$,
so $\omega_Y$ is numerically trivial.
By the classification of surface, $Y$ is either Enriques, abelian or bielliptic.
Since $G$ is local, the map $\epsilon:X\ra Y$ is
a universal homeomorphism, consequently  $Y$ is simply-connected.
Therefore, $Y$ is neither abelian nor bielliptic, and the only remaining possibility
is that $Y$ is an Enriques surface.
Thus $P=\Pic^\tau_{Y/k}$ has order two. It follows that the non-zero homomorphism $\hat{G}\ra P$
is an isomorphism, and $\epsilon:X\ra Y$ must be the K3-like covering.
\qed

\medskip
From now on, we assume that $Y$ is a simply-connected Enriques surface. 
We then have the following well-known facts on the K3-like covering $X$:

\begin{proposition}
\mylabel{facts on X}
The scheme $X$ is integral, the complete local rings $\O_{X,x}^\wedge$ at the closed points $x\in X$
are reduced Zariski singularities, and we have
$$
\omega_X=\O_X,\quad h^0(\O_X)=h^2(\O_X) =1\quadand h^1(\O_X)=0.
$$
Moreover, the scheme $X$ is not smooth.
\end{proposition}

\proof
The first statement follows from Proposition \ref{canonical covering zariski}, 
the second from Corollary \ref{invariants covering}.
Suppose $X$ would be smooth. By the classification of surfaces, $X$ is either an abelian surface
or a K3 surface,  thus the second Betti number is either $b_2=22$ or $b_2=6$. 
However, by the topological invariance of Betti numbers, we have
$b_2(X)=b_2(Y)=10$, contradiction.
\qed

\medskip
An integral  curve $A\subset Y$ with $A^2=-2$ is called
a \emph{$(-2)$-curve}. By the Adjunction Formula, these are the curves on the Enriques surface isomorphic to $\PP^1$.
An integral curve $C\subset Y$ with $C^2=0$ and $\Pic^\tau_{C/k}=\GG_a$ 
is called a \emph{rational cuspidal curve}. These are the curves on $Y$ that are isomorphic
to
$$
\Spec(k[t^2,t^3])\cup\Spec(k[t^{-1}]).
$$
An integral curve $F\subset Y$ with $F^2=0$ and $\Pic^\tau_{F/k}=\GG_m$ is called
a \emph{rational nodal curve}. These are the curves on $Y$ that can be seen as
$$
\Spec(k[[s,t]]/(st))\cup\Spec(k[s^{-1},t^{-1}]/(s^{-1}t^{-1}-1)).
$$

\begin{proposition}
\mylabel{curves images singularities}
Each $(-2)$-curve $A\subset Y$, each rational cuspidal curve $C\subset Y$
and each rational nodal curve $F\subset Y$
passes through the image of $\Sing(X)$.
Each negative-definite curve $E\subset Y$ is an ADE-configuration,
and its fundamental-singular locus is contained in the image of $\Sing(X)$.
\end{proposition}

\proof
To see that $E\subset Y$ is an ADE-configuration, let $Z$ be its fundamental cycle,
and consider the exact sequence
$$
H^1(Y,\O_Y)\lra H^1(Z,\O_Z)\lra H^2(Y,\O_Y(-Z))\lra H^2(Y,\O_Y).
$$
The map on the right is Serre dual to $H^0(Y,\omega_Y)\ra H^0(Y,\omega_Y(Z))$.
Any non-zero global section of $\omega_Y(Z)$ would define a curve $Z'\subset Y$   numerically
equivalent to $Z\subset Y$. Since the latter is negative definite, the only possibility is $Z'=Z$,
whence $\omega_Y\simeq\O_Y$. In this case, we conclude that the inclusion $H^0(Y,\omega_Y)\subset  H^0(Y,\omega_Y(Z))$
is bijective, again because $Z\subset Y$ is negative-definite.
The upshot is that in the above exact sequence, the map on the right is injective,
whence the map on the left is surjective.

If $Y$ is ordinary  then $h^1(\O_Y)=0$, hence $h^1(\O_Z)=0$ and $E\subset Y$ is an ADE-configuration.
If $Y$ is supersingular, we have $h^1(\O_Y)=1$ and $\Pic^\tau_{Y/k}=\alpha_2$.
Seeking a contradiction, we assume that the surjection $H^1(Y,\O_Y)\ra H^1(Z,\O_Z)$ is non-zero, hence
bijective. Using that curves have smooth Picard groups, we infer that $\Pic^\tau_{E/k}=\GG_a$.
Each irreducible component $E_i\subset E$ must be a $(-2)$-curve, whence $E_i\simeq\PP^1$ is smooth.
It follows that there are two irreducible components $E_i\neq E_j$ with intersection number
$n=(E_i\cdot E_j)\geq 2$. Then  the intersection  matrix $\begin{pmatrix} -2&n\\n&-2\end{pmatrix}$ is
not negative-definite, contradicting that $E\subset Y$ is negative-definite.
Again we see that $h^1(\O_Z)=0$, and the curve $E\subset Y$ is an ADE-configuration.

The assertions on the curves $A,F,E$ and the image of $\Sing(X)$
follow  from Propositions \ref{singularity over line}, \ref{singularity over nc}.
and \ref{singularity over fundamental-singular}, respectively.
It remains to treat the rational cuspidal curve $C$, which has $\Pic^\tau_{C/k}=\GG_a$.
If the restriction map $\Pic^\tau_{Y/k}\ra\Pic^\tau_{C/k}$ is trivial,
we may apply Proposition \ref{singularity total space}.
Now suppose that the restriction map is non-zero. Then $\tilde{C}=\epsilon^{-1}(C)$
is a non-trivial $G$-torsor over $C$. It becomes trivial after pulling-back
along the normalization $\nu:\PP^1\ra C$, according to Proposition \ref{pullback trivial}.
We thus get a cartesian diagram
$$
\begin{CD}
\PP^1\oplus\O_{\PP^1}	@>>>	\tilde{C}\\
@VVV				@VVV\\
\PP^1			@>>>	C.
\end{CD}
$$
where the horizontal arrows are birational. Consequently, the Weil divisor $\tilde{C}\subset X$
is of the form $\tilde{C}=2B$, where $B=\tilde{C}_\red$ is an irreducible Weil divisor,
and the morphism $B\ra C$ is birational.
Since $C$ is the rational cuspidal curve, the morphism $B\ra C$ is either the normalization or the identity map.
To see that the letter does not happen, note that the $G$-torsor $\epsilon:X\ra Y$ becomes
trivial after pulling-back to itself. Since $h^0(\O_X)=1$, this means
that  the restriction map $\Pic^\tau_{Y/k}\ra\Pic^\tau_{X/k}$ is trivial. It follows that
the composite map $\Pic^\tau_{Y/k}\ra\Pic^\tau_{B}$ is trivial, contradiction.
Thus we must have $B=\PP^1$. Consider the short exact sequence
$$
0\lra\shN\lra\O_{\tilde{C}}\lra \O_B\lra 0.
$$
The sheaf $\shN$ is torsion free, and therefore invertible when regarded as module over $\O_B$.
By Riemann--Roch, its degree is 
$$
\deg(\shN) = \chi(\shN) - \chi(\O_B) = \chi(\O_{\tilde{C}}) - 2\chi(\O_B) = 2\chi(\O_C) -2\chi(\O_B) =-2.
$$
Seeking a contradiction, we now assume that the surface $X$ is locally factorial along $B$.
Then $B\subset X$ is a Cartier divisor isomorphic to $\PP^1$ with selfintersection number $B^2=-\deg(\shN) = 2$.
Riemann--Roch yields $-2=\deg(K_B)=(K_X+B)\cdot B = 2$, contradiction.
\qed

\medskip
As the referee noticed, an alternative argument may use a genus-one fibration $\varphi:Y\ra\PP^1$
and  Katsura's result (\cite{Katsura 1982}, Proposition 3.2) that for some  non-zero rational 1-form $\omega$ on the projective line,
in suitable coordinates actually given by   $\omega=dt/t$ or $\omega=dt$, the pull-back 
$\varphi^*(\omega_Y)$ to the Enriques surface extends to a global section of $\Omega^1_{Y/k}$. The   zero-locus
of this global 1-form on $Y$ describes the the image of the singular locus of the K3-like covering $X$, and from
this the assertion can be inferred.
In some sense, the above proof is a coordinate-free version for this.

We now furthermore assume that the K3-like covering $X$ has only isolated singularities, that
is, the surface $X$ is normal.
Let $r:S\ra X$ be a resolution of singularities. 
The \emph{geometric genus} of a singularity $x\in X$ is the length 
$p_g(\O_{X,x}) = \length R^1r_*(\O_S)_x$.
A singularity $x\in X$ is called \emph{rational} if its geometric genus is $p_g=0$,
and  we call it  \emph{elliptic} if $p_g=1$. Note that this term has various meanings
in the literature. For example, Wagreich \cite{Wagreich 1970} uses it for singularities
with arithmetic genus $p_a=1$, which is the largest integer of the form $1-\chi(\O_Z)$,
where $Z\subset S$ ranges over the exceptional divisors.
According to loc.\ cit., this includes all singularities with fundamental genus $p_f=1$,
which is defined as $p_f(\O_{X,x})=h^1(\O_Z)$, where $Z\subset S$ is the fundamental divisor
\cite{Artin 1966}. A related class of singularities are the minimally elliptic singularities
studied  by Laufer \cite{Laufer 1977}.
A singularity with Hilbert--Samuel multiplicity
$e(\O_{X,x})=2$ is called a \emph{double point}.

\begin{proposition}
\mylabel{singularities on X}
All rational singularities on the K3-like covering $X$ are rational double points,
and the only possible types are $A_1$,  $E_7^0$, $E_8^0$  and $D_{2n}^0$.
There is at most one non-rational singularity $x\in X$. If present, it is an elliptic double point.
\end{proposition}

\proof
The singularities on the  K3-like covering are Zariski singularities, whence
have multiplicity $e=2$, by Proposition \ref{hilbert-samuel multiplicity}. Furthermore, they 
are Gorenstein, thus all rational singularities must be rational double points.
The only rational double points that are Zariski singularities are of type $A_1, D_{2n}^0, E_7^0,E_8^0$,
according to Proposition \ref{zariski rdp}.
Next, consider the Leray--Serre spectral sequence for the minimal resolution of singularities $r:S\ra X$.
It gives an exact sequence
\begin{equation}
\label{long exact sequence}
0\lra H^1(X,\O_X)\lra H^1(S,\O_S)\lra H^0(X,R^1r_*\O_S)\lra H^2(X,\O_X).
\end{equation}
The term on the left vanishes, and  the term on the right is one-dimensional. Moreover, we have
$K_S=K_{S/X}$, which is a Cartier divisor supported on the exceptional curve with coefficients $\leq 0$.

Now suppose that there is a non-rational singularity $x\in X$. Then $K_{S/X}< 0$,
and Serre Duality gives $H^2(S,\O_S)=0$. Thus the Picard scheme is smooth, and the connected component
$\Pic^0_{S/k}$ is an abelian variety.
On the other hand,  each  integral curve $E\subset S$  contained in the exceptional locus
is a curve of genus zero, 
so the restriction map $\Pic^0_{S/k}\ra\Pic^0_{E/k}$ is trivial. In turn, the Albanese morphism
$S\ra A$ contracts each such $E\subset S$, and thus factors over $Y$. Using that $\Pic_{Y/k}$
is 0-dimensional, 
we conclude that $\dim(A)=0$, and thus $H^1(S,\O_S)=0$. It now follows from  the above exact sequence
that $R^1r_*(\O_S)$ has length $\leq 1$. In turn, there is precisely one such singularity, and its geometric genus is one.
\qed

\medskip
We will later see that if there is an elliptic singularity $x\in X$, there are no further
singularities. In any case, the structure of the smooth surface $S$ depends 
on the nature of the singularities on $X$:

\begin{proposition}
\mylabel{numbers for S}
The smooth surface $S$ is a   K3 surface with $\rho(S)=b_2(S)=22$
if the K3-like covering $X$ contains only rational singularities and the resolution $r:S\ra X$ is minimal,
and is a rational surface with $\rho(S)=b_2(S)\geq 11$ otherwise.
\end{proposition}

\proof
If all singularities are rational, the exact sequence \eqref{long exact sequence}
shows that $h^1(\O_S)=0$. Furthermore, we have $\omega_{S/X}=\O_S$, whence $\omega_S=\O_S$.
By the classification of algebraic surfaces, $S$ is a K3 surface.
It then has Betti number $b_2=22$, whereas our Enriques surface $Y$ and its
K3-like covering have $b_2(X)=b_2(Y)=10$. It follows that
the exceptional divisor for $r:S\ra X$ consist of $12=b_2(S)-b_2(X)$ irreducible components.
It follows that $\rho(S)=b_2(S)$.

If there is a non-rational singularity $x\in X$, then $-K_S=-K_{S/X}$ is an effective Cartier divisor,
and we still have $H^1(S,\O_S)=0$. By the classification of surfaces, $S$ is rational.
In any case, the Picard number satisfies $\rho(S)>\rho(X)\geq \rho(Y)=10$, because there is at least one singularity,
and thus $b_2(S)\geq 11$.
\qed

\medskip
Every Enriques surface $Y$ admits at least one  genus-one fibration $\varphi:Y\ra \PP^1$.
It has one or two   multiple fibers $Y_b=2C$, and their multiplicity is $m=2$.
Two multiple fibers occur  precisely for the classical Enriques surfaces, and we then have
$\omega_Y=\O_Y(C_1-C_2)$, where $C_1,C_2\subset Y$ are the two half-fibers.
Two fibrations $\varphi,\varphi':Y\ra\PP^1$ are called \emph{orthogonal}
if for   respective half-fibers $C,C'\subset Y$ have intersection number  $(C\cdot C')=1$.

For our simply-connected Enriques surface $Y$ whose K3-like covering $X$ is normal,
there are some strong restrictions on the nature of these fibrations.
The following were already observed by Cossec and Dolgachev
(see \cite{Cossec; Dolgachev 1989}, Proposition 5.7.3 and its Corollary):

\begin{theorem}
\mylabel{properties genus-one fibrations}
For every genus-one fibration $\varphi:Y\ra \PP^1$, the following holds:
\begin{enumerate}
\item
The fibration $\varphi:Y\ra \PP^1$ is   elliptic, and not  quasielliptic.
\item
Each singular fiber   $Y_b=\varphi^{-1}(b)$ is of Kodaira type  $\I_n$ for some $1\leq n\leq 9$,
or of Kodaira type $\II$, $\III$ or $\IV$.
\item
There are at least two other elliptic fibrations $\varphi',\varphi'':Y\ra\PP^1$ so that
$\varphi,\varphi',\varphi''$ are mutually orthogonal.
\end{enumerate}
\end{theorem}

\proof
Suppose there is some quasielliptic fibration $\varphi:Y\ra\PP^1$. Then almost all fibers
$Y_a$, $a\in\PP^1$ are rational cuspidal curves. According to Proposition \ref{curves images singularities},
the normal surface $X$ must contain infinitely many singularities, contradiction.

The Kodaira types in (ii) correspond to the curves of canonical type that are reduced.
If there would be another Kodaira type, the corresponding fiber contains 
an ADE-configuration of type $D_4$. Its fundamental cycle is non-reduced, and again
$X$ contains a curve of singularities, contradiction.

The last property holds for all Enriques surfaces with a few exceptions, according
to \cite{Cossec; Dolgachev 1989}, Theorem 3.5.1. The exceptions are called \emph{extra special} in loc.\ cit.,
and each of them contains a genus-one fibration containing a fiber of type $\II^*$,
$\I_4^*$ or $\III^*$, which are non-reduced. This contradicts  Property (ii).
\qed

\medskip
An integral curve $E\subset Y$ is called \emph{non-movable} if $h^0(\O_Y(E))=1$.

\begin{proposition}
\mylabel{preimage elliptic}
Let $E\subset Y$ be a non-movable elliptic curve.
Then the schematic preimage $\epsilon^{-1}(E)\subset X$ is an elliptic curve, and
the K3-like covering $X$ is smooth along this preimage.
\end{proposition}

\proof
The curve $2E$ is movable, and defines an elliptic fibration $\varphi:Y\ra\PP^1$ for which
$E\subset Y_b$ is a half-fiber. Suppose first that $Y$ is classical.
Then there is another half-fiber $C$, and  $\omega_Y=\O_Y(E-C)$ has order two in $\Pic(Y)$.
According to \cite{Cossec; Dolgachev 1989}, Theorem 5.7.2 multiple fibers are not wild,
so the restriction $\omega_Y|E=\O_E(E)$ has order two in $\Pic(E)$.
It follows that the restriction map $\Pic^\tau_{Y/k}\ra\Pic^\tau_{E/k}$ is injective.
Now suppose that $Y$ is supersingular, such that $\omega_Y=\O_Y$. The short exact sequence
$0\ra\O_Y(-E)\ra \O_Y\ra\O_E\ra 0$ yields a long exact sequence
$$
H^1(Y,\O_Y)\lra H^1(E,\O_E)\lra H^2(Y,\O_Y(-E))\lra H^2(Y,\O_Y).
$$
The map on the right is Serre dual to $H^0(Y,\O_Y)\ra H^0(Y,\O_Y(E))$.
The latter is bijective, because the curve $E\subset Y$ is non-movable.
In turn, the map on the left $H^1(Y,\O_Y)\ra H^1(E,\O_E)$ is surjective,
and again we conclude that the restriction map $\Pic^\tau_{Y/k}\ra\Pic^\tau_{E/k}$ is injective.

In both cases, we see that the $G$-torsor $E'=\epsilon^{-1}(E)\ra E$ is non-trivial.
According to Proposition 4.9, the curve $E'$ is integral, with $h^0(\O_{E'})=h^1(\O_{E'})=1$.
This curve is normal. If not, the normalization $E''\ra E'$ must have $h^1(\O_{E''})<1$,
giving an integral surjection $\PP^1\ra E$, contradiction.
Thus $E'$ is a smooth Cartier divisor in the surface $X$, and we conclude that $X$ is smooth
along $E'$.
\qed

\medskip
Let us write $x_1,\ldots,x_m\in X$ for the singularities on the K3-like coverings.
The complete local rings $\O_{X,x_i}$ are
Zariski singularity, thus given by a  formal power series of the form  $g_i=z^2-f_i(x,y)$.
Recall that the \emph{Tjurina number} satisfies 
$$
\tau_i = 2\cdot\length k[[x,y]]/(f_i,\partial f_i/\partial x,\partial f_i/\partial y),
$$
by Proposition \ref{tjurina for isolated zariski}.
Note that this is always an even number.  

\begin{proposition}
\mylabel{theta and tjurina}
We have  $\Theta_{X/k}=\O_X^{\oplus 2}$, 
and the Tjurina numbers for the singularities satisfy $\sum_{i=1}^m \tau_i = 24$.
In particular, if $S$ is a K3 surface, then the exceptional divisor for
the resolution of singularities $r:S\ra X$ has 12 irreducible components.
\end{proposition}

\proof
If $S$ is K3, this is due to Ekedahl,   Hyland and Shepherd-Barron \cite{Ekedahl; Hyland; Shepherd-Barron 2012}, 
and their arguments generalize as follows:
Since $\epsilon:X\ra Y$ is a torsor for the Cartier dual $G=\underline{\Hom}(\Pic^0_{Y/k},\GG_m)$,
we have a short exact sequence
$$
0\lra\O_X\lra\Theta_{X/k}\lra \shL\lra 0
$$
for some invertible sheaf $\shL$. The map on the left  corresponds to 
the $p$-closed vector field defining the $G$-action on $X$.
Using that $\det(\Theta_{X/k})=\omega_X=\O_X$,
we conclude that $\shL\simeq\O_X$. The preceding extension splits, because
we have $\Ext^1(\O_X,\O_X)=H^1(X,\O_X)=0$. Summing up, $\Theta_{X/k}=\O_X^{\oplus 2}$.

Now suppose that $S$ is a K3 surface.
It then has Betti number $b_2=22$, whereas our Enriques surface $Y$ and its
K3-like covering have $b_2(X)=b_2(Y)=10$. It follows that
the exceptional divisor for $r:S\ra X$ consist of $12=b_2(S)-b_2(X)$ irreducible components.
Since the resulting singularities are rational double points, we must have $\sum \tau_i = 24$.

Finally suppose that $S$ is rational.
Here we have Betti number $b_2(S)=10-K_S^2$
and Chern number $c_2(S)=2 + b_2(S) = 12-K_S^2$. The Chern numbers  are related,
according to \cite{Ekedahl; Hyland; Shepherd-Barron 2012}, Proposition 3.12 and Corollary 3.13,
in the following way:
$$
c_2(S) = c_2(X) + \gamma \quadand c_2(\Theta_{X/k}) = c_2(X) - \tau.
$$
Here    $c_2(X)=c_2(L_{X/k}^\bullet)$ is defined as the second Chern class of the cotangent complex
$L_{X/k}^\bullet$, which in our situation is a  complex of length one comprising locally free sheaves
of finite rank.
It can be defined with the help of an embedding $X\subset\PP^n$ into some projective $n$-space.
Moreover, $\tau=\sum \tau_i$, and the other correction term $\gamma$ is given by the Local Noether Formula
$$
-1=-\length R^1g_*(\O_S) = \frac{K_S^2 + \gamma}{12},
$$
from \cite{Ekedahl; Hyland; Shepherd-Barron 2012}, Proposition 3.12. 
In other words $\gamma= -12 - K_S^2$. Combining these equations, one gets  
$$
0=c_2(\Theta_{X/k}) = c_2(S) -\gamma -\tau = (12-K_S^2) + (12 +K_S^2) -\tau = 24-\tau.
$$
Again, we have $\sum \tau_i=24$.
\qed

\section{Normal surfaces with trivial tangent sheaf}
\mylabel{Trivial tangent sheaf}

Let $k$ be an  algebraically closed ground field   of characteristic $p=2$.
In order to understand and construct simply-connected Enriques surface $Y$ and
their K3-like coverings $X$, we now  impose  
eight axiomatic conditions on  certain elliptic fibrations.

Let $h:S\ra\PP^1$ be a smooth elliptic surface,  whose total space $S$ is either
a   rational surface  or a K3 surface.
An   curve $C\subset S$ is called \emph{vertical} if each irreducible component  is contained
in some fiber $h^{-1}(a)$, $a\in\PP^1$. It is    \emph{horizontal} if each irreducible component
dominates $\PP^1$. The same locution is used for fibrations with normal total space.
Let $S\ra S'\ra\PP^1$ be the relative minimal model, obtained by successively contracting
vertical $(-1)$-curves. Note that if $S$ is a K3 surface, there are no $(-1)$-curves
at all, such that $S=S'$.

Let $E\subset S$ be a  vertical negative-definite curve, and  $r:S\ra X$ be its contraction.
Then $X$ is a   proper normal   surface, with an induced elliptic fibration $f:X\ra\PP^1$, satisfying $h=f\circ r$.
We denote by $x_1,\ldots,x_r\in X$ the images of the connected components $E_1,\ldots,E_r\subset E$,
and assume that each local ring $\O_{X,x_i}$ is singular.
Then $r:S\ra X$ is a resolution of singularities, but we do not assume that it is the minimal resolution.
Our     conditions are:

\newcounter{ENr}
\setcounter{ENr}{-1}
\newcommand{\En}[1]{\refstepcounter{ENr}\textbf{(E\arabic{ENr})} \label{#1} }
\newcommand{\eref}[1] {\text{\rm (E\ref{#1})}}

\medskip
\begin{list}{-mm}{\leftmargin2em\itemsep1em}
\item[\En{Reduced}] 
\emph{The    fibers $f ^{-1}(a)\subset X$ are reduced for all closed points $a\in\PP^1$.}

\item[\En{CI}] 
\emph{The singular local rings $\O_{X,x_i}$, $1\leq i\leq r$ are complete intersections
with free tangent modules.}

\item[\En{Tjurina}] 
\emph{The Tjurina numbers $\tau_i$ for the local rings $\O_{X,x_i}$ add up to
$\sum\tau_i =24$.}

\item[\En{Length}] 
\emph{If  $S$ is a rational surface, the skyscraper sheaf $R^1r_*(\O_S)$ has length one;
it is supported in the multiple fiber  if $h:S\ra\PP^1$ has a multiple fiber.}

\item[\En{Nonfixed}] 
\emph{If $S$ is a rational surface, then $\omega_S=\shN(D)$, where
$D$ is a divisor supported on the exceptional divisor $E\subset S$,
and the invertible sheaf $\shN$ is non-fixed.}

\item[\En{Elliptic}] 
\emph{There is a horizontal Cartier divisor $F\subset X$ that is an elliptic curve.}

\item[\En{Zariski}] 
\emph{For each singularity $x_i\in X$, the local rings $\O_{X,x_i}$ are Zariski singularities.}

\item[\En{Cuspidal}] 
\emph{There is a horizontal Cartier divisor $C\subset X$ that is a rational cuspidal curve,
and the singular point on $C$ is a regular point on $X$.}
\end{list}

\medskip
Here the  invertible sheaf $\shN\in\Pic(S)$ is called \emph{non-fixed}
if $h^0(\shN)\geq 1$, and  furthermore $h^0(\shN)=1$ implies $\shN=\O_S$.
Note that \eref{CI} is a consequence of \eref{Zariski}, according to  Proposition \ref{theta for zariski surface}.
However, it is useful to have \eref{CI} as a separate condition, because it will be handy when 
the Enriques surface acquires rational double points.
Our conditions indeed help to analyze and construct simply-connected Enriques surface,
because we have:

\begin{proposition}
\mylabel{conditions necessary}
Suppose $Y$ is a simply-connected Enriques surface whose K3-like covering $X$ is normal,
$r:S\ra X$ is a resolution of singularities, $E\subset S$ is the exceptional divisor, 
and $h:S\ra \PP^1$ is the fibration
induced from some elliptic fibration $\varphi:Y\ra \PP^1$.
Then conditions \eref{CI}--\eref{Zariski} do hold, whereas condition \eref{Cuspidal} does not hold.
\end{proposition}

\proof
Condition \eref{Nonfixed} is true, because $\omega_X=\O_X$, and thus $\omega_S=\O_S(K_{S/X})$,
where the relative canonical divisor $K_{S/X}$ is supported by the exceptional divisor.
Property \eref{Zariski} and thus also \eref{CI} hold  according to Theorem \ref{homological characterization}.
Proposition \ref{theta and tjurina} gives \eref{Tjurina}, and Theorem \ref{properties genus-one fibrations} ensures
\eref{Elliptic}.

Condition \eref{Length} pertains to the case that $S$ is rational.
According to  Proposition \ref{numbers for S} together with Proposition 
\ref{singularities on X}, there is exactly one elliptic singularity $x_1\in X$,
whence $R^1r_*(\O_S)$ has length one.
Now suppose that the relatively minimal rational elliptic surface $S'\ra\PP^1$ has a multiple fiber
$S'_a$. Write $C'$  for the underlying indecomposable curve of canonical type.
Then $K_{S'}=-C'$, according to \cite{Cossec; Dolgachev 1989}, Proposition 5.6.1.
Seeking a contradiction, we assume that the elliptic singularity $x_1\in X$ does not lie on $X_a$.
Without restriction, we may assume that $r:S\ra X$ is the minimal resolution of singularities.
Then $K_S=-E$, where $E\subset S$ is a negative-definite curve mapping to $x_1\in X$.
In turn, $K_{S'}= -E'$, where $E'\subset S'$ is the image of $E$.
Thus the curves $E'\neq C'$ are linearly equivalent. But since $S'_a$ is a multiple fiber,
the indecomposable curve of canonical type  $C'$ is non-movable, contradiction.

Finally, we have to verify that condition \eref{Cuspidal} does not hold.
Suppose to the contrary that there is a Cartier divisor $C\subset X$
that is a rational cuspidal curve such that the local ring $\O_{X,c}$ is regular,
where  $c\in C$ is the singular point.
The Adjunction Formula $\deg(\omega_C) = (K_X+C)\cdot C$ implies $C^2=0$.
Consequently,    the image $D=\epsilon(C)$ on the Enriques surface 
is an integral rational curve with $D^2=0$. It follows that $D$ is   a rational cuspidal curve.
By Proposition \ref{curves images singularities}, there is
a singular point $x\in X$ mapping to $D$. 
But for each closed point $x\in C$, the local ring $\O_{X,x}$ is regular, because
either $x=c$ or $\O_{C,x}$ is regular, contradiction.
\qed

\medskip
The goal of this section is to establish a converse for Proposition \ref{conditions necessary}.
Let us write $\O_X(n)$ for the preimage of the invertible sheaves $\O_{\PP^1}(n)$ under
the elliptic fibration $f:X\ra\PP^1$, and likewise we write $\O_S(n)$ and $\O_{S'}(n)$.
The following  is the key step in producing Enriques surfaces with normal K3-like coverings:

\begin{theorem}
\mylabel{conditions for theta}
Suppose that  $r:S\ra X$ satisfies the conditions  \eref{CI}--\eref{Nonfixed}. Then the following holds:
\begin{enumerate}
\item 
The dualizing sheaf is given by $\omega_X=\O_X$.
\item
The tangent sheaf $\Theta_{X/k}$ is locally free of rank two,
and there is an extension $0\ra \shF^\vee\ra\Theta_{X/k}\ra\shF\ra 0$ for
some coherent subsheaf $\shF$ inside the invertible sheaf $f^*(\Theta_{\PP^1/k})=\O_X(2)$. 
\item
If  also  condition  \eref{Reduced} holds, the inclusion $\shF\subset\O_X(2)$ is an equality,
and we have $\Theta_{X/k}=\O_X(-n)\oplus\O_X(n)$ for some $0\leq n\leq 2$. 
\item
We have $n=0$ if furthermore  \eref{Elliptic} holds,  wheras $n\neq 0$
if we have \eref{Cuspidal} instead.
\end{enumerate}
%
\end{theorem}

\proof
Our first step is to show that $\omega_X=\O_X$.
The case that $S$ is K3 is obvious: Every negative-definite curve on $S$  
produces a rational double point on $X$,  thus our $X$ is Gorenstein with $\omega_X=\O_X$.
The more interesting case is that $S$ is rational: Then $H^1(S,\O_S)=H^2(S,\O_S)=0$,
and the Leray--Serre spectral sequence for 
the resolution of singularities $r:S\ra X$ gives
an exact sequence
$$
0\ra H^1(X,\O_X)\ra H^1(S,\O_S)\ra H^0(X,R^1r_*(\O_S))\ra H^2(X,\O_X)\ra H^2(S,\O_S).
$$
It follows that $h^1(\O_X)=0$ and   $h^2(\O_X)=1$, the latter by condition \eref{Length}.
Serre duality gives $h^0(\omega_X)=1$. 
Suppose that $\omega_X$ is non-trivial. Then $\omega_X=\O_X(C)$ for some Cartier divisor
$C\subset X$ that is not linearly equivalent to any other effective Cartier divisor.
According to condition \eref{Nonfixed}, we have $\omega_S=\shN(D)$, where
$D\subset S$ is supported by the exceptional divisor $E\subset S$ and $\shN$ is non-fixed.
If $\shN=\O_S$, then $\omega_X$ is trivial outside the singularities,
whence everywhere trivial, contradiction. Thus there are  two curves $A\neq B$ on $S$
both giving $\shN$. Write $A=A_1+A_2$ and $B=B_1+B_2$ where $A_2$, $B_2$ are the parts
supported by the exceptional divisor $E\subset S$. Using $h^0(\omega_X)=1$,
we infer that both $A_1,B_1\subset S$ map to $C\subset X$, hence $A_1=B_1$. It follows that 
$A_2\neq B_2$ are linearly equivalent and supported by the exceptional divisor $E\subset S$. This contradicts the
fact that the intersection form $(E_i\cdot E_j)$ is negative-definite.
Summing up, we have $\omega_X=\O_X$.

The tangent sheaf $\Theta_{X/k}$ is locally free of rank two, according to condition \eref{CI}.
Note that the dual of $\Theta_{X/k}$ coincides with the bidual of $\Omega^1_{X/k}$. Consequently
we have $\det(\Theta_{X/k})=\omega_X^\vee=\O_X$.
The next step is to verify that the second Chern class $c_2(\Theta_{X/k})$ vanishes.
The case of K3 surfaces was already treated in \cite{Ekedahl; Hyland; Shepherd-Barron 2012}, Section 3.
We proceed in a similar way: In both cases on has
$$
c_2(\Theta_{X/k}) = c_2(X) - \tau\quadand c_2(X) = c_2(S) - \nu,
$$
where $\tau=\sum\tau_i=24$ is the sum of Tjurina numbers, and $\nu$ is a correction term
$$
\nu =   c_1^2(X) - c_1^2(S) -12\cdot\length R^1r_*(\O_S).
$$
according to the Local Noether Formula
\cite{Ekedahl; Hyland; Shepherd-Barron 2012}, Proposition 3.12. In our situation $c_1^2(X)=c_1^2(\Theta_X)=0$.
If $S$ is K3, then $c_1^2(S) =0$, $c_2(S)=24$ and $\nu=0$, which gives $c_2(\Theta_{X/k})=0$.
If $S$ is rational, we have $c_2(S) = 2 + b_2=12-c_1^2(S)$ and $\nu= c_1^2(S) -12$,
also giving 
$$
c_2(\Theta_{X/k})= c_2(S) -\tau -\nu = (12 - c_1^2(S)) - 24 - (-c_1^2(S) -12) = 0.
$$

Our next step is to regard the tangent sheaf as an extension.
Consider the short exact sequence
$f^*(\Omega^1_{\PP^1/k})\ra \Omega^1_{X/k} \ra\Omega^1_{X/\PP^1}\ra 0$.
The map on the left is injective, because the function field extension $k(\PP^1)\subset k(X)$ is separable
and $X$ has no embedded components. Thus we have a short exact sequence
$$
0\lra f^*(\Omega^1_{\PP^1/k})\lra \Omega^1_{X/k} \lra\Omega^1_{X/\PP^1}\lra 0.
$$
The sheaf $\Omega^1_{X/\PP^1}$ is invertible at each point
$x\in X$ at which the fiber  $X_{f(x)}$ is regular. 
Dualizing, we get an exact sequence $0\ra (\Omega^1_{X/\PP^1})^\vee\ra\Theta_{Y/k}\ra\shF\ra 0$,
for some coherent subsheaf $\shF$ of $f^*(\Theta_{\PP^1/k})=\O_X(2)$.
Both outer terms are invertible in codimenson one, 
and the term on the left is reflexive. From $\det(\Theta_{Y/k})=\O_X$ we infer that
$(\Omega^1_{X/\PP^1})^\vee=\shF^\vee$, and get the desired extension
\begin{equation}
\label{extension sequence}
0\lra\shF^\vee\lra \Theta_{X/k}\lra \shF\lra 0.
\end{equation}
This establishes  assertions (i) and (ii).

Now assume that condition \eref{Reduced} holds, which means that all  closed fibers $f^{-1}(a)$ 
are reduced.
Then $\Omega^1_{X/\PP^1}$ is invertible in codimension one,
hence the inclusion $\shF\subset\O_X(2)$ is an isomorphism
outside finitely many closed points.
It follows that the dual map  $\O_X(-2)\ra\shF^\vee$ is an isomorphism
and that $\shF=\shI\O_X(2)$ for some   sheaf of ideals
$\shI\subset\O_X$ that defines  a finite subscheme $Z\subset X$.
Note that the structure sheaf $\O_Z$ has finite projective dimension,
by the exact sequence \eqref{extension sequence}.
Rewrite this sequence as  $0\ra\shL^\vee\ra\Theta_{Y/k}\ra \shI\shL\ra 0$,
with $\shL=\O_X(2)$.
Then
$$
0=c_2(\Theta_{X/k}^\vee) = -c_1^2(\shL) + \length(\O_X/\shI).
$$
Since $\shL$ comes from $\PP^1$, we have $c_1^2(\shL)=0$ and whence $\shI=\O_X$. 
Summing up, $\shF=\shL$, and the tangent sheaf is an extension of $\shL=\O_X(2)$ by $\shL^\vee=\O_X(-2)$. The extension class
lies in 
$$
\Ext^1(\O_X(2),\O_X(-2)) = H^1(X,\O_X(-4)).
$$
The latter can be computed with the  Leray--Serre spectral sequence, together with the Projection Formula:
We get an exact sequence
$$
0\lra H^1(\PP^1,\O_{\PP^1}(-4))\lra H^1(X,\O_X(-4)) \lra H^0(\PP^1,  R^1f_*(\O_X)(-4)).
$$
The term on the right vanishes: The higher direct image sheaf $R^1f_*(\O_X)$ commutes with base-change.
Since $S\ra\PP^1$ has no wild fibers, according to \cite{Cossec; Dolgachev 1989}, Proposition 5.6.1,
it follows that the map $t\mapsto h^1(\O_{X_t})$
takes the constant values one, and thus $R^1f_*(\O_X)$ is invertible.
The Leray--Serre spectral sequence yields
$$
\chi(R^1f_*(\O_X)) = \chi(\O_{\PP^1}) - \chi(\O_X) = 1 - 2 =-1,
$$
thus $R^1f_*(\O_X)=\O_{\PP^1}(-2)$, and consequently $R^1f_*(\O_X)(-4)$ has no non-zero global
section.
We infer that  the extension of $\O_X(2)$ by $\O_X(-2)$ giving $\Theta_{X/k}^\vee$ comes
from an extension of $\O_{\PP^1}(-2)$ by $\O_{\PP^1}(2)$,
in particular $\Theta_{X/k}$ is the preimage of some locally free sheaf $\shE$ of rank two
on $\PP^1$. But all such sheafs on $\PP^1$ are sums of line bundles, according to
Grothendieck (see for example \cite{Okonek; Schneider; Spindler 1980}). 
The splitting type of $\shE$ is of the form $(-n,n)$ for some
integer $n\geq 0$, because  $c_1(\shE)=0$. Since $\shE$ surjects onto $\O_{\PP^1}(2)$,
only the three possibilities $n=0,1,2$ exists.
This gives assertion (iii).

Now suppose that furthermore condition \eref{Elliptic} holds:
Consider the horizontal Cartier divisor $F\subset X$ that is an elliptic curve,
and let $d>0$ be the degree of the map $F\ra\PP^1$.
We have  an exact sequence
$$
0\lra \O_F(-F)\lra \Omega^1_{X/k}|_F\lra\Omega^1_{F/k}\lra 0.
$$
The scheme $X$ is regular along the regular Cartier divisor $F$, so the above
sheaves are locally free on $F$.
The sheaves $\Omega^1_{X/k}$ and $\Theta_{X/k}^\vee$ have the
same restriction to $F$, which is thus of form $\shM\oplus\shM^\vee$,
where $\deg(\shM)=nd$. Seeking a contradiction, we assume $n\neq 0$,
such that  $\deg(\shM)>0$. This implies that there is only the zero map $\shM\ra\Omega^1_{F/k}=\O_F$.
It follows that $\shM^\vee\ra\Omega^1_{F/k}$ is surjective, thus bijective, contradicting
$\deg(\shM^\vee)<0$.

Finally, suppose that condition \eref{Cuspidal} holds: Consider the horizontal Cartier divisor $C\subset X$ that is a
rational cuspidal curve, 
and let $d>0$ be the degree of the map $C\ra\PP^1$. 
Now we have an exact sequence
$$
0\lra \O_C(-C)\lra \Omega^1_{X/k}|_C\lra\Omega^1_{C/k}\lra 0.
$$
To compute the term on the right, write the rational cuspidal curve with two affine charts as 
$$
C=\Spec k[t^2,t^3] \cup \Spec k[t^{-1}].
$$
On the first chart, the K\"ahler differentials are generated by $dt^2, dt^3$ modulo the
relation $t^4dt^2$. On the second chart, $dt^{-1}$ is a generator. On the overlap,
we have $dt^3=t^2dt= t^4\cdot dt^{-1}$. 
It follows that  the coherent sheaf  $\shN=\Omega^1_{C/k}/\text{\rm Torsion}$
is invertible with $\deg(\shN)=-4$. As in the preceding paragraph, $\Omega^1_{X/k}|_C = \shM\oplus\shM^\vee$
with $\deg(\shM)=nd$. Hence the map $\shM\ra\shN$ vanishes,
and $\shM^\vee\ra\shN$ is surjective, thus bijective. It follows that $nd=4$,
therefore $n\neq 0$.
\qed

\begin{proposition}
\mylabel{p-closed}
Suppose that $r:S\ra X$ satisfies conditions \eref{Reduced}--\eref{Elliptic}.
Then the canonical map of restricted Lie algebras
$H^0(X,\Theta_{X/k})\ra  H^0(\PP^1,\Theta_{\PP^1/k})$ is injective,
and every non-zero vector in $H^0(X,\Theta_{X/k})$ is $p$-closed.
\end{proposition}

\proof
According to  Theorem \ref{conditions for theta}, the 
canonical map for the sheaves of K\"ahler differentials $f^*(\Omega^1_{\PP^1/k})\ra\Omega^1_{X/k}$
induces a short exact sequence
$$
0\lra \O_X(-2) \lra\Theta_{X/k}\lra f^*(\Theta_{\PP^1/k})\lra 0.
$$
In turn, we get an exact sequence
$0\ra \O_{\PP^1}(-2) \ra f_*(\Theta_{X/k}) \ra \Theta_{\PP^1/k}$,
and thus an inclusion $H^0(X,\Theta_{X/k})\subset H^0(\PP^1,\Theta_{\PP^1/k})$.

Since $\Theta_{\PP^1/k}=\O_{\PP^1}(2)$,
the restricted Lie algebra $\ideal h =  H^0(\PP^1,\Theta_{\PP^1/k})$ is 3-dimensional.
It can be seen as a semidirect product $\ideal h = \ideal a\rtimes \ideal b$,
where $\ideal a=k^{\oplus 2}$ with trivial Lie bracket and $p$-map, and $\ideal b = kb$
with $b^{[p]} = b$. The  semidirect product structure comes from the homomorphism
$\ideal b\ra\mathfrak{gl}(\ideal a)$ given by $b\mapsto \id_{\ideal a}$,
compare the discussion in \cite{Schroeer 2007}, Section 3. Explicitely, Lie bracket and $p$-map are given by
$$
[a+\lambda b,a'+\lambda'b]= \lambda a'-\lambda'a\quadand
(a+\lambda b)^{[2]} = \lambda(a+\lambda b).
$$
It is easy to see that every non-zero vector in $\ideal a\rtimes\ideal b$
is $p$-closed, whence the same holds for the restricted Lie algebra $H^0(X,\Theta_X)$.
\qed

\medskip
Now suppose that $S\ra X$ satisfies the condition \eref{Reduced}--\eref{Zariski}.
Then the  restricted Lie algebra $\ideal g=H^0(X,\Theta_{X/k})$ is two-dimensional,
and every non-zero vector is $p$-closed.
Moreover, the singularities  $x_1,\ldots,x_r\in X$ are Zariski.
Let $A_i=\O_{X,x_i}^\wedge$ be the corresponding complete local ring.
Then $A_i\ideal g =\Theta_{A_i/k}$. In turn, Proposition \ref{canonical line}
yields canonical lines $\ideal l_i\subset \ideal g$, for $1\leq i\leq r$.

\begin{theorem}
\mylabel{conditions sufficient}
Suppose $S\ra X$ satisfies conditions \eref{Reduced}--\eref{Zariski}.  
Then $X$ is a K3-like covering. More precisely, 
for each  vector field  $D\in\ideal g= H^0(X,\Theta_{X/k})$ not contained in
the union of the canonical lines $\ideal l_i\subset \ideal g$, the ensuing quotient $Y=X/G$
is a simply-connected Enriques surface. 
The Picard group $\Pic^\tau_{Y/k}$ is the  Cartier dual of $G$,
and the projection $\epsilon:X\ra Y$ is the K3-like covering.
\end{theorem}

\proof
Since the vector field $D$ avoids the canonical lines, the $G$-action is free and
the quotient $Y=X/G$ is smooth, according to Proposition \ref{canonical line}.
Furthermore, $\omega_X=\O_X$ by Theorem \ref{conditions for theta}.
In the resolution of singularities $r:S\ra X$, the smooth surface $S$ is either K3 or rational,
and in both cases the algebraic fundamental group $\pi_1(S)$ vanishes.
It follows that $\pi_1(X)$ vanishes as well.
Now Proposition \ref{enriques quotient} ensures that $Y$ is a simply-connected Enriques surface,
the quotient morphism $\epsilon:X\ra Y$ is the K3-like covering,
and $\Pic^\tau_{Y/k}$ is the Cartier dual of $G$.
\qed

\section{The twistor curves in the  moduli stack}
\mylabel{Twistor curves}

Twistor curves where introduced and studied by Ekedahl, Hyland and Shepherd-Barron \cite{Ekedahl; Hyland; Shepherd-Barron 2012}.
We now  investigate further their remarkable discovery.
Let $X$ be a normal K3-like covering, 
and suppose that there is an elliptic fibration $f:X\ra\PP^1$
whose fibers are reduced.
Combining Proposition \ref{conditions necessary} and Theorem \ref{conditions sufficient} 
one obtains:

\begin{proposition}
\mylabel{k3-like theta}
The tangent sheaf $\Theta_{X/k}$ is isomorphic to $\O_X^{\oplus 2}$,
every vector in the 2-dimensional restricted Lie algebra $\ideal g=H^0(X,\Theta_{X/k})$
is $p$-closed, and for every  point $x\in X$, we have $\O_{X,x}\ideal g=\Theta_{X/k,x}$.
Moreover, the morphism $f:X\ra\PP^1$ induces  an inclusion
of restricted Lie algebras $\ideal g\subset H^0(\PP^1,\Theta_{\PP^1/k})$.
\end{proposition}

\proof
Choose an Enriques surface $Y$ so that $\epsilon:X\ra Y$ is the K3-like covering.
Let $\varphi:Y\ra\PP^1$ be the elliptic fibration induced from $f:X\ra\PP^1$.
According to Proposition \ref{conditions necessary}, conditions 
\eref{CI}--\eref{Zariski} hold. The claim on $\Theta_{X/k}$  now follows from
Theorem \ref{conditions for theta}. Moreover, every vector space basis $D_1,D_2\in\ideal g$ yields an isomorphism
$(D_1,D_2):\O_X^{\oplus 2}\ra \Theta_{X/k}$, whence an $\O_{X,x}$-basis $D_1,D_2\in\Theta_{X/k,x}$
at each  point $x\in X$.
The assertions on the restricted Lie algebra
$\ideal g$ come from Proposition \ref{p-closed}.
\qed

\medskip
Let $x_1,\ldots,x_r\in X$ be the singularities.
The corresponding local rings $A_i=\O_{X,x_i}$ are Zariski singularities,
and we have $A_i\ideal g=\Theta_{A_i/k}$.
According to Proposition \ref{canonical line}, this gives 
canonical lines $\ideal l_i\subset\ideal g$, $1\leq i\leq r$
in the two-dimensional restricted Lie algebra $\ideal g=\Theta_{X/k}$.

We call the projective line  $\PP(\ideal g)\simeq \PP^1$ the \emph{twistor curve}
attached to the K3-like covering $X$, in analogy to   twistor spaces coming from 
the two-spheres of complex structures on   complex hyperk\"ahler manifolds. 
The closed points $t_i\in\PP(\ideal g)$
corresponding to the canonical lines ${\ideal l}_i\subset\ideal g_i$  are called the \emph{boundary points}.
We write 
$$
\PP(\ideal g)^\circ= \PP(\ideal g)\smallsetminus \{t_1,\ldots,t_r\}
$$
for the complementary open subset of  \emph{interior points}.
We saw in the last section that each interior point $t\in\PP(\ideal g)$ of the
twistor line gives a simply connected Enriques surface $Y_t=X/G_t$, where
$G_t$ is the height-one group scheme corresponding to the non-zero vector field $D\in H^0(X,\Theta_{X/k})$
giving the interior point $t\in\PP(\ideal g)$.
It is easy to see that this construction extends to a flat family 
$$
\formal Y\lra\PP(\ideal g)^\circ  
$$
of Enriques surfaces $\formal Y_t=X/G_t$.   This family can be seen as  a morphism
$$
\PP(\ideal g)^\circ\lra \shM_\Enr
$$
into the moduli stack $\shM_\Enr$ of Enriques surface.
We now seek to understand this map more closely, in particular how to verify that 
this map is non-constant.

The elliptic fibration $f:X\ra \PP^1$ on the K3-like covering induces a family of elliptic fibrations
$\formal Y\ra \PP^1\times\PP(\ideal g)^\circ$
over the interior   of the twistor line. Moreover, the nonzero vector field 
$$
D\in\ideal g=H^0(X,\Theta_{X/k})\subset H^0(\PP^1,\Theta_{\PP^1/k})
$$
induce non-zero vector fields on $\PP^1$. In turn, we get a commutative diagram
\begin{equation}
\label{fibration diagram}
\begin{CD}
X	@>\epsilon>>	 X/G\\
@VfVV			@VV\varphi V\\
\PP^1	@>>>		\PP^1/G.
\end{CD}
\end{equation}
The vertical map on the right is the induced elliptic fibration
$\varphi:Y\ra\PP^1$ on the simply-connected Enriques surface $Y=X/G$.
Since $\Theta_{\PP^1/k}=\O_{\PP^1}(2)$,
this induced vector field vanishes along a divisor $A\subset\PP^1$ of degree two.
Geometrically speaking, the   tangent vectors comprising the vector field  $D\in H^0(X,\Theta_{X/k})$
are vertical at all points  $x\in X$ lying over a point $a\in A$.
The divisor $A\subset \PP^1$ of degree $\deg(A)=2$ consists either of one or two closed points.

\begin{proposition}
\mylabel{classical or supersingular}
If $A\subset\PP^1$ consist of two points, then the simply-connected Enriques surface $Y=X/G$ is
classical. If this divisor consists of only one point, $Y$ is supersingular.
\end{proposition}

\proof
Write $\tilde{\ideal g}=H^0(\PP^1,\Theta_{\PP^1/k})$ for the three-dimensional restricted Lie algebra
that contains $\ideal g = H^0(X,\Theta_{X/k})$. 
With $\PP^1=\Spec k[t] \cup\Spec k[t^{-1}]$, we may regard 
$$
D_0=\partial/\partial t, \quad D_1=t\partial/\partial t\quadand D_2= t^2\partial/\partial t
$$
as a vector space basis for $\tilde{\ideal g}$. 
Given a linear combination $D=\sum\lambda_iD_i$, one computes $D(t)=\sum\lambda_it^i$.
A direct computation shows  that  $D^2=0$ holds if and only if 
$\lambda_1=0$. The condition $D^2=0$  means that $G=\alpha_p$, in other words that the
Enriques surface $Y$ is supersingular.
On the other hand, the condition $\lambda_1=0$ means that $D$, regarded as  a section of $\Theta_{\PP^1/k}=\O_{\PP^1}(2)$,
vanishes at a single point.
\qed

\begin{proposition}
\mylabel{multiple fibers}
The fibers $Y_b\subset Y$  over the images  $b\in \PP^1=\PP^1/G$ of the points 
points $a\in A\subset \PP^1$ are precisely the multiple fibers for $\varphi:Y\ra\PP^1$.
\end{proposition}

\proof
Write $B\subset\PP^1/G$ for the image of $A\subset\PP^1$, and consider the complementary
open subsets $U=\PP^1\smallsetminus A$ and $V=\PP^1/G\smallsetminus B$.
The $G$-action on $U$ is free, and the projection $U\ra V$ is a $G$-torsor.
The diagram \eqref{fibration diagram}, which is $G$-equivariant,
yields a morphism
$X_U\ra U\times_V Y_V$
of $G$-torsors over $Y_V$. This must be an isomorphism, by the general fact that
the categories  of torsors are groupoids.

Now let $Y_b$ be some multiple fiber, for some closed point $b\in\PP^1/G$.
Then for each closed point $y\in Y_b$, the local ring has  $\edim(\O_{Y_b,y})=2$.
Consequently, the fiber product
$\PP_1\times_{\PP^1/G} Y_b$ has embedding dimension $\geq 3$ at each closed point.
In light of the preceding paragraph, it follows that $b\not\in V$, because
our K3-like covering $X$ is normal.
Summing up, we have shown that   multiple fibers for $\varphi:X\ra\PP^1$
may occur at most over points $b\in B$.

If $A\subset\PP^1$ consists of only one point, the $Y$ is supersingular 
by Proposition \ref{classical or supersingular},  and there
is precisely one multiple fiber. If $A$  consists of two points, then $Y$ is classical, and there
are precisely two multiple fibers. In both case, we conclude that for all points $b\in B$,
the fiber $Y_b$ must be  multiple.
\qed

\medskip
It follows that if the restricted Lie algebra  $\ideal g\subset H^0(\PP^1,\Theta_{\PP^1/k})$ contains
global vectors fields  $D$ with $D^2=0$ and $D^2\neq 0$, and the former defines an interior point 
in the twistor curve,
then the flat family $\formal Y\ra\PP(\ideal g)^\circ$ of simply-connected Enriques surfaces
contains both supersingular and classical members, whence the morphism 
$\PP(\ideal g)^\circ\ra \shM_\Enr$ is non-constant.
Moreover, one may use the position of the  points $b\in \PP^1/G$ where $Y_b$ is multiple,
relative to other points $b'\in\PP^1/G$ where $Y_{b'}$ is singular,
to deduce further results on the non-constancy of the twistor construction
$\PP(\ideal g)^\circ\ra \shM_\Enr$.

Note that some simply-connected Enriques surface $Y$ admit  non-zero global vector fields.
This is automatically the case if $Y$ is supersingular, according to \cite{Cossec; Dolgachev 1989}, Proposition 1.4.2.
For $Y$ classical it holds if and only if it is a so-called \emph{exceptional Enriques} surface,
according to \cite{Ekedahl; Shepherd-Barron 2004}, Theorem B.
The problem to describe   global 1-forms rather then vector fields was solved in \cite{Katsura 1982}, Section 3.

\section{K3-like coverings and elliptic fibrations}
\mylabel{K3-like coverings fibrations}

Let $k$ be an algebraically closed ground field of characteristic $p=2$,
and $Y$ be a simply-connected Enriques surface, with K3-like covering
$\epsilon:X\ra Y$. The latter is  a  torsor for the Cartier dual $G$ of 
$\Pic^\tau_{X/k}$. Assume  that $X$ has only isolated singularities.
We  now want to study the geometry and the singularities of $X$ in more detail, by choosing an
elliptic fibration $\varphi:Y\ra\PP^1$. 

\begin{proposition}
\mylabel{stein factorization}
The Stein factorization for the composite map $\varphi\circ\epsilon:X\ra \PP^1$ is
the Frobenius map $F:\PP^1\ra\PP^1$.
\end{proposition}

\proof
This could be deduced from the diagram \eqref{fibration diagram}. Let us give a  direct, independent argument:
Write $C$ for the Stein factorization, such that we have a commutative diagram
$$
\begin{CD}
X 	@>\epsilon>>	Y\\
@V\psi VV		@VV\varphi V\\
C	@>>s>		\PP^1.
\end{CD}
$$
The surjection $s:C\ra\PP^1$ is radical,
because this holds for $\epsilon:X\ra Y$. Whence $C=\PP^1$, and the morphism
is a power of the relative Frobenius map $F:\PP^1\ra\PP^1$. Since $\deg(\epsilon)=2$,
we either have $s=F$ or $s=\id$. To rule out the latter is suffices to check
that $\shA=\varphi_*(\epsilon_*(\O_X))$ has rank $\geq 2$ as $\O_{\PP^1}$-module.
Consider first the case that $Y$ is classical, such that $G=\mu_2$.
Then $\epsilon_*(\O_X)=\O_Y\oplus\shL$ for some invertible sheaf $\shL\in\Pic(Y)$ of order two,
whence $\shL=\omega_Y=\O_Y(C_1-C_2)$, where $C_1,C_2\subset Y$ are the two half-fibers.
It immediately follows that $\rank(\shA)\geq 2$.

Finally, suppose that $Y$ is supersingular, such that $G=\alpha_2$.
Then we have a short exact sequence
$0\ra\O_Y\ra\epsilon_*(\O_X)\ra\O_Y\ra 0$,
according to \cite{Ekedahl 1988}, Proposition 1.7. In turn, we get a long exact sequence
$$
0\lra\O_{\PP^1}\lra\shA\lra \varphi_*(\O_Y)\lra R^1\varphi_*(\O_Y).
$$
Write $R^1\varphi_*(\O_Y)=\shL\oplus\shT$ for some invertible sheaf $\shL$ and some torsion sheaf $\shT$.
The Canonical Bundle Formula (\cite{Bombieri; Mumford 1977}, Theorem 2) ensures that $\shL=\O_{\PP^1}(-2)$.
It follows that the coboundary map $\O_{\PP^1}=\varphi_*(\O_Y)\ra R^1\varphi_*(\O_Y)$ has nontrivial kernel,
and again $\rank(\shA)\geq 2$.
\qed

\medskip
Write  $f:X\ra \PP^1$ for the Stein factorization, which sits in the  commutative diagram
$$
\begin{CD}
X	@>\epsilon>>	Y\\
@VfVV			@VV\varphi V\\
\PP^1	@>>_F>		\PP^1.
\end{CD}
$$
We thus get a canonical morphism $X\ra Y^\fpb$ from the K3-like covering to the \emph{Frobenius pullback}
$Y^\fpb=Y\times_{\PP^1}\PP^1$.

\begin{proposition}
\mylabel{normalization}
The morphism $X\ra Y^\fpb$ is the normalization.
Its ramification locus consists   of the preimages $\epsilon^{-1}(Y_b)\subset X$
of the multiple fibers $Y_b$, $b\in \PP^1$.
\end{proposition}

\proof
Both maps $X\ra Y$ and $Y^\fpb\ra Y$ have degree two, whence $X\ra Y^\fpb$ is birational.
It thus must be the normalization, because we assume throughout that $X$ is normal.
For each closed point $y\in Y_b$ lying on a multiple fiber, the embedding dimension is
$\edim(\O_{Y_b,y})=2$, so the local rings for each closed point on the preimage of $Y_b$
in the Frobenius pullback $Y^\fpb$ 
has embedding dimension three. It follows that the preimage $\epsilon^{-1}(Y_b)\subset X$ belongs to the
ramification locus, which is the locus on $X$ defined by the conductor ideal. On the other hand, $Y^\fpb$ is smooth at each point mapping to a point
$y\in Y$ where the morphism $\varphi:Y\ra\PP^1$ is smooth. 
Since all   fibers of $\varphi:Y\ra\PP^1$ are reduced by Theorem  \ref{properties genus-one fibrations}, 
it follows that  the ramification locus for $X\ra Y^\fpb$ consists precisely of the preimages of the multiple fibers.
\qed

\medskip
Next, consider the \emph{jacobian fibration} $J\ra\PP^1$ attached to the elliptic fibration 
$\varphi:Y\ra\PP^1$, which is a relatively minimal elliptic fibration endowed with a section $O\subset J$.
According to \cite{Cossec; Dolgachev 1989}, Theorem 5.7.2 the smooth surface $J$ is rational. Moreover,
all such rational elliptic surfaces $J\ra\PP^1$ were classified by Lang \cite{Lang 2000}.
In some sense, we completely understand such jacobian fibrations.
Our strategy throughout is to relate the rational elliptic surface $J\ra\PP^1$ 
to the fibration $f:X\ra\PP^1$ on the K3-like covering.

\begin{proposition}
\mylabel{fibers jacobian}
For each closed point $b\in \PP^1$, the fibers $Y_b$ and $J_b$ have the same Kodaira type. 
In particular, the rational elliptic surface $J\ra\PP^1$ has only reduced
fibers, and the possible Kodaira types
for  the singular fibers
are  $\I_n$ with $1\leq n\leq 9$, $\II$, $\III$ or $\IV$.
\end{proposition}

\proof
If $b\in \PP^1$ is a closed point whose fiber $Y_b$ is non-multiple, then
$J\ra \PP^1$ admits a section over the formal completion $\Spec(\O_{\PP^1,b}^\wedge)$, 
which yields an isomorphisms
$$
J\otimes_{\PP^1}\Spec(\O_{\PP^1,b}^\wedge)  \simeq Y\otimes_{\PP^1}\Spec(\O_{\PP^1,b}^\wedge).
$$
For the multiple fibers $Y_b$, one can at least say that 
$J_b$ and $Y_b$ have the same Kodaira type,
according to  a general result of 
Liu, Lorenzini and Raynaud \cite{Liu; Lorenzini; Raynaud 2005}, Theorem 6.6.
The second assertion follows from Theorem  \ref{properties genus-one fibrations}.
\qed

\medskip
To proceed, consider for the rational elliptic surface $J\ra\PP^1$ the Frobenius pullback 
$$
X'=J^\fpb = J\times_{\PP^1}\PP^1.
$$
According to Proposition \ref{frobenius pullback zariski} below, this is a normal surface, having only isolated
Zariski singularities, and the  dualizing sheaf is $\omega_X=\O_X$.
Clearly, the elliptic fibration $X'\ra\PP^1$ admits a section.
Write $\eta\in \PP^1$ for the generic point.

\begin{proposition}
\mylabel{birational sections}
We have $X'_\eta\simeq X_\eta$ if and only if the fibration $f:X\ra\PP^1$ admits
a section.
\end{proposition}

\proof
The condition is obviously necessary. Conversely, suppose that $f:X\ra\PP^1$ admits a section.
By definition, the generic fiber of $J\ra\PP^1$ is the jacobian for the generic fiber of $Y\ra\PP^1$.
In turn,  $X'_\eta=(J^{(2/\PP^1)})_\eta$ is the jacobian for $(Y^{(2/\PP^1)})_\eta$.
The latter coincides with $X_\eta$, by Proposition \ref{normalization}.
Thus $X'_\eta\simeq X_\eta$.
\qed

\medskip
We need to understand this condition better, in order to exploit the connection 
between the K3-like covering $X$ of the Enriques surface $Y$ 
and the Frobenius pullback $X'$ of the rational elliptic surface $J$.
Let us call a curve $A\subset Y$ on the Enriques surface  a \emph{two-section} if 
$A$ is an irreducible curve that is  horizontal with respect to the elliptic fibration $\varphi:Y\ra\PP^1$,
of relative degree two. We say that $A\subset Y$ is a \emph{radical two-section}
if the surjection $A\ra\PP^1$ is radical, hence a universal homeomorphism.

\begin{proposition}
\mylabel{two-section = section}
Let  $A\subset Y$ be a radical two-section for $\varphi:Y\ra\PP^1$. Then
$B=\epsilon^{-1}(A)_\red$ is a section for $f:X\ra\PP^1$.
Conversely, the image $A=\epsilon(B)$ for any section $B\subset X$  
is a radical two-section. Moreover, such  $A$ and  $B$ exists if and only if $X'_\eta\simeq X_\eta$.
\end{proposition}

\proof
The first assertion is obvious, because $X$ is the normalization of $Y^{(2/\PP^1)}$
and the two morphisms $F:\PP^1\ra\PP^1$ and $\varphi:A\ra\PP^1$ coincide generically.
Conversely, suppose that $B\subset X$ is a section. Its image $A=\epsilon(B)$
is an integral Cartier divisor, and the composition
$B\ra A\ra\PP^1$ is radical of degree two.
We thus have either $\deg(B/A)=1$ and $\deg(A/\PP^1)=2$, or $\deg(B/A)=2$ and  $\deg(A/\PP^1)=1$.
The latter case is impossible, because the elliptic fibration $\varphi:Y\ra \PP^1$ has no sections.
Thus $A\subset Y$ must be  a two-section.
The last statement follows from Proposition \ref{birational sections}.
\qed

\medskip 
A curve $A\subset Y$ is called \emph{rational} if it is reduced, and its
normalization is isomorphic to the projective line $\PP^1$.

\begin{proposition}
\mylabel{radical = rational}
A two-section $A\subset Y$ is radical if and only if the curve $A$ is rational.
\end{proposition}

\proof
The condition is clearly necessary. Conversely, suppose that the two-section $A\subset Y$
is rational, and let $\PP^1\ra A$ be the normalization map. The induced $G$-torsor
$X_{\PP^1}=\PP^1\times G$ is trivial. It follows that the Weil divisor $\epsilon^{-1}(A)\subset X$
is non-reduced, thus its reduction $B$ yields a section for $f:X\ra \PP^1$.
According to Proposition \ref{two-section = section}, the image $A=\epsilon(B)$ is a radical two-section.
\qed

\begin{proposition}
\mylabel{two-section = half-fiber}
Suppose there is another elliptic fibration $\varphi':Y\ra\PP^1$ that is orthogonal
to $\varphi:Y\ra\PP^1$, and has a  singular half-fiber $C\subset\varphi'^{-1}(b')$.
Then some irreducible component $A\subset C$ is a radical two-section
for $\varphi:Y\ra\PP^1$.
\end{proposition}

\proof
By the classification of singular fibers, each irreducible component $A\subset C$ is a rational curve.
Let $F=\varphi^{-1}(b)$ be a fiber, such that $(C\cdot F)=2$. Since $F$ is movable, 
there is some irreducible component $A\subset C$ with $1\leq (A\cdot F)\leq 2$.
We must have $(A\cdot F)=2$, because $\varphi:Y\ra\PP^1$ admits no section.
The assertion now follows form Proposition \ref{radical = rational}. 
\qed

\medskip
For the applications we have in mind, it is permissible to replace the chosen
elliptic fibration $\varphi:Y\ra\PP^1$ by another, more suitable one:

\begin{proposition}
\mylabel{radical two-section}
Suppose that the simply-connected Enriques surface $Y$ contains a $(-2)$-curve $A\subset Y$,
or a rational cuspidal curve $C$ that is not movable. Then
there is an elliptic fibration $\varphi:Y\ra\PP^1$ admitting a radical two-section.
\end{proposition}

\proof
Suppose there is a $(-2)$-curve $A\subset Y$. Then the desired elliptic fibration exists.
In characteristic $p\neq 2$, Cossec showed that  
the desired genus-one fibration exists (\cite{Cossec 1985}, Theorem 4). 
This carries over to characteristic $p=2$, according to Lang (\cite{Lang 1988}, Theorem A.3.
The fibration is indeed elliptic by Theorem \ref{properties genus-one fibrations}.

Now suppose that there is a non-movable rational cuspidal curve $C\subset Y$.
Then it is a half-fiber of some genus-one fibration $\psi:Y\ra\PP^1$.
Again by Theorem \ref{properties genus-one fibrations}, there is an orthogonal fibration $\varphi:Y\ra\PP^1$.
In both cases, the two-sections are radical, by  Proposition \ref{radical = rational}.
\qed

\medskip 
We do not know if the conditions of the preceding proposition holds for
\emph{all} simply-connected Enriques surfaces. Note that there are Enriques surfaces $Y$ without $(-2)$-curves, according to
\cite{Lang 1983}, Theorem 4.3. 
In any case, the exceptions must be restricted in the following sense:

\begin{proposition}
\mylabel{no radical two-section}
Suppose  the simply-connected Enriques surface $Y$ has no elliptic fibration admitting  a radical two-section.
Then for every  elliptic fibration $\varphi:Y\ra \PP^1$, the following holds:
The singular fibers have Kodaira type $\I_1$ or $\II$, and only the following configurations are possible:
$$
1^{12},\quad \II+1^8,\quad \II+1^6,\quad \II+1^5,\quad \II+\II+\II,\quad \II+\II,\quadand   \II.
$$
Moreover, the half-fibers are smooth.  
\end{proposition}

\proof
In light of Proposition \ref{radical two-section} and Proposition \ref{two-section = half-fiber}, the Enriques surface
$Y$ contains neither $(-2)$-curves nor non-movable rational cuspidal curves.
In turn, only $\I_0$, $\I_1$  or $\II$ are possible Kodaira types, and the half-fibers are smooth.
According to Proposition \ref{fibers jacobian}, the Kodaira types of the fibers $Y_b$ and $J_b$ coincide, for all closed points
$b\in \PP^1$.
Looking at Lang's classification \cite{Lang 2000} of rational elliptic surfaces,
we obtain the list of possible configurations.
\qed

\medskip
Let us now assume that $\varphi:Y\ra\PP^1$ admits a radical two-section.
Then we have an isomorphism $X_\eta=X'_\eta$, an in particular
the normal surfaces $X$ and $X'$, which both have trivial canonical class,
are birational. Using terminology from the \emph{minimal model program} in dimension three and higher,
one may regard $X$ as a \emph{flop} of $X'$, because their canonical class  neither became more negative
nor more positive.

Choose a common resolution of singularities $X'\stackrel{r'}{\leftarrow} S\stackrel{r}{\rightarrow} X$.
The induced projection $h:S\ra X'\ra\PP^1$ is an elliptic fibration on the smooth surface $S$,
and we write $S\ra J'$ for the relative minimal model. Locally, $J'$ is obtained
by applying the Tate Algorithm to the  Frobenius base-change of the  Weierstra\ss{}  equation for $J$.
Summing up, we have the following commutative diagram: 
$$
\xymatrix{
	&	& S\ar[ddll]\ar[dd]_{r'}\ar[dr]^r \\
	&	&				& X\ar[rr]^\epsilon\ar[dd]_(.3)f	& 			& Y\ar[dd]^\varphi\\
J'	&	& X'\ar[dr]\ar[rr]|\hole	&				& J\ar[dr]		& \\
	&	&  				&  \PP^1\ar[rr]_F		& 			& \PP^1
}
$$
If $X$ and $X'$ have only rational singularities, then their minimal resolution of singularities
coincides with the relatively minimal model of $X_\eta\simeq X_\eta'$, whence both 
are given by $S$. It follows that  $S$ is a K3-surface, and 
the morphism $S\ra J'$ is
an isomorphism. In this case, $X'\leftarrow S\rightarrow X$ are both crepant resolutions.
On the other hand, if $X$ has an elliptic singularity, such that $S$ and hence $J'$  is  a rational surface,
we have $K_{J'}=-F$ for some fiber $F\subset J'$. Now one should regard $J'$ as
a \emph{flip} of both $X'$ and $X$, because the canonical class became more negative.
We observe:

\begin{lemma}
\mylabel{dual graph}
For all closed points $a\in \PP^1$, with Frobenius image $b=F(a)\in\PP^1$,
the four curves $Y_b$, $X_a$, $X'_a$ and $J_b$ have the same dual graph.
\end{lemma}

\proof
The curves $X_a$ and  $Y_b$ have the same dual graph, because the K3-like covering  $\epsilon:X\ra Y$ is
a universal homeomorphism. The same argument applies to the curves $X'_a$ and $J_b$.
Finally, the curves $Y_b$ and $J_b$ have the same dual graph
according \cite{Liu; Lorenzini; Raynaud 2005}, Theorem 6.6.
\qed

\medskip
The singularities on the K3-like covering $X$ of the Enriques surface $Y$ are closely related to the
singularities on the Frobenius pull-back $X'$ of the rational elliptic surface $J$, 
at least outside the multiple fibers:

\begin{proposition}
\mylabel{formal isomorphism non-multiple}
For each closed point $a\in \PP^1$ whose Frobenius image $b\in\PP^1$ has
a non-multiple fiber  $Y_b$, there is  an isomorphism of two-dimensional schemes
$$
X\times_{\PP^1} \Spec(\O_{\PP^1,a}^\wedge) \simeq X'\times_{\PP^1}\Spec(\O_{\PP^1,a}^\wedge).
$$
\end{proposition}

\proof
Write $R=\O_{\PP^1,b}^\wedge$ for the complete local ring of the point $b\in\PP^1$.
Since $Y_b$ is non-multiple, the induced projection $Y\times_{\PP^1}\Spec(R)\ra\Spec(R)$
admits a section, which in turn 
induces an identification $Y\times_{\PP^1}\Spec(R) = J\times_{\PP^1}\Spec(R)$.
In light of Proposition \ref{normalization}, taking the fiber product with the Frobenius morphism $F:\PP^1\ra\PP^1$
yields the assertion.
\qed

\section{Ogg's Formula and Frobenius base-change}
\mylabel{Ogg's formula}

In this section   let $k$ be an algebraically closed ground field of arbitrary characteristic $p>0$,
and $J\ra B$ be a    smooth elliptic surface that  is jacobian and relatively minimal.
For simplicity, we assume that  $B$ is a smooth proper connected curve, such that  $J$ is a smooth proper connected surface, 
although the analysis also applies to  local situations as well. 

Let us write $B^{(p)}$ for the scheme $B$, endowed with the new structure morphism
$B^{(p)}\ra\Spec(k)\stackrel{F}{\ra}\Spec(k)$, obtained by transport of structure with
the absolute Frobenius of $k$. Then the absolute Frobenius on $B$ becomes the \emph{relative
Frobenius morphism} $B^{(p)}\ra B$. The cartesian square
$$
\begin{CD}
J^{(p/B)}	@>>> 	J\\
@VVV		@VVV\\
B^{(p)} @>>>	B\\
\end{CD}
$$
defines an integral proper connected surface $J^{(p/B)}$ endowed with an elliptic fibration.
We call $J^{(p/B)}$ the \emph{Frobenius pullback}. Let us first record:

\begin{lemma}
\mylabel{frobenius pullback zariski}
The singularities on the Frobenius pullback $J^{(p/B)}$ are Zariski singularities.
\end{lemma}

\proof
Fix a closed point $x\in J^{(p/B)}$, and write $a\in B^{(p/B)}$, $y\in J$ and $b\in B$ for its images.
Choose a uniformizer $\pi\in\O_{B,b}$, and write $z=\pi$ for the resulting uniformizer $z\in\O_{B^{(p)},a}$.
Then the morphism of complete local $k$-algebras
$$
k[[z]]\lra \O_{B^{(p)},b}^\wedge,\quad \sum\lambda_iz^i\longmapsto\sum\lambda_i^pz^i
$$
is bijective. Using this identification, we may regard the extension $\O_{B,b}^\wedge\subset\O_{B^{(p)},a}^\wedge$
as $k[[\pi]]\subset k[[\pi,z]]/(z^p-\pi)$.
In turn, we get
$$
\O_{J^{(p/B)},x}^\wedge = \O_{J,y}^\wedge[z]/[z^p-\pi],
$$
where the uniformizer $\pi\in\O_{B,b}^\wedge$ becomes a non-unit $\pi\in \O_{J,y}^\wedge$.
Since the local ring $\O_{J,y}$ is regular, the local ring $\O_{J^{(p/B)},x}^\wedge$ is a Zariski singularity.
\qed

\medskip
Let $S\ra J^{(p/B)}$ be  a resolution of singularities, which inherits an elliptic fibration 
$S\ra B$, and write  $S\ra J'$ for the
contraction to the relative minimal model. Then we have  a new
smooth surface $J'$ endowed with an elliptic fibration $J'\ra B^{(p)}$, which is jacobian and relatively minimal.
We call $J'$ the \emph{smooth Frobenius pullback}.
Note that in general, the Kodaira dimension of the surfaces $J$ and $J'$ are different.

Fix a section $O\subset J$ for the jacobian elliptic fibration $J\ra B$, and let $W\ra\PP^1$
be the ensuing Weierstra\ss{}  fibration. The normal surface $W$ is obtained by contracting
all vertical irreducible curves disjoint to the section.
Locally at each point $b\in B$,
it is given by some minimal Weierstra\ss{}  equation
$$
y^2+a_1xy+a_3y = x^3+a_2x^2+a_4x+a_6 \qquad a_i\in\O_{B,b}.
$$
We denote by $v_b\geq 0$ the \emph{valuation of the discriminant} $\Delta_b\in\O_{B,b}$ 
for any such minimal Weierstra\ss{}  equation, 
by $m_b\geq 1$ the \emph{number of irreducible components} in the fiber $J_b$, 
and by $\delta_b\geq 0$ the \emph{wild part of the
conductor} for the Galois representation on the geometric generic $l$-torsion points.
Here $l\neq p$ is any prime different from the characteristic.
If the fiber $J_b$ is of semistable, that is of type $\I_m$, then $\delta_b=0$ and $v_b=m$.
If the fiber is unstable,  Ogg's Formula  gives
\begin{equation}
\label{ogg formula}
v_b = 2+\delta_b + (m_b-1).
\end{equation}
For a proof, see the original paper \cite{Ogg 1967}  and also  \cite{Saito 1988}, which
containes a complete proof and generalizations.

The closed point $b\in B $ may also be seen as a closed point on the scheme $B^{(p)}$.
To avoid confusion, we denote it by $a\in B^{(p)}$.
This gives numerical invariants    $v_a,\delta_a,m_a$ for the smooth Frobenius pullback $J'\ra B$.
Now observe that we may obtain  $J'$ from the base-changed Weierstra\ss{}  fibration  
$W^{(p/B)}$ as follows: Run  the Tate Algorithm \cite{Tate 1972}
until the Weierstra\ss{}  equation becomes minimal and then do the minimal resolution of singularities.
This gives another invariant: The \emph{length of the Tate Algorithm}
$\lambda_a\geq 0$, that is, the number  of repetitions that are necessary to finish   the Tate Algorithm.
The case $\lambda_a=0$ means that the Weierstra\ss{}  equation is already minimal.

If the fiber is semistable, we have $\lambda_a=0$ and $m_a=pm_b$.
Indeed, a local computation shows that 
on $W^{(p/B)}$ we get rational double points of type $A_{p-1}$ lying over the
singular points of the semistable fiber $J_b$.
The situation is more interesting at the unstable fibers:

\begin{lemma}
\mylabel{unstable fiber}
Suppose the fiber $J_b$ is unstable. Then $J'_a$ is unstable as well,
and its numerical invariants are given by the formulas
$$
v_a=pv_b-12\lambda_a,\quad \delta_a=\delta_b\quadand m_a=pm_b+(p-1)(\delta_b+1) -12\lambda_a.
$$
In characteristic two, this means $v_a=2v_b-12\lambda_a$ and $m_a=2m_b+\delta_b+1-12\lambda_a$.
\end{lemma}

\proof
Among all function field extensions that achieve semistable reduction, there
is a smallest one, and this smallest one is a Galois extension, which follows from  the N\'eron--Ogg--Shafarevich Criterion
(see for example \cite{Bosch; Luetkebohmert; Raynaud 1990}, Section 7.4, Theorem 5).
In particular, unstable fibers stay unstable under purely inseparable field extensions.
Each round of the Tate Algorithm reduces the valuation of the discriminant by 12.
The wild part of the conductor depends on a Galois representation, and is thus not
affected by purely inseparable extensions. The   equation for the number  $m_a\geq 1$ of irreducible components 
is now a consequence of Ogg's Formula \eqref{ogg formula}.
\qed

\medskip
By construction, the Weierstra\ss{}  model $W$ is normal, and the fibers
$W\ra B$ contain at most one non-smooth point, which is a rational double point.
In turn,   the  base-change $W^{(p/B)}$ is normal, and contains fiber-wise at most
one singularity.

\begin{proposition}
\mylabel{singularity rational}
We have $\lambda_a>0$ if and only if the normal surface $W^{(p/B)}$ contains a non-rational singularity
lying over the point $a\in B^{(p)}$.  
\end{proposition}

\proof
There is nothing to prove if $W^{(p/B)}$ is smooth along the fiber. Now suppose that a singularity
$x\in W^{(p/B)}$ is present, such that the fiber $C=W_a^{(p/B)}$ is a rational cuspidal curve.
Let $S\ra W^{(p/B)}$ be the minimal resolution of this singularity.
It the singularity is rational, then $S$ is the smooth model at $a\in B^{(b)}$, 
thus the Weierstra\ss{}  equation was minimal, whence $\lambda_a=0$.
Conversely, suppose that the singularity is non-rational.  Then $W^{(p/B)}$
cannot be a Weierstra\ss{}  model, so the Weierstra\ss{}  equation was not minimal,
that is, $\lambda_a>0$.
\qed

\begin{proposition}
\mylabel{rational elliptic}
If both $J$ and $J'$ are rational surfaces, then we have $\sum\lambda_a=p-1$,
where the sum runs over all closed points $a\in B^{(p)}=\PP^1$.
\end{proposition}

\proof
Let us write $v=\sum v_b$, $v'=\sum v_a$  and $\lambda=\sum\lambda_a$ for the sums of   local invariants.
Since $J$ and $J'$ are rational, $v=v'=12$ must hold, according to \cite{Lang 2000},
Lemma 0.1. On the other hand, we have $v'=pv-12\tau$, which forces $\tau=p-1$. 
\qed

\medskip
In particular, in characteristic $p=2$ there is precisely one fiber that
does not stay minimal, and it becomes minimal after one round of the Tate Algorithm.
We shall analyze this particular situation in the next section.

\section{Geometric interpretation of the Tate Algorithm}
\mylabel{Tate algorithm}

The \emph{Tate Algorithm}  turns arbitrary Weierstra\ss{} equations over   discrete valuation rings
into  minimal Weierstra\ss{}  equations \cite{Tate 1972}. This algorithm is of paramount importance
for the arithmetic theory of elliptic curves, as well as the structure of elliptic surface. 
The goal of this section 
is to describe the geometry behind the algorithm. This is perhaps  well-known, but we could 
not find a suitable reference.

Let $R$ be a discrete valuation ring with residue field $k=R/\maxid_R$ and field of
fractions $F=\Frac(R)$, and choose a uniformizer $\pi\in R$.
Suppose we have a Weierstra\ss{}  equation
$$
y^2 +a_1xy + a_3y = x^3 + a_2x^2 + a_4x + a_6
$$
with coefficients $a_i\in R$ and discriminant $\Delta\neq 0$.
Setting $y=X/Z$, $x=X/Z$ and multiplying with $Z^3$ gives its homogenization
$$
Y^2Z + a_1XYZ + a_3YZ^2 = X^3 + a_2 X^2Z + a_4XZ^2 + a_6Z^3,
$$
which defines a relative cubic $J\subset \PP^2_R=\Proj R[X,Y,Z]$
endowed with a section $O\subset J$ given by the equations $X=Z=0$,
such that $\O_J(1)=\O_J(3O)$.

\begin{theorem}
\mylabel{tate geometric}
Let $S\ra J$ be the blowing-up with center given by $\pi=X=Z=0$,
and $S\ra J'$ be the blowing-down of the strict transform of the closed fiber $J_0$.
Then $J'\ra\Spec(R)$ is a relative cubic   given by the Weierstra\ss{}  equation
\begin{equation}
\label{nonminimal Weierstrass}
y^2 + a_1'xy + a'_3y = x^3 + a_2'x^2 + a_4'x + a'_6
\end{equation}
whose  coefficients are $a'_i=a_i\pi^i$.
\end{theorem}

\proof
Write the closed fiber of the blowing-up as $S_0=C\cup E$, where $C$ is the strict transform of $J_0$,
and $E$ is the $(-1)$-curve.
In light of  the exact sequence
$$
1\lra \O_{S_0}^\times\lra \O_C^\times\oplus\O_E^\times\lra \kappa(s)^\times\lra 1,
$$
where $s\in C\cap E$ is the rational intersection point, 
we easily see that the  restriction map $\Pic(S_0)\ra\Pic(C)\oplus\Pic(E)$
is injective. Consider the invertible sheaf $\shL=\O_S(3O+2E)$, which has 
$(\shL\cdot C)=2$ and $(\shL\cdot E)=1$.
Using Riemann--Roch and the Theorem on Formal Functions, one easily infers that such $\shL$ 
is globally generated, with $R^1h_*(\shL)=0$ and $h_*(\shL)$ free of rank three. Here $h:S\ra\Spec(R)$
is the structure morphism.  
We thus get a morphisms $S\ra \PP^2_R$. Its image $X$ is an effective Cartier divisor. Base-changing to 
the generic fiber, we see that  $X\subset\PP^2_R$ is a relative cubic.

Me may compute the  cubic equation for  $X\subset\PP^2_R$ by regarding $X$ as
a closure of some affine cubic.
Consider first the affine open subset
$D_+(Y)\subset\PP^2_R$. Here the relative cubic $J$ is given by the dehomogenized cubic equation 
$$
\frac{Z}{Y} + 
a_1\frac{X}{Y}\frac{Z}{Y} + 
a_3\left(\frac{Z}{Y}\right)^2 
=
\left(\frac{X}{Y}\right)^3 + 
a_2 \left(\frac{X}{Y}\right)^2\frac{Z}{Y} + 
a_4\frac{X}{Y}\left(\frac{Z}{Y}\right)^2 + 
a_6\left(\frac{Z}{Y}\right)^3.
$$
The blowing-up $S\ra J$ is covered by 
two affine charts, and we look at the $\pi$-chart given by the indeterminates $ X/\pi Y,Z/\pi Y$ over $R$,
subject to the relation
\begin{multline*}
\frac{Z}{\pi Y} + 
a_1\pi\frac{X}{\pi Y}\frac{Z}{\pi Y} + 
a_3\pi\left(\frac{Z}{\pi Y}\right)^2 
=\\
\pi^2\left(\frac{X}{\pi Y}\right)^3  + 
a_2\pi^2\left(\frac{X}{\pi Y}\right)^2 \frac{Z}{\pi Y}  + 
a_4\pi^2 \frac{X}{\pi Y} \left(\frac{Z}{\pi Y}\right)^2 + 
a_6\pi^2\left(\frac{Z}{\pi Y}\right)^3.
\end{multline*}
Via the substitutions $X/\pi Y=U/V $ and $Z/\pi Y=W/V$, we may view this 
as a closed   subscheme inside $D_+(V)\subset \PP^2_R=\Proj R[U,V,W]$.
Its closure in $\PP^2_R$ is given by the homogeneous cubic equation
\begin{equation}
\label{homogeneous cubic}
 V^2W + a_1\pi UVW + a_2\pi VW^2 = \pi^2 U^3 + \pi^2a_2U^2W + \pi^2a_4UW^2 + \pi^2a_6W^3.
\end{equation}
The closed fiber is thus given $\pi=V^2W=0$,  which is the union of the line $W=0$
and the double line $V^2=0$  in $\PP^2_k$. Indeed, since $S\smallsetminus C$ coincides with the $\pi$-chart of the
blowing-up,   $S\smallsetminus C \ra X$ is an open embedding, and the above closure equals $X$.

Now consider the automorphism
$$
f=\begin{pmatrix}
\pi^{-2} \\
& \pi^{-2}\\
&& 1
\end{pmatrix}
\in\PGL_3(F) =\Aut(\PP^2_F).
$$
Over the function field, it gives a new cubic $f(J_F)\subset\PP^2_F$. Applying $f$
to the cubic equation \eqref{homogeneous cubic} and multiplying  by $\pi^4$, we get the homogeneous equation
$$
V^2W + a_1'UVW + a_3'VW = U^3 + a_2'U^2W + a_4'UW^2 + a_6'W^3,
$$
where $a_i'=\pi^ia_i$. Dehomogenizing the above homogeneous  equation with $x=U/W$ and $y=V/W$ 
gives the desired Weierstra\ss{}  equation \eqref{nonminimal Weierstrass}.
Thus the   relative cubic $J'\subset \PP^2_R$ defined by the above homogeneous equation
is the schematic closure of  $f(J_F)$.
Note that the closed fiber $J'_k$ is a rational cuspidal curve,
whence $J'$ is normal. 

Clearly, the rational map $f:\PP^2_R\dashrightarrow \PP^2_R$ is defined on the generic fiber $\PP^2_K$
and the standard affine open subsets $D_+(U)$ and $D_+(V)$, whence it is defined everywhere
except at the point $z=(0:0:1)\in\PP^2_k$.
In turn, the rational map $f:J\dashrightarrow J'$ is defined outside $z\in J$.
The latter is nothing but the 
intersection of the line   $\pi=W=0$ and the double line $\pi=V^2=0$ on the closed fiber.
Every section $(a:b:0)$ for $J\ra\Spec(R)$ is mapped to itself under the rational
map $f:\PP^2_R\dashrightarrow\PP^2_R$, whence $f:J\dashrightarrow J'$ is a bijective morphism on the complement
of the double line $\pi=V^2=0$. Now choose a modification  $S'\ra J$ with center $z\in J$ so that
the rational map $J\dashrightarrow J'$ comes from a morphism $S'\ra J'$.
Let $E\subset S'_0$ be an irreducible component different from the strict transform
of the line $\pi=W=0$. By Zariski's Main Theorem,  it cannot map to a regular point on $J'_0$, because
these are already images of points form the line  $\pi=W=0$.
In turn, it must map to the singular point of the rational cuspidal curve $J'_0$.
Thus the morphism $S'\ra J'$ factors over $J$, and furthermore  the induced morphism $J\ra J'$
contracts the double line $\pi=V^2=0$.
It follows that $S\ra J'$ is the contraction of $C\subset S$.
\qed

\medskip
This observation gives some information on elliptic singularities, which turns out very useful.
Let us call a Weierstra\ss{}  equation 
\begin{equation}
\label{almost minimal equation}
y^2 +a'_1xy + a'_3y = x^3 + a'_2x^2 + a'_4x + a'_6
\end{equation}
\emph{almost minimal}  if its discriminant $\Delta\in R$ is nonzero,
and the Weierstra\ss{}  equation  becomes minimal after one round in the Tate Algorithm.

\begin{proposition}
\mylabel{almost minimal}
Every almost minimal Weierstra\ss{}  equation
defines an elliptic singularity, for which the exceptional divisor on the minimal resolution
of singularities coincides with the closed fiber on the regular minimal model $J\ra\Spec(R)$.
The intersection numbers on the exceptional divisor and the fiber divisor coincide,
except that for one component that has
multiplicity one in the fiber, the intersection number $0$ or $-2$ become
 $-1$ respective  $-3$ on the exceptional divisor.
\end{proposition}

\proof
Running through the Tate Algorithm, we   may assume
that $\pi^i| a'_i$. Now Theorem \ref{tate geometric} describes how the singularity on $J'$
arises from $J$ via the birational correspondence $J\leftarrow S\rightarrow J'$.
The statement on the exceptional divisor and the intersection numbers   follow.  
It remains to check that the singularity is indeed elliptic.
Write the closed fiber as $S_0=E+C+D$, where $E$ is the $(-1)$-curve,
$C$ is the irreducible component that intersects $E$, and $D$ be the union of the remaining
irreducible components. The contraction $r':S\ra J'$ factors over the contraction $S'\ra J'$ of $D$,
which introduces  only rational double points. In light of the Leray--Serre spectral sequence,
it  suffices to show that for 
the induced map $g:S'\ra J'$ the higher direct image  $R^1g_*(\O_{S'})$ has length one.

Changing notation, we write the closed fiber $S'_0=E+C$, where $E=\PP^1$
and $C$ is a rational cuspidal curve. The intersection point $E\cap C$ is smooth in $S'$,
whence both $E,C\subset S'$ are Cartier divisors, with intersection matrix
$N=(\begin{smallmatrix}-1&1\\1&-1\end{smallmatrix})$, so $\deg\O_C(-C)=1$. By the Theorem on Formal Functions,
we have to show that $h^1(\O_{nC})=1$ for all $n\geq 1$.
This is obvious for $n=1$. We proceed by induction, using the long exact sequence
attached to the short exact sequence
$$
0\lra \O_C(-nC)\lra \O_{(n+1)C}\lra \O_{nC}\lra 0,
$$
together with the fact that $\shL=\O_C(-nC)$ has degree $n\geq 1$, whence $H^1(C,\shL)=0$
by Riemann--Roch.
\qed

\medskip
Note that   blowing-ups of  elliptic singularities are  usually non-normal.
Furthermore, the sequence of normalized blowing-ups of reduced singular points
usually produces a non-minimal resolution. This makes it often so difficult
to compute a resolution of singularities explicitly. The above result
avoids all these complications and makes it unnecessary to cope with explicite
equations.
For our goals, the cases of Kodaira  type $\II$ is  particularly important:

\begin{corollary}
\mylabel{almost minimal II}
Every almost minimal Weierstra\ss{}  equation 
of type $\II$ defines an elliptic singularity, for which the exceptional divisor $E$ on the 
minimal resolution of singularities is a rational cuspidal curve   with intersection matrix
$N=(-1)$.
\end{corollary}

\medskip
This is one of the simplest cases of Wagreich's classification of 
elliptic double points \cite{Wagreich 1970}, Theorem 3.5.
It belongs to the family of elliptic double points, whose exceptional divisors
are symbolically described by   Wagreich in the following way:
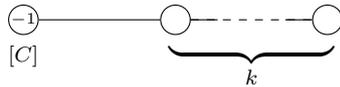
\begin{center}
\begin{tikzpicture}
\node (a) at (1,1) 	[circle, draw, label=below:{$\scriptstyle{[C]}$}] 	{};
\node (a') at (1,1) 	[] 							{$\scriptscriptstyle{-1}$};
\node (b) at (3,1)	[circle, draw] 						{};
\node (b') at (3.7,1)								{};
\node (c') at (4.3,1)								{};
\node (c) at (5,1)	[circle, draw] 						{};
\node (d) at (4,1)	[below]        {$\underbrace{\phantom{mmmmmm}}_{k}$};
\draw (a) -- (b);
\draw [dashed] (b) -- (c);
\draw (b) -- (b');
\draw (c) -- (c');
\end{tikzpicture}
\end{center}
Here the  chain of length $k\geq 0$ to the right consists of $(-2)$-curves.
Our singularity  is the special case $k=0$.
We simply say that it is the \emph{elliptic singularity obtained by contracting
a cuspidal $(-1)$-curve}.
One  may also describe it
in terms of the minimal good resolution,
where the exceptional divisor  has  simple normal crossings.
This is obtained from the minimal resolution $S$ by three further
blowing-ups, such that $E$ has four irreducible components, which 
are smooth rational curves. The dual graph takes the following form:
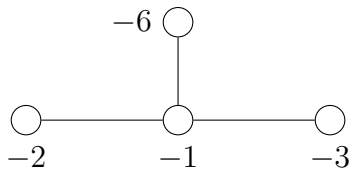
\begin{center}
\begin{tikzpicture}
[node distance = 2cm]
\node (a) at (1,1) 	[circle, draw, label=below:$-2$] {};
\node (b) [right of=a] 	[circle, draw, label=below:$-1$] {};
\node (c) [right of=b]	[circle, draw, label=below:$-3$] {};
\node (d) [above of=b]	[below=.5cm, circle, draw, label=left:$-6$] {};
\draw (a)--(b)-- (c); 
\draw (b) -- (d);
\end{tikzpicture}
\end{center}

\section{Frobenius pullback of rational elliptic surfaces}
\mylabel{Frobenius pullback}

Let $k$ be an algebraically closed ground field of arbitrary characteristic $p>0$.
Suppose that   $f:J\ra\PP^1$ is a rational elliptic surface that is relatively minimal and jacobian.
In this section, we  examine the Frobenius pullback $J^{(p/\PP^1)}$, which is defined by the cartesian diagram
$$
\begin{CD}
J^{(p/\PP^1)}	@>>>	J\\
@VgVV			@VVfV\\
\PP^1	@>>F>		\PP^1,
\end{CD}
$$
where $F:\PP^1\ra\PP^1$ is the relative Frobenius.
Here we collect some elementary general facts, although eventually we are mainly interested in characteristic two.
First note that the singularities on $J^{(p/B)}$ are Zariski, according to Proposition \ref{frobenius pullback zariski}.
In particular, the Frobenius pullback is locally of complete intersection, hence Gorenstein,
and its tangent sheaf is locally free.
Let us write $\O_J(n)$ and  $\O_{J^{(p/\PP^1)}}(n)$ for the invertible sheaves $f^*\O_{\PP^1}(n)$ and  $g^*\O_{\PP^1}(n)$,
respectively.

\begin{proposition}
\mylabel{dualizing sheaf}
The  dualizing sheaf  is $\omega_{J^{(p/\PP^1)}}=\O_{J^{(p/\PP^1)}}(p-2)$.
\end{proposition}

\proof
By assumption, $f:J\ra\PP^1$ is a rational elliptic surface that
is relatively minimal and jacobian.
According to \cite{Cossec; Dolgachev 1989}, Proposition 5.6.1 we have $\omega_J=\O_J(-1)$.
Using $\omega_{\PP^1}=\O_{\PP^1}(-2)$
one gets $\omega_{J/\PP^1}=\O_J(1)$. Relative dualizing sheaves commute with
flat base change.
It follows that the relative dualizing sheaf 
for  $J^{(p/\PP^1)}\ra\PP^1$ is the preimage of $F^*\O_{\PP^1}(1)=\O_{\PP^1}(p)$,
and the result follows.
\qed

\begin{proposition}
\mylabel{normal reduced fibers}
The singular locus on the surface $J^{(p/\PP^1)}$ corresponds
to the locus on the surface $J$ where the fibration
$J\ra\PP^1$ is not smooth. In particular, $J^{(p/\PP^1)}$ is normal
if and only if $J\ra\PP^1$ has only reduced fibers.
\end{proposition}

\proof
Clearly,  $ J^{(p/\PP^1)}$ is smooth over each $y\in J$ at which the fiber $F=J_{f(y)}$ is smooth.
Now suppose that the local ring has  $\edim(\O_{F,y})\geq 2$.
The preimage of the fiber in $J^{(p/\PP^1)}$ is a scheme isomorphic to $F'=F\otimes_kk[T]/(T^p)$,
which has $\edim(\O_{F',y})\geq 3$. In turn, the Frobenius pullback acquires a singularity there.
\qed

\begin{proposition}
\mylabel{cohomology vanishes}
We have $H^1(J^{(p/\PP^1)}, \O_{J^{(p/\PP^1)}}) =0$.
\end{proposition}

\proof
The  Leray--Serre spectral  sequence for $g:J^{(p/\PP^1)}\ra\PP^1$  gives an exact sequence
\begin{equation}
\label{cohomology exact sequence}
H^1(\PP^1,\O_{\PP^1})\lra H^1(J^{(p/\PP^1)}, \O_{J^{(p/\PP^1)}}) \lra H^0(J^{(p/\PP^1)}, R^1g_*\O_{J^{(p/\PP^1)}}),
\end{equation}
and the term on the left vanishes. 
Forming higher direct images commutes with flat base-change.
So to understand the term on the right, we merely have to compute the coherent sheaf $R^1f_*(\O_J)$.
Choose a section $O\subset J$. Using   $\omega_J=\O_J(-1)$ and Adjunction Formula, one gets
$(O\cdot O)=-1$. The ensuing short exact sequence $0\ra\O_J\ra\O_J(O)\ra\O_{\PP^1}(-1)\ra 0$
induces an exact sequence
$$
\O_{\PP^1}(-1)\lra R^1f_*(\O_J)\lra R^1f_*(\O_J(O)).
$$
The terms on the left is invertible, the term in the middle has rank one,
and the term on the right vanishes, consequently $R^1f_*(\O_J)=\O_{\PP^1}(-1)$.
Thus $R^1g_*\O_{J^{(p/\PP^1)}}=\O_{\PP^1}(-p)$, whose group of global sections vanishes.
The assertion thus follows from the  exact sequence \eqref{cohomology exact sequence}.
\qed

\medskip
Since $J$ is a smooth rational surface, a minimal model is either some Hirzebruch surface or the projective plane.
In the latter case, $J$ factors over the blowing-up of the projective plane, by $\rho=10$.
It follows that in all cases there is some \emph{ruling} $r:J\ra\PP^1$, that is, fibrations whose generic fiber
is a projective line.  

\begin{proposition}
\mylabel{ruling degree two}
For each closed fiber $D=r^{-1}(t)$ of some ruling, the induced
map $f:D\ra\PP^1$ has degree two.
\end{proposition}

\proof
Regarding this degree as an intersection number, we easily reduce to the 
case that the fiber $D$ is smooth. Then $D=\PP^1$. Now the Adjunction Formula gives
$-2=(D\cdot\omega_S)=-(D\cdot F)$, and the assertion follows.
\qed

\medskip
Now suppose that the elliptic fibration $f:J\ra\PP^1$ has only reduced fibers, such that $J^{(p/\PP^1)}$ is normal.
Let $r:J\ra\PP^1$ be some ruling. Consider the Stein factorization $J^{(p/\PP^1)}\ra C$ given by
the commutative diagram
$$
\begin{CD}
J^{(p/\PP^1)}	@>>>	 J\\
@VVV			@VVrV\\
C		@>>>	\PP^1.
\end{CD}
$$

\begin{proposition}
\mylabel{fibrations from rulings}
The morphism $C\ra\PP^1$ is an isomorphism, and 
the generic fiber  of $J^{(p/\PP^1)}\ra C$ is a regular curve $R$ with  
$h^1(\O_R)=p-1$. If $p\geq 3$, the geometric generic fiber is a rational curve
with two singularities.
In characteristic two,  $J^{(p/\PP^1)}\ra C$ is a quasielliptic fibration.
\end{proposition}

\proof
Let $D=r^{-1}(t)$ be a smooth fiber, thus $D\simeq\PP^1$. The induced morphism $f:D\ra\PP^1$ is affine
of degree two, so $D$ is the relative spectrum of some coherent $\O_{\PP^1}$-algebra
whose underlying module is $\shA=\O_{\PP^1}\oplus\O_{\PP^1}(-1)$.
Likewise, the relative  Frobenius $F:\PP^1\ra\PP^1$ is  the relative spectrum
of some  $\shB=\O_{\PP^1}\oplus\O_{\PP^1}(-1)^{\oplus (p-1)}$.
In turn, the preimage of $D$ in $J^{(2/\PP^1)}$ is the spectrum of $\shC=\shA\otimes\shB$,
and we see $h^0(\shC)=1$ and $h^1(\shC)=p-1$.
For $p=2$, the double covering $\PP^1=D\ra\PP^1$ is ramified at two points,
so the preimage of $D$ has precisely two singularities.

In light of $h^0(\shC)=1$, the  Theorem on Formal Functions
ensures that the composite mapping $J^{(p/\PP^1)}\ra \PP^1$
equals its own Stein factorization, thus $C=\PP^1$. 
Its generic fiber is normal, because this holds for $J^{(p/\PP^1)}$, 
and its arithmetic genus equals $h^1(\shC)=p-1$.
In characteristic $p=2$, this means that the fibration is quasielliptic.
In odd characteristics, the non-smooth locus of the generic fiber for $J^{(p/\PP^1)}\ra C$ consist of two points,
because this holds for almost all closed fibers.
\qed

\section{Lang's classification and Frobenius pullback}
\mylabel{Lang classification}

For   rational  elliptic surfaces over the complex numbers, 
the possible configurations of Kodaira types were classified by
Persson \cite{Persson 1990} and Miranda \cite{Miranda 1990}.
This classification was extended by Lang \cite{Lang 2000}  to  characteristic $p=2$,
including for each case  examples of global Weierstra\ss{}  equations. These equations show that
all numerically possible cases indeed exist. This builds on the normal forms for
unstable fibers obtained in \cite{Lang 1994}, Section 2A.
Lang's results will be the key
to understand and construct simply-connected Enriques surfaces and their K3-like coverings.

He introduced  very useful short-hand notation for configurations of Kodaira types:
For example, $\II + 21^6$ stands for  the configuration comprising one fiber of type $\II$, one fiber of type $\I_2$,
and six fibers of type $\I_1$. 
It turns out that there are 147 possible  configurations of Kodaira types. Here we are only interested in
cases where all fibers are reduced, for which there are 110 configurations.
There are 35 cases where all fibers of $J$ are semistable, 58 cases with    semistable fibers and
precisely one unstable fiber, and 17 cases with only unstable fibers.
The latter is equivalent to the condition that the global $j$-invariant is $j(J/\PP^1)=0$.

Lang  used  alpha-numerical symbols for  
normal forms of  \emph{local} Weierstra\ss{}  equations for unstable fibers.
Let us call them the \emph{Lang types}. One may tabulate the Lang types for reduced fibers as follows:
$$
\begin{array}[t]
{l 			| c	  c	  c	| c	  c	| c	| c	  c	  c	| c	  c	  c	| c}
\text{Lang type}	& {\rm 1A}& {\rm 1B}& {\rm 1C}& {\rm 2A}& {\rm 2B}& {\rm 3}& {\rm 9A}	& {\rm 9B}& {\rm 9C}& {\rm 10A}	& {\rm 10B}& {\rm 10C}& {\rm 11}\\ 
\hline &&&&&&&&&&&&&\\[-2ex]
\text{Kodaira type}	& \II	& \II	& \II 	& \III	& \III 	& \IV	& \II	& \II	& \II	& \III	& \III	& \III	& \IV\\
m			& 1	& 1	& 1 	& 2	& 2 	& 3	& 1	& 1	& 1	& 2	& 2	& 2	& 3\\
v			& 4	& 6	& 7 	& 4	& 6 	& 4	& 4	& 8	& 12	& 4	& 8	& 12	& 4\\
\delta			& 2	& 4	& 5 	& 1	& 3 	& 0	& 2	& 6	& 10	& 1	& 5	& 9	& 0\\
\end{array}
$$
\vspace{3em}
\centerline{Table \stepcounter{table}\arabic{table}: \text{The reduced unstable Lang types}}
\vspace{-3.5em}

\noindent
As in Section \ref{Ogg's formula},  the integer $m\geq 1$ is the number of irreducible components in the fiber, 
$v\geq 1$ is the valuation of a minimal Weierstra\ss{}  equation,
and $\delta\geq 0 $ is the wild part of the conductor, which satisfies $\delta= v -m-1$ according
to Ogg's Formula. 

In what follows, we work over an algebraically closed ground field $k$ of characteristic $p=2$.
Fix a rational elliptic surface $J\ra \PP^1$ that is jacobian, with only reduced fibers.
Hence the Frobenius pullback $X'=J^\fpb$ is normal.
Choose a section $O\subset J$, and let $W\ra\PP^1$ be the resulting Weierstra\ss{}  fibration.

In what follows, fibers of Lang type  9C play an exceptional role.
They have Kodaira type $\II$ and numerical invariants $v=12$, $m=1$ and $\delta=10$.
If present, there are no other singular fibers, hence the contraction $J\ra W$ to
the Weierstra\ss{}  model is an isomorphism, and the whole of $J$ is
given by the global Weierstra\ss{}  equation 
\begin{equation}
\label{9C}
y^2+t^3\gamma_0 y = x^3+t\gamma_1x^2+ t\gamma_3x + t\gamma_5,
\end{equation}
according to \cite{Lang 1994}, Section 2A. Here  $\gamma_i\in k[t]$ are polynomials of degree $\leq i$, and furthermore
$t\nmid\gamma_0,\gamma_5$. We regard the base of the fibration as the projective line
$\PP^1=\Spec k[t]\cup\Spec k[t^{-1}]$. The singular fiber is located over $t=0$.

\begin{theorem}
\mylabel{nonrational singularity}
The Frobenius pullback $X'=J^\fpb$ contains a non-rational singularity $x_1'\in X'$ if and only
if $J\ra\PP^1$ contains a singular fiber of Lang type {\rm 9C}, and the coefficient
$\gamma_3\in k[t]$ in the above Weierstra\ss{}  equation is divisible by   $t\in k[t]$.
In this case, the singularity is an elliptic double point,  
the exceptional divisor on the minimal resolution
of singularities is a rational cuspidal curve with intersection matrix $N=(-1)$,
and there are no other singularities on $X'$.
\end{theorem}

\proof
We start by collecting  some general observations.
Clearly a singularity $x\in J^\fpb$ is  rational if  
the local Weierstra\ss{}  equations for the corresponding fiber in the Weierstra\ss{}  
model $W\ra\PP^1$ of $J\ra\PP^1$ remains minimal after Frobenius pullback.
This automatically holds for the semistable fibers. 
Let $J_a\subset J$ be an unstable fiber.
The invariants for the minimal model of $J^{(2/\PP^1)}$ satisfy $2v-12\lambda = 2+\delta + (m'-1)$,
where $\lambda\geq 0$ is the number of repetitions in the Tate Algorithm, as explained
in Proposition \ref{unstable fiber}. Thus
$12\lambda\leq 2v-1-\delta$. 
Using Table 1, we see that the right hand side is $<12$, and thus $\lambda =0$ and the singularities on $X'$ are rational, 
unless the fiber is of Lang type 9C or 10C,
in which case $\lambda=1$ may occur.
 
According to \cite{Lang 2000}, Lemma 0.1, we have $\sum_{a\in\PP^1}v_a=12$.
Fibers of Lang type 9C and 10C have $v_a=12$, whence $J\ra\PP^1$ has exactly one
singular fiber. Their Kodaira types are $\II$ or $\III$, whence the projection $J\ra\PP^1$
is smooth outside a single point, so $\Sing(X')$ contains exactly one point.
Let $x'_1\in X'$ be the unique singularity.
To  verify our assertion, it suffices to treat these two special cases.
 
Suppose first  that $J\ra\PP^1$ contains a fiber of Lang type  10C.
By \cite{Lang 2000}, the Weierstra\ss{}  equation is of the form
\begin{equation}
\label{10C}
y^2+t^3\gamma_0 y = x^3+t\gamma_1x^2+ t\gamma_3x + t^2\gamma_4,
\end{equation}
where   $\gamma_i\in k[t]$ are polynomials of degree $\leq i$, subject to 
$t\nmid\gamma_0,\gamma_3$.  The singular fiber is again located over $t=0$.
Going through the Tate Algorithm (\cite{Tate 1972}, Section 7) and using $t\nmid\gamma_3$, we see that the 
Frobenius pullback of the Weierstra\ss{} equation remains minimal, and
the singular fiber on $X'$ has Kodaira type $\I_n^*$ for some $n\geq 0$.
We have $v_a=2v_b=24$ and $\delta_a=\delta_b=9$, so  Ogg's Formula $v_a=2+\delta_a+(m_a-1)$
gives $m_a=14$, thus the fiber is of type $\I_9^*$, and the singularity is a rational double point
of type $D_{12}$.
%
%

Now suppose that $J\ra\PP^1$ contains a fiber of Lang type  9C 
as in \eqref{9C}. Again going through the Tate Algorithm, we see that the 
the singular fiber on $X'$ has Kodaira type $\I_n^*$ for some $n\geq 0$ if and only if $t\nmid \gamma_3$.
In this situation we  argue as in the preceding paragraph to 
deduce that the fiber is of Kodaira type $\I_8^*$, and the singularity is again a rational double point
of type $D_{12}$. For $t\mid\gamma_3$, the Weierstra\ss{} equation is not minimal.
This establishes the equivalence of the two conditions. 
%
%
%

Now suppose that there is a non-rational singularity $x'_1\in X'$, hence $J\ra\PP^1$ contains
a fiber of Lang type  9C. The local Weierstra\ss{} equation for $X'$ must have $v'=2\cdot 12-12$,
whence is almost minimal, and the singularity is elliptic, by Proposition \ref{almost minimal}.
The minimal model of $X'$ has   $v'=2\cdot v-12=12$. Ogg's Formula $12=2+(m'-1)+10$ yields
$m'=1$, so the Kodaira type is $\II$.
The statement on the exceptional divisor now follows from Corollary \ref{almost minimal II}.
\qed

\medskip
We now turn to the global properties:
Let $r':S\ra X'$ be the minimal  resolution of singularities, and write
$h:S\ra \PP^1$ for the induced elliptic fibration.

\begin{lemma}
The smooth surface $S$ is either a K3 surface or a rational surface.
The latter holds if and only if $X'$ contains a non-rational singularity.
\end{lemma}

\proof
According to Proposition \ref{dualizing sheaf} and Proposition \ref{cohomology vanishes}, 
we have $\omega_{X'}=\O_{X'}$ and $h^1(\O_{X'})=0$. 
Now suppose that $X'$ has only rational singularities. For the minimal resolution of singularities,
this means $\omega_S=\O_S$ and $h^1(\O_S)=0$. Furthermore, the rational elliptic surface $J$
has Betti number $b_2(J)=10$, thus $b_2(S)\geq 11$. From the classification of surfaces we infer
that $S$ is a K3 surface.

Now suppose that $X'$ contains a non-rational singularity $x'_1\in X'$.
Then  $K_{S/X'} = -D'$ for some negative-definite curve $D'\subset S$ contracted
by $r:S\ra X'$. Consequently the plurigenera $P_m(S)=h^0(\omega_S^{\otimes n})$ vanish for all $n\geq 1$. 
In particular, we have $h^2(\O_S)=h^0(\omega_S)=0$, so the Picard scheme $\Pic_{S/k}$ is smooth.
Now consider the  Albanese map $S\ra \Alb_{S/k}$.
This map factors over $X'$, because  the exceptional curve for the resolution of singularities $r:S\ra X'$ 
are rational curves by Proposition \ref{components rational}. 
Since $X'\ra J$ is a universal homeomorphism, and $J$ is covered by rational
curves, the Albanese map must be zero, and we infer $h^1(\O_S)=0$.
The Castelnuovo Criterion  now ensures that the surface $S$ is rational.
\qed

\medskip
Thus we are precisely in the  situation studied in Section \ref{Trivial tangent sheaf},
and determine which of the seven condition \eref{CI}--\eref{Cuspidal} from Section \ref{Trivial tangent sheaf} do hold:

\begin{proposition}
\mylabel{conditions fpb}
The elliptic fibrations $h':S\ra\PP^1$ and    $r':S\ra X'$ satisfy  conditions
\eref{CI}--\eref{Cuspidal}, except  condition \eref{Elliptic},
and the  tangent sheaf is given by    $\Theta_{X'}=\O_{X'}(2)\oplus \O_{X'}(-2)$.
\end{proposition}

\proof
According to Proposition \ref{fibrations from rulings}, each ruling  $r:J\ra\PP^1$ induce
a quasielliptic fibration $X'\ra\PP^1$, hence \eref{Cuspidal} holds.
The singularities on the Frobenius pullback $X'$ are Zariski singularities by Lemma \ref{frobenius pullback zariski},
so condition \eref{Zariski} and in particular  \eref{CI} hold.

By assumption, $J\ra\PP^1$ admits a section. Thus the same holds for $h':S\ra X'$,
hence the latter admits no multiple fiber.
If $S$ is rational, then 
by Theorem \ref{nonrational singularity} the surface contains a single singularity, which is
an elliptic singularity, so $R^1r'_*(\O_S)$ has length one.
In other words, condition \eref{Length} holds.
Furthermore, we have $\omega_S=\O_S(-C')$ for some negative-definite curve
supported on the exceptional divisor for $r':S\ra X'$, so condition \eref{Nonfixed} holds.

It remains to check condition \eref{Tjurina} on the total Tjurina number.
Let $v_b,m_b,\delta_b$ be the numerical invariants for
$J\ra\PP^1$, and $v_a,m_a,\delta_a$ be the ensuing invariants
for the minimal model $J'\ra\PP^1$ of the Frobenius pullback $X'$. Recall that $\delta_a=\delta_b$.
Since $J$ is a rational surface, we have $\sum v_b=12$, according to 
\cite{Lang 2000}, Lemma 0.1.
Suppose first that $X'$ has only rational singularities.
Then $v_a=2v_b$, and the total number of exceptional divisors for the 
resolution of singularities $J'\ra X'$ is $\sum (m_a-m_b)$.
In light of Proposition \ref{tjurina rdp}, we have to verify
$\sum(m_a-m_b)=12$.
If the fiber $J_b$ is semistable, the same holds for $J_a$, and  Ogg's Formula gives
$m_a-m_b = v_a - v_b =v_b$.
If the fiber $J_b$ is unstable, also $J_a$ is unstable, and likewise get
$$
m_a-m_b = (v_a -1 - \delta_a) - (v_b -1 - \delta_b) = 2v_b-v_b = v_b.
$$
In any case, we obtain $\sum(m_a-m_b)=\sum v_b=12$.

Now suppose that $X'$ contains an elliptic singularity.
Then $J\ra\PP^1$ is given by the Weierstra\ss{}  equation \eqref{9C},
and its Frobenius pullback is described by 
$$
y^2+t^6\gamma_0^2 y = x^3+t^2\gamma_1^2x^2+ t^2\gamma_3^2x + t^2\gamma_5^2,
$$
where $t\nmid\gamma_0$.
The Tjurina ideal   is generated by the given  Weiererstra\ss\ polynomial $y^2 +t^6\gamma_0^2 y-\ldots$, together with
the partial derivatives $t^6\gamma_0^2$ and $x^2+t^2\gamma_3^2$.
Using that  $\gamma_0$ is a unit in $k[[t]]$, and regarding the generators as integral equations in 
$t$, $x$ and $y$ of degree $d=6,2,2$, we infer that 
$$
k[[t,x,y]]/( t^6, x^2+t^2\gamma_3^2,y^2 +t^6\gamma_0^2 y-\ldots)
$$
has length $\tau =6\cdot 2\cdot 2=24$. 

Finally, we compute the tangent sheaf. Choose a ruling $J\ra \PP^1$, and let
$X'\ra\PP^1$ be the resulting quasielliptic fibration.
For each smooth fiber $D\subset J$, the preimage $C\subset X'$ has degree two
with respect to the elliptic fibration $X'\ra\PP^1$, according to Proposition \ref{ruling degree two}.
Now Theorem \ref{conditions for theta} gives $\Theta_{X'/k}=\O_{X'}(-2)\oplus\O_{X'}(2)$.
\qed

\section{K3-like coverings with elliptic double point}
\mylabel{K3-like with edp}

In this section, we shall construct   simply-connected Enriques surfaces $Y$
whose K3-like covering $X$ is normal and contains an elliptic double point.
The idea is simple: We start with a rational elliptic surface $J\ra\PP^1$
containing a singular fiber of Lang type  9C, chosen
so that the  Frobenius pullback $X'=J^\fpb$ stays rational. 
It comes with an induced elliptic fibration $X'\ra\PP^1$.
Its relatively minimal smooth model $J'$ can be seen as  a flip $X'\leftarrow S\rightarrow J'$,
and the desired K3-like covering  arises as a  flop $X'\leftarrow S\rightarrow X$.
The crucial step is the construction of the flop $X$, starting from the flip $J'$,
so that the flop   removes all quasielliptic fibrations, introduce another elliptic fibration,
and yet does not change the nature of the Zariski singularities.

For this, we need as key ingredients Shioda's theory of  Mordell--Weil lattices
\cite{Shioda 1990}  and their classification for rational elliptic surface, which is due 
to Oguiso and Shioda \cite{Oguiso; Shioda 1991}.
We start by recalling this. Let  
$J\ra B$ be an  elliptic surface over a smooth proper curve $B$, with smooth total space
and endowed  with zero section $0\subset J$. We also assume that it is relatively minimal.
The N\'eron--Severi group $\NS(J)=\Pic(J)/\Pic^0(J)$ is endowed with
the intersection pairing $(P\cdot Q)$.
The group of sections for $J\ra B$ is called the \emph{Mordell--Weil  group} $\MW(J/B)=\Pic^0(J_\eta)$.
The \emph{trivial sublattice}
$T=T(J/B)$ inside the  N\'eron--Severi group $\NS(J)$ 
is the subgroup generated by the zero-section $O\subset J$, together
with all irreducible components $\Theta\subset J_b$, $b\in B$ that are disjoint from the zero-section.
Thus we may identify $\MW(J/B)$,  up to its torsion subgroup, as a subgroup of
the orthogonal complement $T^\perp\subset\NS(J)\otimes\QQ$. 
Consequently, it  acquires a non-degenerate bilinear $\QQ$-valued form 
$$
\langle P,Q\rangle=-(\bar{P}\cdot\bar{Q})
$$ 
called the \emph{height pairing},
which was extensively studied by Shioda \cite{Shioda 1990}. The sign is   a customary convention,
introduced to have a positive-definite height pairing,
and $P\mapsto \bar{P}$ denotes   the  orthogonal projection onto $T^\perp$.
One also calls $\MW(J/B)$ the \emph{Mordell--Weil lattice}.
The subgroup comprising all sections $P\subset J$ that pass through the
same fiber components as the zero section $O\subset J$ is called the \emph{narrow Mordell--Weil lattice}.

The \emph{explicit formula} for  the height pairing is
$$
\langle P,Q\rangle = \chi(\O_J) + (P\cdot O) + (Q\cdot O) - (P\cdot Q) - \sum\contr_b(P,Q),
$$
according to \cite{Shioda 1990}, Theorem 8.6.
The summands in the middle are the usual intersection numbers, 
and the \emph{local contributions} $\contr_b(P,Q)\in\QQ_{\geq 0}$
are certain  rational numbers tabulated in \cite{Shioda 1990}, page 229. They reflect how the sections $P,Q$ pass through
the  irreducible components of the fiber $J_b$, in relation to the zero-section $O\subset J$. 
The corresponding quadratic form can be written as
$$
\langle P,P\rangle = 2\chi(\O_J) + 2(P\cdot O) - \sum\contr_b(P),
$$
where $\contr_b(P)=\contr_b(P,P)$. 

For each point $b\in B$, we denote by $\Phi_b$ the group of components in the N\'eron model
of $J\otimes\O_{B,b}$. Its order is the number of irreducible components in $S_b$ with multiplicity
$m=1$ in the fiber. Clearly,  the narrow Mordell--Weil lattice is the subgroup of $\MW(J/B)$
of elements $P$ with $P\cong O$ in $\Phi_b$ for all $b\in B$.
Note that if $P\cong O$ in $\Phi_b$, then the local contributions $\contr_b(P,Q)$
in the explicit formula for the height pairing vanish.

Suppose now that $J\ra\PP^1$ be a jacobian rational elliptic surface.
Then we have $\chi(\O_J)=1$.
Oguiso and Shioda \cite{Oguiso; Shioda 1991} tabulated the possible isomorphism classes of   Mordell--Weil
lattices $\MW(J/\PP^1)$ in terms of the trivial sublattices $T=T(J/\PP^1)$, a list containing 74 cases.
The  arguments are purely lattice-theoretical, and apply to  arbitrary characteristics.

Let us give an example: If the configuration of Kodaira types is $\III+ 32^21$, the trivial lattice becomes
$A_1^{\oplus 3}\oplus A_2$, which appears as No.\ 23 in the Oguiso--Shioda list; we then
read off that the  Mordell--Weil lattice is given by
$$
\MW(J/\PP^1) = A_1^\vee \oplus \frac{1}{6}\begin{pmatrix}2 & 1\\1&2\end{pmatrix}.
$$
Here the second summand is a Gram matrix, and the first summand is the dual $A_1^\vee=(1/2)$ for
the ADE-lattice $A_1=(2)$. 
The following observation will be useful:

\begin{proposition}
\mylabel{disjoint section}
Suppose that $J\ra\PP^1$ has only reduced fibers. Then 
there is  a section $P\subset J$ that is disjoint from the zero-section $O\subset J$
and has $P\not\cong O$ in the group of components $\Phi_b$ for some point $b\in \PP^1$.
Furthermore, every torsion section $P\neq O$ is disjoint from the zero-section.
\end{proposition}

\proof
According to \cite{Oguiso; Shioda 1991}, Theorem 2.5, the Mordell--Weil group is generated by   sections
that are disjoint from the zero-section.
So we first check that the Mordell--Weil group is non-zero. According to the Oguiso--Shioda list, only No.\ 62 
has $\MW(J/\PP^1)=0$. In this case, the vertical lattice is $T=E_8$. 
It follows that $J\ra\PP^1$ has a fiber of Kodaira type $\II^*$, which is nonreduced, contradiction.
Going through the Oguiso--Shioda list, one observes that in all other cases the
narrow Mordell--Weil lattice is a proper sublattice, so the first assertion follows.
The last  statement is \cite{Oguiso; Shioda 1991}, Proposition 5.4.
\qed

\medskip
In what follows, $k$ denotes an algebraically closed ground field of characteristic $p=2$.
Choose a Weierstra\ss{}  equation of Lang type 9C, as given in  \eqref{9C}, with $t\mid\gamma_3$.
The resulting rational elliptic surface $J\ra\PP^1$ then contains only one
singular fiber, which is located over the origin $b\in\PP^1$ and of Lang type  9C, and the
Frobenius pullback $X'=J^\fpb$ is a normal rational surface containing a unique
singularity, which is an elliptic double point obtained by contracting 
a cuspidal $(-1)$-curve, and is located over the point $a\in \PP^1$
corresponding to $b$. 
Note that the exceptional divisor on the minimal resolution of singularities
is a rational cuspidal curve with self-intersection $-1$.

Let $r':S\ra X'$ be the non-minimal resolution obtained by blowing-up further a closed regular point
on the rational cuspidal curve. Write $h:S\ra\PP^1$ for the induced elliptic fibration.
The fiber    $S_a=h^{-1}(a)$ has three irreducible components $D_1,D_2,D_3$,
where $D_2$ is the rational cuspidal curve and $D_1,D_3$ are smooth rational curves.
The dual graph, with the selfintersection numbers,  is 
\begin{center}
\begin{tikzpicture}
\tikzstyle{vertex}=[circle, draw, inner sep=0pt, minimum size=2ex]
\node[vertex] 	(a) at (2,1) 	[label=below:{$-1$}, label=above:{$D_1$}] 	{};
\node[vertex] 	(b) at (4,1) 	[label=below:{$-2$}, label=above:{$D_2$}] 	{};
\node[vertex] 	(c) at (6,1) 	[label=below:{$-1$}, label=above:{$D_3$}] 	{};
\node 		(d) at (10,1) 	{with $D_2$ rational cuspidal.};
\draw (a) -- (b) --  (c);
\end{tikzpicture}
\end{center}

\begin{proposition}
\mylabel{picard number non-minimal}
The smooth rational surface $S$ has Picard number $\rho(S)=12$.
\end{proposition}

\proof
The elliptic fibration $h:S\ra\PP^1$ is relatively minimal, except at the fiber $S_a=D_1+D_2+D_3$.
The curves $D_1,D_3$ are  disjoint $(-1)$-curves.
Contracting them yields the relatively minimal model. 
The latter has Picard number ten, by  \cite{Cossec; Dolgachev 1989}, Proposition 5.6.1,
so our surface $S$ has Picard number twelve.
\qed

\medskip
We choose the indices for the irreducible components $S_a=D_1+D_2+D_3$  so that $r':S\ra X'$ contracts $E'=D_2+D_3$,
and that $S\ra J'$ contracts the $(-1)$-curves $D_1+D_3$. Since $J'$ is a rational elliptic surface,
we have $K_{J'}=-J'_a$ whence $K_S=-D_2$.
Now let  $r:S\ra X$ be the contraction of $E=D_1+D_2$.
This yields another  normal surface $X$ with $\omega_X=\O_X$, containing
a unique elliptic singularity obtained by contracting a cuspidal $(-1)$-curve.

\begin{theorem}
\mylabel{existence K3-like II}
The normal surface $X$, with its  elliptic singularity obtained by contracting a cuspidal $(-1)$-curve,  
is a  K3-like covering.
\end{theorem}

\proof
In light of Theorem \ref{conditions sufficient}, we have to verify that
the elliptic fibration $h:S\ra\PP^1$ and the contraction $r:S\ra X$ of the exceptional divisor $E=D_1+D_2$ satisfies
conditions \eref{CI}--\eref{Zariski} of Section \ref{Trivial tangent sheaf}.
Recall that $x_1\in X$ and $x'_1\in X'$ are the unique singular points.
It follows from  Lemma \ref{isomorphic singularities} below 
that the complete local rings $\O_{X,x_1}^\wedge$ and $\O^\wedge_{X',x'_1}$ are isomorphic.
Thus all conditions but \eref{Elliptic}  follow from 
Proposition \ref{conditions fpb}.

It remains to show that condition \eref{Elliptic} holds, that is, we must construct on  $X$ 
another elliptic fibration.
For this we first consider the flip $J'$, which is a smooth rational elliptic surface,
whose sole singular fiber is of type $\II$. Choose a zero-section $O\subset J'$,
and consider the resulting Mordell--Weil lattice $\MW(J'/\PP^1)$.
Since the  trivial lattice is $T=0$,   the Mordell--Weil lattice
is $E_8$, according to \cite{Oguiso; Shioda 1991}. For each section $P\subset J'$ corresponding to a root $\alpha\in E_8$,
the height pairing gives
$2=\langle P,P\rangle = 2 +2 (P\cdot 0)$, whence $P$ is disjoint from the zero-section.
Now let $\alpha,\beta\in E_8$ be two adjacent simple roots, and let
$P,P'\subset J'$ be the sections corresponding to the roots $\alpha,\alpha+\beta\in E_8$,
such that $\alpha\cdot(\alpha+\beta)=2-1=1$. 
In terms of the height pairing, this means
$$
1=\langle P,P'\rangle = 1 +(P\cdot O) + (P'\cdot O) - (P\cdot P') = 1-(P\cdot P'),
$$
thus $O,P,P'\subset J'$ are pairwise disjoint. From this we infer that
at least one of the resulting three strict transforms on $S$
is a section for $h:S\ra\PP^1$   passing through the irreducible component $D_2\subset S_a=h^{-1}(a)$.
After changing the zero-section $O\subset J$, we may assume that the strict transform $Q\subset S$
of the zero-section satisfies $(Q\cdot D_2)=1$.
The dual graph takes the form
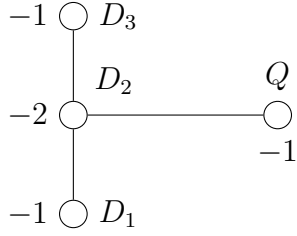
\begin{center}
\begin{tikzpicture}
[node distance=2cm]
\tikzstyle{vertex}=[circle, draw, inner sep=0mm, minimum size=2ex]
\node[vertex] (b) at (2,3) 	[label=left:{$-2$}, label=45:{$D_2$}] 		{};
\node[vertex] (a) [below of=b] 	[above=.5cm, label=left:{$-1$}, label=right:{$D_1$}] 	{};
\node[vertex] (c) [above of=b] 	[below=.5cm, label=left:{$-1$}, label=right:{$D_3$}] 	{};
\node[vertex] (d) [right of=b] 	[right=.5cm, label=below:{$-1$}, label=above:{$Q$}] 	{};
\draw (a) -- (b) --  (c);
\draw (b) -- (d);
\node (d) [right of=d]		 [right=.5cm] {with $D_2$   rational cuspidal.};
\end{tikzpicture}
\end{center}

Consider the negative-semidefinite curve  $C=D_1+D_2+ Q\subset S$ of arithmetic genus $h^1(\O_C)=1$.
Then  some $\O_S(nC)$, with $n\geq 1$
induces another genus-one fibration $g:S\ra\PP^1$.

We claim that the curve $C$ itself is not movable. Seeking a contradiction,
we assume it were. Then the curve $C\subset S$ is \emph{d-semistable},
which means that the normal sheaf is $\O_C(-C)=\O_C$, such that
$\O_{D_2}(-D_2-D_1-Q)=\O_{D_2}$. In light of the   fibration
$h:S\ra\PP^1$ defined by the movable curve $D_1+D_2+D_3$, we have $\O_{D_2}(-D_2)=\O_{D_2}(D_1+D_3)$.
It follows that $\O_{D_2}(D_3-Q)=\O_{D_2}$. By Riemann--Roch, this does not
hold for the two points $D_3\cap D_2\neq Q\cap D_2$ on the  rational cuspidal curve $D_2$.
We conclude that $C$ is not movable.

However, the curve $2C $ is movable. To see this, set $\shL=\O_S(C)$. Since $\Pic^0(C)=k$ is 2-torsion  and
the sheaf $\shL_C$ is numerically trivial, we must have $\shL^{\otimes 2}_C=\O_C$. Consider the short exact sequence
$0\ra\shL\ra\shL^{\otimes 2}\ra \O_C\ra 0$
and the ensuing long exact sequence
$$
H^0(S,\shL^{\otimes 2})\lra H^0(C,\O_C)\lra H^1(S,\shL).
$$
It thus suffices to check that the term on the right vanishes.
By Serre Duality,  $h^2(\shL)=h^0(\omega_S\otimes\shL^\vee)=0$,
because $K_S=-D_2$. Riemann--Roch gives
$$
1-h^1(\shL)=\chi(\shL)= C\cdot (C + D_2)/2 + \chi(\O_S) = 1,
$$
thus $h^1(\shL)=0$.
Summing up, $2C\subset S$ is movable, and thus defines a genus-one fibration
$g:S\ra\PP^1$. The intersection number
$$
C\cdot (D_1+D_2+D_3) = C\cdot D_3 = D_2\cdot D_3= 1
$$
shows that $2C$ and thus also all other   fibers of $g:S\ra\PP^1$ have degree two with
respect to the original elliptic fiber $h:S\ra\PP^1$.

According to Proposition \ref{picard number non-minimal}, the rational surface $S$ has Picard number $\rho(S)=12$.
By construction, the irreducible components $D_1,Q\subset C$ are disjoint $(-1)$-curves on $S$.
Contracting them thus gives the relatively minimal model for the genus-one fibration $g:S\ra\PP^1$.
We conclude that all $g$-fibers but $C\subset S$ are minimal, contain at most two irreducible components,
and are thus of Kodaira type $\II$, $\III$ or $\I_n$ with $n\geq 2$.

It remains to check that our genus-one fibration is elliptic.
Seeking a contradiction, we assume that $g:S\ra\PP^1$ is quasielliptic.
Now we use the fact that the fibers of $h:S\ra\PP^1$, except the sole singular fiber $h^{-1}(a)=D_1+D_2+D_3$,   
contain no $(-2)$-curves.
It follows that every component of every fiber of $g:S\ra\PP^1$   is horizontal with
respect to $h:S\ra\PP^1$, with the exception of $D_1$ and $D_2$. Moreover,  every fiber other than $2C$
has Kodaira type $\II$ or $\III$, and its irreducible components  necessarily pass through $D_3\subset  h^{-1}(a)$,
because $D_1$ and $D_2$ are contained in the fiber $2C$.

Since $S$ is rational, all other fibers beside $2C$ are non-multiple (\cite{Cossec; Dolgachev 1989}, Proposition 5.6.1),
whence there type coincides with the type of the corresponding fiber on the
associated jacobian quasielliptic fibration.
According to  Ito's analysis of jacobian quasielliptic fibrations in characteristic two 
(\cite{Ito 1994}, Proposition 5.1), there is   at least one reducible fiber, which in our situation
must be of type $\III$. In fact, the Mordell--Weil lattice for a quasielliptic fibration is 2-torsion,
whence the trivial lattice has rank ten, so there are eight fibers of Kodaira type $\III$.
Each irreducible component $R,R'\subset S$ of such a  fiber is a section for $h:S\ra\PP^1$,
both pass through $D_3\subset S$, and they meet tangentially somewhere.
Consequently, their  images define two sections $P,P'\subset J'$ with
intersection matrix $(\begin{smallmatrix}-1&3\\3&-1\end{smallmatrix})$.
Obviously, $P,P'$ are disjoint from the zero-section $O\subset J'$, which is the image of $Q\subset S$.
To derive a contradiction, we compute the Gram matrix for the height paring on $J'$.
As noted above, the sections $P,P'\neq O$ are  disjoint from the zero-section.
Thus $\langle P,P\rangle = 2 + 2(P\cdot O) = 2$,   likewise $\langle P',P'\rangle =2$, and finally
$$
\langle P,P'\rangle = 1 + (P\cdot O)+(P'\cdot O) - (P\cdot P')= -2.
$$
So the Gram matrix with respect to the height pairing is $(\begin{smallmatrix}2&-2\\-2&2\end{smallmatrix})$,
which has zero determinant. On the other hand, the Mordell--Weil lattice is non-degenerate,
contradiction. Thus $g:S\ra\PP^1$ and the induced fibration on $X$ are elliptic  
rather than quasielliptic fibrations.
This means that the  condition \eref{Elliptic} indeed holds.
\qed

\medskip
In the preceding proof, we have used the following observation:

\begin{lemma}
\mylabel{isomorphic singularities}
Let $R$ be a discrete valuation ring, and $J\ra\Spec(R)$ a relatively minimal
genus-one fibration, endowed with two disjoint sections $P_1,P_2\subset J$.
Let $S\ra J$ be the blowing-up whose center consists of the two closed points
$z_i\in P_i$, and let $E_i\subset S$ be their preimage.
Write  $S\ra X_i$ for the  blowing-down so that $E_i$ surjects onto the closed fiber of $X_i$,
and let $x_i\in X_i$ be the resulting singular point. Then the two local rings
$\O_{X_1,x_1}$ and $\O_{X_2,x_2}$ are   isomorphic.
\end{lemma}

\proof
Let $f:J_\eta\ra J_\eta$ be the isomorphism that sends the the generic fiber of $P_1$
to the generic fiber of $P_2$. By relative minimality, it extends to 
an isomorphism $f:J\ra J$. It has $f(z_1)=z_2$, and by the universal property
of blowing-ups induced an isomorphism $f:S\ra S$ with $f(E_1)=E_2$.
In turn, it maps the exceptional locus for $S\ra X_1$ to the exceptional locus for $S\ra X_2$.
By the universal property of contractions \cite{EGA II}, Lemma 8.11.1,
it yields the   isomorphism between the two local rings in question.
\qed

\medskip
As the referee   pointed out, the resulting simply-connected Enriques surfaces $Y=X/G$
must be supersingular, rather then classical.
Indeed, the elliptic fibration $J\ra\PP^1$ has constant $j$-invariant zero,  according to Lang's classification
\cite{Lang 2000}, and all fibers but $J_b$ are smooth.
In turn, this also holds for the induced fibrations on $J'$, $S$ and $X$.
If $Y$ would be classical, the induced elliptic fibration  $\varphi:Y\ra\PP^1$ would have 
two multiple fibers that are tame, and one of which has   the supersingular elliptic curve
$F\subset Y$ as  reduction. This implies $\O_F(F)\simeq\O_F$, and the fiber must be wild, contradiction.

\section{Uniqueness of the elliptic double point}
\mylabel{Uniqueness edp}

Let $Y$ be a simply-connected Enriques surface over an algebraically closed
ground field $k$ of characteristic $p=2$, with normal
K3-like covering. The goal of this section is to establish the
following uniqueness result:

\begin{theorem}
\mylabel{unique elliptic singularity}
If there is an   elliptic singularity $x_1\in X$, then
there are no other singularities on $X$.
\end{theorem}

\proof
Seeking a contradiction, we assume that there is another singularity $x_2\in X$.
Let $y_1,y_2\in Y$ be the images of the two singularities $x_1,x_2\in X$.

First, we show that there are no elliptic fibrations $\varphi:Y\ra\PP^1$ admitting
 radical two-sections $A\subset Y$. Seeking a contradiction, we assume that such a fibration exists.
Let $J\ra\PP^1$ be its jacobian fibration, which is a rational elliptic surface,
and write $X'=J^\fpb$ for the Frobenius pullback.
By Proposition \ref{two-section = section}, 
the normal surfaces $X'$ and $X$ are birational. It follows that the surface $X'$ is rational,
whence contains an elliptic singularity. According to Theorem \ref{nonrational singularity},
$J\ra\PP^1$ must contain a fiber of Lang type 9C, and this is the only singular fiber $J_b$.
In turn, the Frobenius base-change $X'=J^\fpb$ has but one singularity, lying on the fiber $X'_a$,
where $a\in\PP^1$ corresponds to $b\in\PP^1$ under the relative Frobenius morphism.

To proceed, recall that if $2F\subset Y$ is a multiple fiber with $F$ smooth, the K3-like covering $X$
is smooth along $\epsilon^{-1}(F)$, according to Proposition \ref{preimage elliptic}.
If the fiber $Y_b$ is not a multiple fiber, then all half-fibers must be smooth,
and by Proposition \ref{formal isomorphism non-multiple} there is but one singularity on $X$, namely the point 
corresponds to the singular point in $Y_b$,
contradiction. Thus $Y_b$ is a multiple fiber. According to Proposition \ref{fibers jacobian}, the
two fibers $Y_b$, $J_b$   have the same Kodaira symbol $\II$.
Write $ Y_b=2F$, with  a rational curve  $F$.
By the above discussion, the two singularities 
$x_1,x_2\in X$ map to $F$.
This analysis applies not only to $\varphi:Y\ra\PP^1$, but also to all other elliptic fibrations
admitting a radical two-section.

According to Theorem \ref{properties genus-one fibrations}, there is another elliptic fibration $\varphi':Y\ra\PP^1$
with $(F\cdot F')=1$. Then the rational curve $F$ is a radical two-section for the new fibration $\varphi'$,
according to Proposition \ref{radical = rational}.
In turn, we have $x_1,x_2\in F'$.
Now $x_1,x_2\in F\cap F'$ contradict the intersection number $(F\cdot F')=1$.
Consequently, there is no elliptic fibration on $Y$ admitting a radical two-section.
According to Proposition \ref{radical two-section}, there are no $(-2)$-curves on $Y$.

Now choose again an elliptic fibration $\varphi:Y\ra\PP^1$.
Then all fibers are irreducible, according to Proposition \ref{radical two-section}.
By the result of Liu, Lorenzini and Raynaud \cite{Liu; Lorenzini; Raynaud 2005}, Theorem 6.6
the same holds for the associated jacobian fibration $J\ra\PP^1$.
Suppose  for all  elliptic fibrations on $Y$ the jacobian would have only one singular fiber $J_b$.
Arguing as above, this fiber $Y_b=2F$ is multiple, and we find two orthogonal elliptic fibrations with
$(F\cdot F')=1$ and thus $y_1,y_2\in F\cap F'$, contradiction.
Consequently, we may choose an elliptic fibration $\varphi:Y\ra\PP^1$ so
that its jacobian $J\ra\PP^1$ contains at least two singular fibers.

Since all fibers are irreducible, the  trivial lattice $T\subset \NS(J)$ is zero, and
the Mordell--Weil group $\MW(J/\PP^1)$ is isomorphic to $E_8$.
It follows that the Picard group of the generic fiber $Y_\eta$ is a free abelian group of rank nine.
Since $\epsilon:X\ra Y$ is a universal homeomorphism, the same holds for
the generic fiber $X_\eta$. Now let $S\ra X$ be the resolution of singularities,
and $S\ra S'$ be the relative minimal model over $\PP^1$.
Then $S'\ra\PP^1$ is a rational elliptic surface, which is not necessarily jacobian.
In turn, $\Pic(S')$ has rank ten.
However, we just saw $\Pic(S'_\eta/\eta)$ has rank nine, whence
all fibers of $S'\ra\PP^1$ are irreducible.
Again by the result of Liu, Lorenzini and Raynaud \cite{Liu; Lorenzini; Raynaud 2005}, Theorem 6.6
the same holds for the smooth model $J'\ra\PP^1$ of the Frobenius pullback $X'=J^\fpb$.
If a singular fiber $J_a'$ is semistable, then also  the corresponding
fiber $J_b$ is semistable, and $m_a=2m_b\geq 2$, contradiction. Thus all singular fibers $J'_a$
are of Kodaira type $\II$, with numerical invariant $m_a=1$.
According to Theorem \ref{nonrational singularity}, the presence of two 
singular fibers on $J$ ensures that all singularities on the Frobenius pullback $X'=J^\fpb$
are rational. By Lemma \ref{unstable fiber}, each singular fiber $J'_a$ has numerical invariant
$m_a=2m_b+\delta +1 \geq 3$, contradiction.
\qed

\section{K3-like coverings with rational double points}
\mylabel{K3-like with rdp}

Let $k$ be an algebraically closed ground field of characteristic $p=2$.
In this section, we construct K3-like coverings $X$ that are normal with
only rational singularities.
The main result is that for each but possibly six of the 110 configurations
of Kodaira types from Lang's classification  \cite{Lang 2000} of  rational elliptic surfaces
with reduced fibers
actually appears as a configuration on some K3-like covering $X$ and
the ensuing Enriques surfaces $Y=X/G$.

Indeed, we start with   a rational elliptic surface $J\ra\PP^1$ that is jacobian with
only reduced fibers and assume that the Frobenius pullback $X'=J^\fpb$
has only rational singularities. In other words, the minimal resolution
of singularities $r':S\ra X'$ is a K3 surface, endowed with an elliptic fibration $S\ra\PP^1$.
The  case leading to an elliptic singularity was already
treated in Section \ref{K3-like with edp}.

We write  $E'\subset S$ for the exceptional divisor of the minimal resolution of singularities   $r':S\ra X'$.
Let $E\subset S$ be another vertical negative-definite curve, and $r:S\ra X$ be its contraction. We now seek
to choose this $E$ so that the resulting flop
$X'\leftarrow S\rightarrow X$ yields a K3-like covering $X$.
Roughly speaking, the new contraction must 
introduce another elliptic fibration, destroy all quasielliptic fibrations,
yet keep the Kodaira type of the singular fibers.
To achieve these conflicting tasks, we introduce for each closed point $a\in \PP^1$ the following
conditions relating the vertical negative-definite curves 
$E_a$ and $ E'_a$ on the K3-surface $S$:
\newcounter{MNr}
\setcounter{MNr}{0}
\newcommand{\Mu}[1]{\refstepcounter{MNr}\textbf{(M\arabic{MNr})} \label{#1} }
\newcommand{\mref}[1] {\text{\rm (M\ref{#1})}}

\medskip
\begin{list}{-mm}{\leftmargin2em\itemsep1em}
\item[\Mu{Semistable}] 
\emph{If the fibers $X'_a$ and $S_a$  are semistable, then $E_a$ either coincides   with $E'_a$, or it is 
the strict transform of $X'_a$.}

\item[\Mu{II}] 
\emph{If the fiber  $X'_a$ has type $\II$  and $S_a$  has type $\I_n^*$ for some $n\geq 0$, then either $E_a$   coincides
with $E_a'$, or it  is the union
of all irreducible components $\Theta\subset S_a$ except the terminal component
nearest to  the strict transform of $X'_a$.}

\item[\Mu{III}] 
\emph{If the fiber  $X'_a$ has type $\III$  and $S_a$  has type $\I_n^*$ for some $n\geq 0$, then either $E_a$   coincides
with $E_a'$, or it is the union
of all irreducible components $\Theta\subset S_a$ except the two   terminal components   contained in $E_a'$. }
%

\item[\Mu{Other}] 
\emph{In all other cases, $E_a$ coincides with $E'_a$.}
\end{list}

\medskip
By abuse of notation, we here say that the fiber $X_a$ on the normal surface $X$
has a certain Kodaira type if this is the Kodaira type of $J_a$ on the smooth surface $J$.
Moreover, a \emph{terminal irreducible component} on $S_a$ means a component corresponding to 
a   vertex in the dual graph that is adjacent to only one other vertex.
Note that in the boundary case $n=0$ of condition \mref{II},  the curve $E_a$ is the union of the irreducible
components $\Theta\subset S_a$ except an arbitrary terminal component.  

\begin{definition}
\mylabel{mutations}
In the preceding situation, we say that $X'\leftarrow S\rightarrow X$ is a \emph{mutation}
if the conditions \mref{Semistable}--\mref{Other} on the vertical negative-definite curves  
$E,E'\subset S$ holds for all closed points $a\in\PP^1$.
A mutation is a \emph{good mutation} if  in addition
there is a horizontal Cartier divisor $F\subset X$ that is an elliptic curve.
\end{definition}

Good mutations lead to the desired  flops:

\begin{proposition}
\mylabel{mutations K3-like coverings}
Suppose  $X'\leftarrow S\ra X$ is a good mutation of the Frobenius pullback $X'=J^\fpb$. 
Then the normal surface $X$ is a K3-like covering with only rational singularities.
The induced elliptic fibration $\varphi:Y\ra\PP^1$ on the simply-connected
Enriques surface $Y=X/G$ has the same configuration of Kodaira types
as the original rational elliptic surface $J\ra\PP^1$.
Furthermore, the induced fibration $f:X\ra\PP^1$ induces an inclusion 
$H^0(X,\Theta_{X/k})\subset H^0(\PP^1,\Theta_{\PP^1/k})$.
\end{proposition}

\proof
In our situation, $S$ is a K3-surface.
In order to apply Theorem \ref{conditions sufficient}, it suffices 
to check that conditions \eref{Reduced}, \eref{Tjurina}, \eref{Elliptic}
and \eref{Zariski}
from Section \ref{Trivial tangent sheaf} hold.
Condition \eref{Elliptic} holds by definition on good mutations.
According to Proposition \ref{conditions fpb},
the Frobenius pullback $X'=J^\fpb$   satisfies 
conditions \eref{Tjurina} and  \eref{Zariski}.

The fibers on $J$  and hence on $X'$ are reduced by assumption, and this obviously
also transfers to  the fibers on $X$, so condition \eref{Reduced}  holds.
It remains to check that a mutation $X'\leftarrow S\rightarrow X$ 
produces on $X$ only Zariski singularities, and that the total Tjurina does not change.
This is obvious if we are in the case \mref{Semistable}, when $S_a$ is semistable.
It remains to verify the cases \mref{II} and \mref{III}.
Then $S_a$ is a fiber of Kodaira type $\I^*_n$ for some $n\geq 0$.
The following argument shows that the ensuing singularities on $X'$ and $X$ are formally isomorphic.
Set $R=\O_{\PP^1,a}$ and $F=\Frac(R)$. For each terminal component $\Theta,\Theta'\subset S_a$,
we may choose formal sections $A,A'\subset S\otimes \hat{R}$ passing through $\Theta$ and $\Theta'$,
respectively.
Since the generic fiber $U=S\otimes\hat{F}$ is an elliptic curve, there is a translation automorphism
$\tau:U\ra U'$ sending the $\hat{F}$-rational point $A\cap U$ to $A'\cap U$. By the functoriality of the relative minimal
model, this extends to an automorphism of $S\otimes\hat{R}$ sending $\Theta$ to $\Theta'$.
In turn, the singularity arising from contractions of all components but $\Theta$ is
formally isomorphic to the singularity given by the contraction of all components but $\Theta'$.
This settles the case \mref{II}, and the case \mref{III} follows analogously.
\qed

\medskip
It is very easy to specify   mutations. The challenge, however,  is to introduce another elliptic
fibration as well. Rather than giving an elliptic curve, we shall construct its  degenerate fibers.
This boils down to finding   sections   for $J'\ra\PP^1$ with suitable intersection behavior.
Using the following, we shall infer that the resulting fibration  is elliptic rather then quasielliptic:

\begin{proposition}
\mylabel{good mutations}
Let $X'\leftarrow S\rightarrow X$ be a mutation with respect to the vertical negative-definite $E\subset S$.
Suppose that there is curve of canonical type $F\subset S$ with the following properties:
\begin{enumerate}
\item The type of the curve $F$  is either $\I_n$ or $\I_{2n+1}^*$ or $\IV$ or $\IV^*$.
\item It has degree $>0$ with respect to the elliptic fibration $S\ra\PP^1$.
\item Each irreducible component of $E$ is either contained in or disjoint from $F$.
\end{enumerate} 
Then $X'\leftarrow S\rightarrow X$ is a good mutation.
\end{proposition}

\proof
The curve of canonical type $F\subset S$ induces a genus-one fibration $g:S\ra\PP^1$ on the K3 surface $S$.
Condition (i) ensures that the new fibration is elliptic, in light of  \cite{Rudakov; Safarevic 1978}, 
the Proposition on page 150.
Condition (ii) implies that this fibration is different from our original elliptic fibration $S\ra\PP^1$.
By the last condition, the new elliptic fibration on $S$ induces a new elliptic fibration on $X$.
One of its closed fibers is an elliptic curve $E\subset X$ that is horizontal with respect
to the old elliptic fibration $X\ra\PP^1$. Thus the mutation $X'\leftarrow S\rightarrow X$ is good.
\qed

\medskip
We now see that in most cases for $J\ra\PP^1$, one finds a good mutation, which yields
a normal K3-like covering. The following are the configurations of Kodaira symbols
for $J$ where the approach does not work at  the moment:
\begin{equation}
\mylabel{bad kodaira symbols}
\III,\, \II+5,\, \III+6,\, \IV+\IV+\IV,\,\IV+\II,\, \IV+\III.
\end{equation}
The main result of this section is:

\begin{theorem}
\mylabel{good mutations exist}
Suppose that the rational elliptic surface $J\ra\PP^1$ has a configuration of Kodaira symbols
different from the ones in \eqref{bad kodaira symbols}.
Then there exists a good mutation $X'\leftarrow S\rightarrow X$, and the resulting $X$ is a K3-like
covering.
\end{theorem}

\proof
The second statement follows from Theorem \ref{mutations K3-like coverings}.
The strategy for the first statement  is to use Lang's classification for $J\ra\PP^1$ together with 
the Oguiso--Shioda list for the Mordell--Weil lattice 
$\MW(J/\PP^1)$, in order  to guess some vertical negative-definite curve $E\subset S$ that
yields a  mutation, as well as 
to produce suitable sections $P,P',Q,\ldots\subset J$
that allow to construct a curve of canonical type $F\subset S$ 
for which Proposition \ref{good mutations} applies. Throughout, $O\subset J$ denotes
the zero-section. Recall from  Proposition \ref{disjoint section} that in any case there is a section 
 disjoint from the zero-section $O$. 
To simplify  notation, we now write  $a,b\in\PP^1$ for closed points, and 
$S_a,S_b,X'_a,X'_b$ and $J_a,J_b$  for the   fibers of the fibrations $S\ra\PP^1$, $X'\ra\PP^1$ and $J\ra\PP^1$.
We have to examine five  cases:

\newcounter{CasesForJ}
\setcounter{CasesForJ}{0}

\medskip
{\bf Case (\stepcounter{CasesForJ}\roman{CasesForJ}):}
\emph{Suppose there are  at least two semistable fibers $J_a$ and $J_b$.}
Then  define $E_a\subset S_a$ and $E_b\subset S_b$
as the strict transform of the   fibers $X'_a$ and $X'_b$, respectively.
Let $E\subset S$ be the curve obtained from $E'$ by replacing $E'_a+E'_b$ by $E_a+ E_b$.
This yields a mutation  $X'\leftarrow S\ra X$.
Now choose  a section $P\subset J$ disjoint from the zero-section $O\subset J$.
Together with suitable chains of irreducible components inside $J_a$ and $J_b$, they form
a cycle of projective lines. Their   preimage $F\subset S$ becomes a cycle
of $(-2)$-curves, that is, a curve of canonical type $\I_n$ for some $n\geq 4$.
By construction, each irreducible component of $E$ is either disjoint from or contained in $F$. 
Thus Proposition \ref{good mutations} applies.

\medskip
{\bf Case (\stepcounter{CasesForJ}\roman{CasesForJ}):}
\emph{Suppose there is exactly one semistable fiber $J_b$.}
According to \cite{Lang 2000}, Section 1 and 2, there are four possibilities for the configuration of singular fibers:
$$
\II+8,\quad \III+ 8, \quad \II+5,\quad \III+6.
$$
The respective Lang types of the unstable fiber $J_a$ are  $1A,2A,1C,2B$, and the latter two configurations
are excluded from consideration.
Write the semistable fiber $J_b=\sum\Theta_i$ with cyclic indices in the natural way, that is,
$(\Theta_0\cdot O)=1$ and $(\Theta_i\cdot\Theta_{i+1})=1$.

Suppose that the configuration is $\II+8$.
The fiber $J_a$ has numerical invariants
$v=4$, $m=1$, $\delta=2$.   On $S$  the invariants become
$v'=8$ and $m'=5$, whence $S_a$ has Kodaira type $\I^*_0$, and  
the corresponding singularity on $X'$ is of type $D_4$.
Fix some terminal irreducible component in $S_a$ different from the strict transform of $X'_a$,
and define $E_a$ as the union of the remaining components $\Theta\subset S_a$.
Let $E\subset S$ be the curve obtained from $E'$ by replacing $E'_a$ by $E_a$.
This  yields a mutation  $X'\leftarrow S\rightarrow X$.
The $E_a$, the strict transforms of the zero-section $O$ and $\Theta_0$, 
and the   two components in the
semistable fiber $S_b$ adjacent to the strict tranform of $\Theta_0$ form curve $F\subset S$ of canonical type $I^*_3$,
and each irreducible component of $E$ is either disjoint from or contained in $F$.
The case that the configuration is $\III+8$ is analogous and left to the reader.


\medskip
\emph{From now on, we  assume that all singular fibers $J_a$ are unstable.}

\medskip
{\bf Case (\stepcounter{CasesForJ}\roman{CasesForJ}):} 
\emph{Suppose there is a fiber $J_a$ of Kodaira type $\II$, whose Lang type is   9A or  9C.}
This means $v=4r$, $m=1$ and $\delta=4r-2$, with respective values $r=1, 3$.
In turn, the invariants on  $S$ become $v'=8r$ and $m'=4r+1$.
The only possibility is that the fiber $S_a$ has type $\I^*_{4r-4}$,
and the singularity on  $X'=J^{(2/\PP^1)}$ is of type $D_{4r}$.
Fix the terminal component of $S_a$ nearest to the strict
transform of $X'_a$, and let $E_a\subset S_a$ be the union of the remaining irreducible
components.
Let $E\subset S$ be the curve obtained from $E'\subset S$ by replacing
$E'_a$ by $E_a$. This yields a mutation $X'\leftarrow S\rightarrow X$.
Let $P\subset J$ be a section disjoint from the zero-section $O\subset J$.
Together with the $4r$ components of  $E_a$, they induce a  canonical curve $F\subset S$
of type  $\I^*_{4r-3}$. By construction, each irreducible component from $E$ is either disjoint
from or contained in $F$.

\medskip
{\bf Case (\stepcounter{CasesForJ}\roman{CasesForJ}):} 
\emph{Suppose there is a  fiber $J_a$ of Kodaira type $\III$ with Lang type 10A.}
The fiber $J_a$ then has invariants $v=4$, $m=2$ and $\delta=1$. In turn, 
the surface $S$ gets a fiber with $v'=8$ and $m'=6$,
which must be of type $\I_1^*$, and the singularity on $X'$ is of type $D_4$.
Let $E_a$ be the union of the two terminal components in $X'_a$ disjoint from $E_a'$.
Define $E\subset S$ as the curve obtained from $E'\subset S$ by replacing
$E'_a$ by $E_a$.  Its  contraction yields a mutation $X'\leftarrow S\rightarrow X'$.

According to Lang's classification, we have the following seven possibilities for the configuration
of Kodaira types on $J$:
$$
\III+\II+\II,\,
\III+\III+\II,\,
\III+\III+\III,\,
\III+\II+\IV,\,
\III+\IV+\IV,\,
\III+\II,\,
\III+\III.
$$
The corresponding trivial lattices are of the form $A_1^{\oplus n}$ with $n=1,2,3$
and $A_1\oplus A_2^{\oplus m}$ with $m=1,2$.
In the Oguiso--Shioda classification, this are No.\ 2,4,7,6,20.
One sees that in all cases,
the narrow Mordell--Weil lattice contains a vector $\alpha$ with $\alpha^2=2$.
In turn, there is a section $P\subset J$ intersecting the same fiber components as the
zero-section $O\subset J$, with
$$
2=\langle P,P\rangle = 2+ 2(P\cdot O).
$$
It follows that $O,P$ are disjoint. Consequently $O,P$ and the four components of $E_a$ induce a curve
$F\subset S$ of canonical type $\I_1^*$.
By construction, each irreducible component of $E$ is either contained in or disjoint from $F$.

\medskip
{\bf Case (\stepcounter{CasesForJ}\roman{CasesForJ}):} 
\emph{Suppose there is a fiber $J_a$ of Kodaira type $\IV$.}
This fiber has Lang type 11, with invariants $v=4$, $m=3$, $\delta=0$.
Whence the fiber  $S_a$ has  invariants $v'=8$, $m'=7$. 
Their possible Kodaira types are $\IV^*$ or $\I^*_2$.
The respective singularities on $X'$ are $D_4$ or $A_4$.
The latter is not Zariski by Theorem \ref{zariski rdp},
hence the singularity is $D_4$ and the Kodaira symbol is $\IV^*$.
Those containing a fiber of Lang type $9A, 9C, 10A$  were treated above.
According to Lang's Classification, the only remaining possibilities
are $\IV+\II\quadand \IV+\III$, which we have excluded from consideration.

\medskip
This completes the proof:
Cases (i)   applies if all fibers are semistable, according to \cite{Lang 2000}, Section 1.
It also covers all but four cases with $j(J/\PP^1)\neq 0$, by  \cite{Lang 2000}, Section 2.
The  extra cases  are treated in   (ii). The remaining cases have global
$j$-invariant $j(J/\PP^1)=0$, which were analysed in \cite{Lang 2000}, Section 3. Then all fibers are unstable,
and their is a fiber $J_a$ of Lang type 9A, 9C, 10A or 11;
these are treated in (iii)--(v). The exceptions that contain a fiber of Lang type $10C$ or $11$
are excluded from consideration.
\qed

\medskip
Another way to construct Enriques surface $Y$ from an rational elliptic surface $J\ra\PP^1$
is via Ogg--Shafarevich theory, as exposed in \cite{Cossec; Dolgachev 1989}, Chapter V, 
Section 4.


\begin{thebibliography}{ccccc}

\bibitem{Artin 1966}
M.\ Artin:
On isolated rational singularities of surfaces.
Am.\ J.\ Math.\ 88 (1966), 129--136.

\bibitem{Artin 1977}
M.\ Artin:
Coverings of the rational double points in characteristic $p$.
In: W.\ Baily, T.\ Shioda (eds.), 
Complex analysis and algebraic geometry, pp.\ 11--22.
Iwanami Shoten, Tokyo, 1977. 

\bibitem{Bayer; Eisenbud 1995}
D.\ Bayer, D. Eisenbud:
Ribbons and their canonical embeddings.
Trans.\ Am.\ Math.\ Soc.\ 347 (1995), 719--756. 

\bibitem{Blass 1982}
P.\ Blass:
Unirationality of Enriques surfaces in characteristic two.
Compositio Math.\ 45 (1982),   393--398. 

\bibitem{Blass; Lang 1987}
P.\ Blass, J.\ Lang:
Zariski surfaces and differential equations in characteristic $p>0$.
Marcel Dekker, New York, 1987.

\bibitem{Bombieri; Mumford 1976}
E.\ Bombieri, D.\ Mumford:
Enriques' classification of surfaces in char.\ $p$,  III.
Invent.\ Math.\ 35  (1976), 197--232.

\bibitem{Bombieri; Mumford 1977}
E.\ Bombieri, D.\ Mumford:
Enriques' classification of surfaces in char.\ $p$,  II.
In: W.\ Baily, T.\ Shioda (eds.),
Complex analysis and algebraic geometry, pp.\ 23--42.
Cambridge University Press, London, 1977.

\bibitem{Bosch; Luetkebohmert; Raynaud 1990}
S.\ Bosch, W.\ L\"utkebohmert, M.\ Raynaud:
N\'eron models.
Springer, Berlin, 1990.

\bibitem{AC 8-9}
N.\ Bourbaki:
Alg\`ebre commutative. Chapitre 8--9.
Masson, Paris, 1983.

\bibitem{Cossec 1985}
F.\ Cossec:
On the Picard group of Enriques surfaces.
Math.\ Ann.\ 271 (1985),  577--600. 

\bibitem{Cossec; Dolgachev 1989}
F.\ Cossec, I.\ Dolgachev:
Enriques surfaces I.
Birkh\"auser, Boston, MA, 1989.

\bibitem{Demazure; Gabriel 1970}
M.\ Demazure, P.\ Gabriel:
Groupes alg\'ebriques.
Masson, Paris, 1970.

\bibitem{SGA 3b}
M.\ Demazure, A.\ Grothendieck (eds.):
Sch\'emas en groupes II (SGA 3 Tome 2).
Springer, Berlin, 1970.


\bibitem{Ekedahl 1988}
T.\ Ekedahl:
Canonical models of surfaces of general type in positive characteristic. 
Inst.\ Hautes \'Etudes Sci.\ Publ.\ Math.\ 67 (1988), 97--144.

\bibitem{Ekedahl; Shepherd-Barron 2004}
T.\ Ekedahl, N. Shepherd-Barron:
On exceptional Enriques surfaces.
Preprint, math.AG/0405510.

\bibitem{Ekedahl; Hyland; Shepherd-Barron 2012} 
T.\ Ekedahl, J.\ Hyland, N. Shepherd-Barron, 
Moduli and periods of simply connected Enriques surfaces.
Preprint, math.AG/1210.0342.

\bibitem{EGA II}
A.\ Grothendieck:
\'El\'ements de g\'eom\'etrie alg\'ebrique II:
\'Etude globale \'el\'ementaire de quelques classes de morphismes.
Publ.\ Math., Inst.\ Hautes \'Etud.\ Sci.\ 8 (1961).

\bibitem{EGA IVc}
A.\ Grothendieck:
\'El\'ements de g\'eom\'etrie alg\'ebrique IV: \'Etude locale des
sch\'emas et des morphismes de sch\'emas.
Publ.\ Math., Inst.\ Hautes \'Etud.\ Sci.\  28 (1966).

\bibitem{EGA IVd}
A.\ Grothendieck:
\'El\'ements de g\'eom\'etrie alg\'ebrique IV: \'Etude locale des
sch\'emas et des morphismes de sch\'emas.
Publ.\ Math., Inst.\ Hautes \'Etud.\ Sci.\   32 (1967).

\bibitem{SGA 1}
A.\ Grothendieck:
Rev\^etements \'etales et groupe fondamental (SGA 1).
Soci\'et\'e Math\'ematique de France, Paris, 2003. 

\bibitem{SGA 2}
A.\ Grothendieck:
Cohomologie locale des faisceaux coh\'erents et th\'eor\`emes de Lefschetz locaux et globaux (SGA 2).
North-Holland Publishing Company, Amsterdam, 1968.

\bibitem{Hartshorne 1994}
R.\ Hartshorne:
Generalised divisors on Gorenstein schemes.
K-Theory 8 (1994), 287--339.

\bibitem{Ito 1994}
H.\ Ito:
The Mordell--Weil groups of unirational quasi-elliptic surfaces in characteristic $2$.  
Tohoku Math.\ J.\   46  (1994),   221--251.

\bibitem{Katsura 1982}
T.\ Katsura:
A note on Enriques' surfaces in characteristic 2.
Compositio Math.\ 47 (1982),  207--216. 

\bibitem{Katsura; Kondo 2015}
T.\ Katsura, S.\ Kondo:
1-dimensional family of Enriques surfaces in characteristic 2 covered by the supersingular K3 surface with Artin invariant 1.
Pure Appl.\ Math.\ Q.\ 11 (2015), 683--709. 

\bibitem{Kiehl; Kunz 1965}
R.\ Kiehl, E.\ Kunz:
Vollst\"andige Durchschnitte und $p$-Basen. 
Arch.\ Math.\ 16 (1965), 348--362. 

\bibitem{Lang 1983}
W.\ Lang:
On Enriques surfaces in characteristic p. I. 
Math.\ Ann.\ 265 (1983),  45--65. 

\bibitem{Lang 1988}
W.\ Lang:
On Enriques surfaces in characteristic p. II.
Math.\ Ann.\ 281 (1988),  671--685. 

\bibitem{Lang 1994}
W.\ Lang:
Extremal rational elliptic surfaces in characteristic $p$. II. 
Surfaces with three or fewer singular fibres.  
Ark.\ Mat.\  32  (1994), 423--448.

\bibitem{Lang 2000}
W.\ Lang:
Configurations of singular fibres on rational elliptic surfaces in characteristic two. 
Comm.\ Algebra  28  (2000),  5813--5836.

\bibitem{Laufer 1977}
H.\ Laufer:
On minimally elliptic singularities,  
Am.\ J.\ Math.\ 99  (1977), 1257--1295.

\bibitem{Liedtke 2015}
C.\ Liedtke:
Arithmetic moduli and lifting of Enriques surfaces. 
J.\ Reine Angew.\ Math.\ 706 (2015), 35--65. 

\bibitem{Liu; Lorenzini; Raynaud 2005}
Q.\ Liu, D.\ Lorenzini, M.\ Raynaud:
N\'eron models, Lie algebras, and reduction of curves of genus one.
Invent.\ Math.\ 157 (2004),  455--518.

\bibitem{Magma}
W.\ Bosma, J.\ Cannon, C.\ Playoust:
The Magma algebra system. I. The user language.
J.\ Symbolic Comput., 24 (1997), 235--265. 

\bibitem{Miller 2003}
C.\ Miller:
The Frobenius endomorphism and homological dimensions.  
In: L.\ Avramov, M.\ Chardin, M.\ Morales and C.\ Polini (eds.), Commutative algebra, pp.\ 207--234.
Contemp.\ Math.\ 331. Amer.\ Math.\ Soc., Providence,  2003.

\bibitem{Miranda 1990}
R.\ Miranda:
Persson's list of singular fibers for a rational elliptic surface. 
Math.\ Z.\ 205 (1990), 191--211.

\bibitem{Mumford 1966}
D.\ Mumford:
Lectures on curves on an algebraic surface.
Princeton University Press, Princeton, 1966.

\bibitem{Mumford 1970}
D.\ Mumford:
Abelian varieties.
Tata Institute of Fundamental Research Studies in Mathematics 5.
Oxford University Press,  London, 1970.

\bibitem{Ogg 1967}
A.\ Ogg:
Elliptic curves and wild ramification.
Amer.\ J.\ Math.\ 89 (1967), 1--21. 

\bibitem{Oguiso; Shioda 1991}
K.\ Oguiso, T.\ Shioda:
The Mordell--Weil lattice of a rational elliptic surface.
Comment.\ Math.\ Univ.\ St.\ Paul.\ 40 (1991), 83--99.

\bibitem{Okonek; Schneider; Spindler 1980}
C.\ Okonek, M.\ Schneider, H.\  Spindler:
Vector bundles on complex projective spaces. 
Birkh\"auser, Boston, Mass., 1980. 

\bibitem{Persson 1990}
U.\ Persson:
Configurations of Kodaira fibers on rational elliptic surfaces.
Math.\ Z.\ 205 (1990),  1--47. 

\bibitem{Raynaud 1970}
M.\ Raynaud:
Sp\'ecialisation du foncteur de Picard.
Publ.\ Math., Inst.\ Hautes \'Etud.\ Sci.\ 38 (1970), 27--76.

\bibitem{Restuccia; Schneider 2003}
G.\ Restuccia, H.-J.\ Schneider:
On actions of infinitesimal group schemes.
J.\ Algebra 261 (2003), 229--244.

\bibitem{Rudakov; Safarevic 1978}
A.\ Rudakov, I.\ Safarevic:
Supersingular $K3$ surfaces over fields of characteristic $2$.
Math.\ USSR, Izv.\ 13 (1979), 147--165.

\bibitem{Saito 1988}
T.\ Saito:
Conductor, discriminant, and the Noether formula of arithmetic surfaces.
Duke Math.\ J.\ 57 (1988), 151--173. 

\bibitem{Schroeer 2007}
S. Schr\"oer:
Kummer surfaces for the selfproduct of the cuspidal rational curve.
J.\ Algebraic Geom.\ 16 (2007), 305--346.

\bibitem{Schroeer 2008}
S.\ Schr\"oer:
Singularities appearing on generic fibers of morphisms between smooth schemes.
Michigan Math.\ J.\ 56 (2008), 55--76.

\bibitem{Shioda 1990}
T.\ Shioda:
On the Mordell--Weil lattices.
Comment.\ Math.\ Univ.\ St.\ Paul.\ 39 (1990), 211--240. 

\bibitem{Strade; Farnsteiner  1988}
H.\ Strade, R.\  Farnsteiner:
Modular Lie algebras and their representations. 
Marcel Dekker, New York, 1988.

\bibitem{Tate 1972}
J.\ Tate:
Algorithm for determining the type of a singular fiber in an elliptic pencil.  
In: B.\ Birch, W.\ Kuyk (eds.),
Modular functions of one variable IV,  pp. 33--52. 
Springer, Berlin, 1975.

\bibitem{Tate; Oort 1970}
J.\ Tate,  F.\ Oort:
Group schemes of prime order.  
Ann.\ Sci.\ \'Ec.\ Norm.\ Sup\'er.\   3  (1970), 1--21.

\bibitem{Wagreich 1970}
P.\ Wagreich:
Elliptic singularities of surfaces.  
Amer.\ J.\ Math.\  92  (1970), 419--454.

\bibitem{Wang 2013}
X.\ Wang:
Connected Hopf algebras of dimension $p^2$. 
J.\ Algebra 391 (2013), 93--113.

\end{thebibliography}
\end{document}